\input amstex
\documentstyle{amsppt}
\input xy
\xyoption{all}
\topmatter
\title
Kac-Moody Algebras, the Monstrous Moonshine, Jacobi Forms and
Infinite Products
\endtitle
\leftheadtext{Jae-Hyun Yang}
\rightheadtext{Jacobi Forms and Infinite Products}
\author  Jae-Hyun Yang
\endauthor
\magnification=\magstep 1
\pagewidth{12.5cm} \pageheight{18.20cm}
\baselineskip =7mm
\thanks{This work was in part supported by TGRC-KOSEF.
\endgraf
Mathematics Subject Classification (1991)\,: Primary 11F30, 11F55,
17B65, 17B67, 20C34, 20D08.}
\endthanks
\endtopmatter
\document
\NoBlackBoxes
\head  {\bf Table of Contents}\endhead

\def\ex{\par\smallpagebreak}
\vskip 0.5cm
\indent \ 1. Introduction\par\smallpagebreak
\indent \ \ \ \ \  Notations\par \smallpagebreak
\indent \ 2. Kac-Moody Lie Algebras\par\smallpagebreak
\indent \ \ \ \ \ Appendix\,:\ Generalized Kac-Moody Algebras\ex
\indent \ 3. The Moonshine Conjectures and the Monster Lie Algebra\ex
\indent \ \ \ \ \ Appendix\,:\ The No-Ghost Theorem\ex
\indent \ 4. Jacobi Forms\par
\indent \ 5. Infinite Products and Modular Forms\ex
\indent \ 6. Final Remarks\ex
\indent \ \ \ \ 6.1. The Fake Monster Lie Algebras\ex
\indent \ \ \ \ 6.2. Generalized Kac-Moody Algebras of the Arithmetic
Type\ex
\indent \ \ \ \ 6.3. Open Problems\ex
\indent \ Appendix A\,: Classical Modular Forms\ex
\indent \ Appendix B\,: Kohnen Plus Space and Maass Space\ex
\indent \ Appendix C\,: The Orthogonal Group $O_{s+2,2}({\Bbb R})$\ex
\indent \ Appendix D\,: The Leech Lattice $\Lambda$\ex
\indent \ References

\vskip 0.5cm
\head {\bf 1.\ Introduction}\endhead
\indent
Recently R.\,E. Borcherds obtained some quite interesting results in
[Bo6-7]. First he solved the Moonshine Conjectures made by Conway and
Norton\,(\,[C-N]\,). Secondly he constructed automorphic forms on the
orthogonal group $O_{s+2,2}({\Bbb R})$ which are modular products and
then wrote some of the well-known meromorphic modular forms as infinite
products. Modular products roughly mean infinite products whose exponents
are the coefficients of certain nearly holomorphic modular forms.
The theory of Jacobi forms plays an
important role in his second work in [Bo7]. More than 10 years ago
Feingold and Frenkel\,(\,[F-F]\,)
realized the connection between the theory of a
special hyperbolic Kac-Moody Lie algebra of the type
$HA_1^{(1)}$ and that of Jacobi
forms of degree one and then generalized the results of H. Maass to
higher levels. So far the relationship between the theory of Jacobi forms
of higher degree and that of other hyperbolic Kac-Moody algebras has not
been developed yet. The work of Borcherds in [Bo7] gives a light on the
possibility for the relationship between them. This fact urged me to
write a somewhat supplementary or expository note on Borcherds' recent works
which is useful for my research on Jacobi forms although I am not an
expert in the theory of Kac-Moody Lie algebras and lattices.
I hope that this note will be useful for the readers who are interested
in these interesting subjects. I learned
a lot about these subjects while I had been preparing this article.
I had given a lecture on some of these materials at the 4th Symposium
of the Pyungsan Institute for Mathematical Sciences held at the Wonkwang
University in September, 1995.
\vskip 0.052cm As mentioned above, the purpose of this paper is to
give a survey of Borcherds' recent results in [Bo6-7] to the core.
This article is organized as follows. In section 2, we collect
some of the well-known properties of Kac-Moody Lie algebras, e.g.,
the Weyl-Kac character formula, the root multiplicity and so on.
In the appendix, we discuss the generalized Kac-Moody Lie algebras
introduced by Borcherds roughly. In section 3, we give a sketchy
survey on the Moonshine Conjectures solved by
Borcherds\,(\,[Bo6]\,). We discuss the monster Lie algebra and the
no-ghost theorem. This section is completely based on the article
[Bo6]. In section 4, we review some of the theory of Jacobi forms
and discuss singular Jacobi forms briefly. We present Borcherds'
construction of nearly holomorphic Jacobi forms by making use of
the concept of affine vector systems. In section 5, we give a
brief history of infinite products and present the work of
Borcherds that expressed modular forms in the Kohnen ``plus" space
of weight $1/2$ as infinite products. For instance, we write the
well-known modular forms like the discriminant function
$\Delta(\tau)$, the modular invariant $j(\tau)$ and the Eisenstein
series $E_k(\tau)\,(\,k\geq 4,\,k:even\,)$ as infinite products
explicitly. In the final section, we discuss the fake monster Lie
algebras and Kac-Moody Lie algebras of the arithmetic hyperbolic
type defined by V.\,V. Nikulin\,(\,[N5]\,). As an example, we
explain the generalization of Maass correspondence to higher
levels which was done by A.\,J. Feingold and I.\,B.
Frenkel\,(\,[F-F]\,). Finally we also give some open problems
which have to be investigated. In the appendix A, we collect some
of the well-known properties of classical modular forms. In the
appendix B, we briefly discuss the Kohnen ``plus" space and the
Maass ``Spezialschar" which are essential for the understanding of
the works in [Bo7] and [F-F]. In the appendix C, we discuss the
geometrical aspect of the orthogonal group $O_{s+2,2} ({\Bbb R})$
briefly. In the final appendix, we collect some of the well-known
properties of the Leech lattice $\Lambda.$ Also we briefly discuss
the Jacobi theta functions.
\par\smallpagebreak\indent
Finally I would like to give my deep thanks to TGRC-KOSEF for its
financial support on this work. I also would like to give my hearty
thanks to my Korean colleagues for their interest in this work.
\par\smallpagebreak\noindent
\def\BZ{\Bbb Z}
\def\BQ{\Bbb Q}
\def\BR{\Bbb R}
\def\BC{\Bbb C}
\par
\noindent
{\bf
Notations:} \ \ We denote by $\BZ,\,\BQ,\,\BR$ and $\BC$ the ring of integers,
the field of rational numbers, the field of
real numbers, and the field of complex numbers respectively.
$\BZ^+$ and $\BZ_+$ denote the set of all positive integers and the set of
nonnegative integers respectively.
For a positive integer $n,\ \Gamma_n:=Sp(n,\BZ)$ denotes the Siegel modular
group of degree $n$.
The symbol
``:='' means that the expression on the right is the definition of that on the
left.
$F^{(k,l)}$ denotes
the set of all $k\times l$ matrices with entries in a commutative ring $F$.
For a square matrix $A\in F^{(k,k)}$ of degree $k$, $\sigma(A)$ denotes the
trace of $A$.
For $A\in F^{(k,l)}$ and $B\in F^{(k,k)},$ we set $B[A]=\,^t\!ABA.$ For
any $M\in F^{(k,l)},\ ^t\!M$ denotes the transpose matrix of $M$. $E_n$ denotes
the identity matrix of degree $n$.
We denote by $\Lambda$ the Leech lattice. $\Pi_{1,1}$ denotes the unique
unimodular even integral Lorentzian lattice of rank $2.$ We denote by $G$
the MONSTER group. For $g\in G,\ T_g(q)$ denotes the Thompson series of
$g$. $M$ and $V^{\sharp}$ denote the monster Lie algebra and the moonshine
module respectively. $\eta(\tau)$ denotes the Dedekind eta function.
$\tau(n)$ denotes the Ramanujan function. Usually $\rho$ denotes the Weyl
vector. We denote by $[\Gamma_n,k]\,(\,\text{resp.}\,[\Gamma,k]_0\,)$
the complex vector space of all Siegel modular forms\,(\,resp.\, cusp
forms\,) on $H_n$ of weight $k$ with respect to $\Gamma_n$. We denote
by $[\Gamma_2,k]^M$ the Maass space or the Maass Spezialschar.
\par\medpagebreak
\vskip 0.3cm

\define\g{\frak g}
\define\h{\frak h}
\define\a{\alpha}
\redefine\v{\vee}
\redefine\l{\lambda}
\redefine\D{\Delta}
\define\m{\underline m}

\head {\bf 2. Kac-Moody Lie Algebras}\endhead
\par
In this section, we review the basic definitions and properties of Kac-Moody Lie
algebras.\par
An $n\times n$ matrix $A=(a_{ij})$ is called a {\it generalized Cartan
matrix} if it satisfies the following conditions:
(i) $a_{ii}=2$ for $i=1,\cdots,n;$ (ii) $a_{ij}$ are nonpositive integers for
$i\neq j;$ (iii) $a_{ij}=0$ implies $a_{ji}=0.$ An indecomposable
generalized Cartan matrix is said to be {\it of finite type} if all its
principal minors are positive, {\it of affine type} if all its proper principal
minors are positive and det$\,A=0,$ and is said to be {\it of indefinite type} if
$A$ is of neither finite type nor affine type. $A$ is said to be {\it of
hyperbolic type} if it is of indefinite type and all of its proper principal
submatrices are of finite type or affine type, and to be {\it of almost
hyperbolic type} if it is of indefinite type and at least one of its proper
principal submatrices is of finite or affine type.\par
A matrix $A$ is called {\it symmetrizable} if there exists an invertible
diagonal matrix $D=\text{diag}(q_1,\cdots,q_n)$ with $q_i>0,\;q_i\in\Bbb Q$ such that
$DA$ is symmetric. If $A$ is an $n\times n$ matrix of rank $l,$ then {\it a
realization} of $A$ is a triple $(\h,\Pi,\Pi^\v),$ where $\h$ is a complex
vector space of dimension $2n-l,$ $\Pi=\{\a_1,\cdots,\a_n\}$ and
$\Pi^\v=\{\a_1^\v,\cdots,\a_n^\v\}$ are linearly independent subsets of $\h^*$ and $\h$
respectively such that $\a_j(\a_i^\v)=a_{ij}$ for $1\leq i,j\leq n.$
\definition {Definition 2.1} \ The {\it Kac-Moody Lie algebra} $\g(A)$ associated
with the generalized Cartan matrix $A$ is the Lie algebra generated by the
elements $e_i,\;f_i\;(i=1,2,\cdots,n)$ and $\h$ with the defining relations
$$\align
&[h,h']=0\qquad \text{for all}\;h,h'\in\h,\\
&[e_i,f_j]=\delta_{ij}\a_i^\v\qquad\text{for}\; 1\leq i,j\leq n,\\
&[h,e_i]=\a_i(h)e_i,\;\;[h,f_i]=-\a_i(h)f_i\;\quad\text{for}\;
i=1,2,\cdots,n,\\
&(ad\,e_i)^{1-a_{ij}}(e_j)=\;\;(ad\,f_i)^{1-a_{ij}}(f_j)=0\;\quad\text{for}\;
i\neq j.
\endalign$$
The elements of $\Pi$ (resp. $\Pi^\v$) are called {\it simple roots}
(resp. {\it simple coroots}) of $\g (A).$
\enddefinition
Let $Q:=\sum_{i=1}^n\Bbb Z\a_i,\,Q_+:=\sum_{i=1}^n\Bbb Z_+\a_i$ and
$Q_-:=-Q_+.$  $Q$ is called the {\it root lattice}.
For $\a:=\sum_{i=1}^n k_i\a_i\in Q$ the number
$\text{ht}(\a):=\sum_{i=1}^n k_i$ is called the {\it height} of $\a.$ We define a
{\it partial ordering} $\geq $ on $\h^*$ by $\l\geq\mu$ if $\l-\mu\in Q_+.$
For each $\a\in\h^*,$ we put
$$\g_\a:=\{X\in\g(A)\mid [h,X]=\a(h)X\;\text{for all}\;h\in\h\}.$$
If $\g_\a\neq 0,\;\a$ is called a {\it root} and $\g_\a$ is called the
{\it root space} attached to $\a.$ The number $\text{mult}\,\a:=\dim\,\g_\a$ is
called the {\it multiplicity} of $\a.$ The Kac-Moody Lie algebra $\g(A)$ has
the following root space decomposition with respect to $\h$:
$$\g(A)=\sum_{\a\in Q}\g_\a\qquad\text{(direct sum).}\tag2.1$$
A root $\a$ with $\a>0$ (resp. $\a<0$) is called {\it positive}
(resp. {\it negative}). All roots are either positive or negative. We denote
by $\D,\D^+$ and $\D^-$ the set of all roots, positive roots and negative
roots respectively.
\definition {Definition 2.2}
\ Let $\g(A)$ be a symmetrizable Kac-Moody Lie algebra associated with a
symmetrizable generalized Cartan matrix $A=(a_{ij}).$ A $\g(A)$-module $V$ is
$\h-diagonalizable$ if $V=\oplus_{\mu\in\h^*}V_\mu,$ where $V_\mu$ is the
weight space of weight $\mu$ given by
$$V_\mu:=\{v\in V\mid hv=\mu(h)v\;\text{for all}\;h\in\h\}\neq 0.$$
The number $\text{mult}_V\mu:=\dim V_\mu$ is called the {\it multiplicity} of weight
$\mu.$ When all the weight spaces are finite-dimensional, we define the
{\it character} of $V$ by
$$\text{ch}\,V:=\sum_{\mu\in\h^*}(\dim V_\mu)e^\mu=\sum_{\mu\in\h^*}
(\text{mult}_V\mu)e^\mu.
\tag2.2$$
An $\h$-diagonalizable $\g(A)$-module $V$ is said to be {\it integrable} if
all the Chevalley generators $e_i,f_i(i=1,2,\cdots,n)$ are locally nilpotent
on $V.$ A $\g(A)$-module $V$ is called
a {\it highest weight module} with highest weight $\Lambda\in\h^*$ if there exists
a nonzero vector $v\in V$ such that
(i) $e_i v=0$ for all $i=1,2,\cdots,n;$ (ii) $hv=\Lambda(h)v$ for all $h\in\h;$
and (iii) $U(\g(A))\; v=V.$
A vector $v$ is called a {\it highest weight vector.} Here $U(\g(A))$ denotes
the universal enveloping algebra of $\g(A).$
\enddefinition
Let $\frak n_+$ (resp. $\frak n_-$) be the subalgebra of $\g(A)$ generated by
$e_1,\cdots,e_n$ (resp. $f_1,\cdots,f_n$). Then we have the {\it triangular
decomposition}
$$\g(A)=\frak n_-\oplus\h\oplus\frak n_+\qquad\text{(direct sum of vector spaces)}.$$
Let $\frak b_+:=\h+\frak n_+$ be the {\it Borel subalgebra} of $\g(A).$
For a fixed
$\Lambda\in\h^*,$ we let $\Bbb C(\Lambda)$ be the one-dimensional $\frak b_+$-module
with the $\frak b_+$-action defined by
$$\frak n_+\cdot 1=0\;\;\text{and}\;\;h\cdot 1=\Lambda(h)1\;\;
\text{for all}\;\;h\in\h.$$
The induced module
\par\smallpagebreak
$$M(\Lambda):=U(\g(A))\otimes_{U(\frak b_+)}\Bbb C(\Lambda)\tag2.3 $$
\par\smallpagebreak\noindent
is called the {\it Verma module} with highest weight $\Lambda.$ It is known
that
every $\g(A)$-module with highest weight $\Lambda$ is a qutoient of $M(\Lambda)$
and $M(\Lambda)$ contains a unique proper maximal submodule $M'(\Lambda).$
\par
We put
\par\smallpagebreak
$$L(\Lambda):=M(\Lambda)/M'(\Lambda).\tag2.4$$
\par\smallpagebreak
Then we can show that $L(\Lambda)$ is an irreducible $\g(A)$-module.
\par\smallpagebreak
We set
$$\align
 P:&=\{\l\in\h^*\mid <\l,\a_i^\v>\in\Bbb Z\;\;\text{for}\;\;i=1,\cdots,n\},\\
 P_+:&=\{\l\in P\mid <\l,\a_i^\v>\geq 0\;\;\text{for}\;\;i=1,\cdots,n\},\\
 P_{++}:&=\{\l\in P\mid <\l,\a_i^\v>\,>0\;\;\text{for}\;\;i=1,\cdots,n\}.
\endalign$$
\par\smallpagebreak
The set $P$ is called the {\it weight lattice.} Elements from $P$
(resp. $P_{+}$ or $P_{++}$)
are called {\it integral weights} (resp. {\it dominant} or {\it regular
dominant integral weights}). We observe that $Q\subset P$ and
$P_{++}\subset P_{+}\subset P.$ If $\Lambda$ is an element
of $P_{+},$ i.e., $\Lambda$ is a dominant integral weight, then
$L(\Lambda)$ is {\it integrable} (cf. [K], p. 171)
and the {\it Weyl-Kac character formula} for $L(\Lambda)$ is given by
\par\smallpagebreak
$$\text{ch}\,L(\Lambda)\,=\,\frac{\sum_{w\in W}\epsilon(w)e^{w(\Lambda+\rho)-\rho}}
{\prod_{\a\in\D^+}(1-e^{-\a})^{\text{mult}\,\a}}\tag2.5$$
\par\smallpagebreak
Here $\epsilon(w):=(-1)^{\ell(w)}=\det_{\h^*}w$ for $w\in W,\;W$ the
{\it Weyl group} of $\g(A)$ and $\rho$ is an element of $\h^*$ such that
$<\rho,\a_i^\v>=1$ for $i=1,\cdots,n.$ We recall that $W\subset \,\text{Aut}(\h^*)$ is
the subgroup generated by the reflections $\sigma_i(\l):=\l-\l(\a_i^\v)\a_i\;\,
(1\leq i\leq n).$\par
We set $\Lambda=0$ in (2.5). Since $L(0)$ is the trivial module over $\g(A),$
we obtain the following denominator identity or denominator formula:
\par\smallpagebreak
$$\prod_{\a\in\D^+}(1-e^{-\a})^{\text{mult}\,\a}=\sum_{w\in W}\epsilon(w)e^{w(\rho)-
\rho}.\tag2.6$$
\par\smallpagebreak
Substituting (2.6) into (2.5), we obtain another form of the Weyl-Kac charater
formula:
\par\smallpagebreak
$$\text{ch}\,L(\Lambda)=\frac{\sum_{w\in W}\epsilon(w)e^{w(\Lambda+\rho)}}
{\sum_{w\in W}\epsilon(w)e^{w(\rho)}}.\tag2.7$$
\par\smallpagebreak
Of course, in the case when $\g(A)$ is a finite dimensional semisimple Lie
algebra, then (2.7) is the classical Weyl character formula and (2.6) is
the Weyl denominator identity. We remark that an integrable highest weight
module $L(\Lambda)$ over $\g(A)$ is unitarizable and conversely if
$L(\Lambda)$ is irreducible, then $\Lambda\in P_+$(cf. [K], p.196).
\par\smallpagebreak
Let $A=(a_{ij})$ be a generalized Cartan matrix. We associate to $A$ a graph
$\Cal S(A)$ called the Dynkin diagram of $A$ as follows. If $a_{ij}a_{ji}\leq4$ and
$|a_{ij}|\geq|a_{ji}|,$ the vertices $i$ and $j$ are
connected by $|a_{ij}|$ lines, and these lines are equipped with an arrow
pointing toward $i$ if $|a_{ij}|>1.$ If $a_{ij}a_{ji}>4,$ the vertices $i$ and
$j$ are connected by a bold-faced line
equipped with an ordered pair of integers $|a_{ij}|,\,|a_{ji}|.$ We list some
of hyperbolic Kac-Moody Lie algebras.
\ex
\vskip 5mm

\SelectTips{xy}{10}
\objectmargin={-0.7pt}

$HA_1^{(1)}\ :\qquad\qquad \xygraph{!~:{@{<=>}}
    \circ ([]!{+(0,-.2)} {\alpha_{-1}}) - [r]!{-(0.1,0)}
    \circ ([]!{+(0,-.2)} {\alpha_0}) : [r]!{-(0.1,0)}
    \circ ([]!{+(0,-.2)} {\alpha_1})}$

\vskip 1cm $HA_l^{(1)},\ l\geq 2\ :\quad
\xygraph{
    \circ ([]!{+(0,-0.2)} {\alpha_1})( - []!{+(1.8,.5)} \circ( - [r]!{+(-1,0.7)}\circ([]!{+(0.35,0)} {\alpha_{-1}})([]!{+(0,-0.9)} {\alpha_{0}})) - []!{+(1.77,-0.5)}) - [r]!{-(0.1,0)}
    \circ ([]!{+(0,-0.2)} {\alpha_2}) - [r]!{-(0.1,0)}{\,\cdots\,} - [r]!{-(0.1,0)}
    \circ ([]!{+(0,-0.2)} {\alpha_{l - 1}}) - [r]!{-(0.1,0)}
    \circ ([]!{+(0,-0.2)} {\alpha_{l }})
)}
$

\vskip 1cm $HB_l^{(1)},\ l\geq 3\ :\quad \xygraph{!~:{@{=>}}
    \circ ([]!{+(0,-.2)} {\alpha_{-1}}) - [r]!{-(0.1,0)}
    \circ ([]!{+(0,-.2)} {\alpha_0}) - [r]!{-(0.1,0)}
    \circ ([]!{+(0,-.2)} {\alpha_2}) (
        - [u]!{-(0,0.3)} \circ ([]!{+(.25,0)} {\alpha_1}), - [r]!{-(0.1,0)}
    \circ ([]!{+(0,-.2)} {\alpha_3}) - [r]!{-(0.1,0)} {\,\cdots\,} - [r]!{-(0.1,0)}
    \circ ([]!{+(0,-.2)} {\alpha_{l-1}}) : [r]!{-(0.1,0)}
    \circ ([]!{+(0,-.2)} {\alpha_l})}$

\vskip 1cm $HC_l^{(1)},\ l\geq 2\ :\quad \xygraph{!~:{@{<=}}
    \circ ([]!{+(0,-.2)} {\alpha_{-1}}) - [r]!{-(0.1,0)}
    \circ ([]!{+(0,-.2)} {\alpha_0}) :@{=>} [r]!{-(0.1,0)}
    \circ ([]!{+(0,-.2)} {\alpha_1}) - [r]!{-(0.1,0)}
    \circ ([]!{+(0,-.2)} {\alpha_2}) - [r]!{-(0.1,0)} {\,\cdots\,} - [r]!{-(0.1,0)}
    \circ ([]!{+(0,-.2)} {\alpha_{l-1}}) : [r]!{-(0.1,0)}
    \circ ([]!{+(0,-.2)} {\alpha_l})}$

\vskip 1cm $HD_l^{(1)},\ l\geq 4\ :\quad \xygraph{!~:{@{=>}}
    \circ ([]!{+(0,-.2)} {\alpha_{-1}}) - [r]!{-(0.1,0)}
    \circ ([]!{+(0,-.2)} {\alpha_0}) - [r]!{-(0.1,0)}
    \circ ([]!{+(0,-.2)} {\alpha_2}) (
        - [u]!{-(0,0.3)} \circ ([]!{+(.25,0)} {\alpha_1}), - [r]!{-(0.1,0)}
    \circ ([]!{+(0,-.2)} {\alpha_3}) - [r]!{-(0.1,0)} {\,\cdots\,} - [r]!{-(0.1,0)}
    \circ ([]!{+(0,-.2)} {\alpha_{l-2}})(
        - [u]!{-(0,0.3)} \circ ([]!{+(.25,0)} {\alpha_l}), - [r]!{-(0.1,0)}
    \circ ([]!{+(0,-.2)} {\alpha_{l-1}})}$

\vskip 1cm $HF_4^{(1)}\ :\qquad\qquad \xygraph{!~:{@{=>}}
    \circ ([]!{+(0,-.2)} {\alpha_{-1}}) - [r]!{-(0.1,0)}
    \circ ([]!{+(0,-.2)} {\alpha_0}) - [r]!{-(0.1,0)}
    \circ ([]!{+(0,-.2)} {\alpha_1}) - [r]!{-(0.1,0)}
    \circ ([]!{+(0,-.2)} {\alpha_2}) : [r]!{-(0.1,0)}
    \circ ([]!{+(0,-.2)} {\alpha_3}) - [r]!{-(0.1,0)}
    \circ ([]!{+(0,-.2)} {\alpha_4})}$

\vskip 1cm $HG_2^{(1)}\ :\qquad\qquad \xygraph{!~:{@3{->}}
    \circ ([]!{+(0,-.2)} {\alpha_{-1}}) - [r]!{-(0.1,0)}
    \circ ([]!{+(0,-.2)} {\alpha_0}) - [r]!{-(0.1,0)}
    \circ ([]!{+(0,-.2)} {\alpha_1}) : [r]!{-(0.1,0)}
    \circ ([]!{+(0,-.2)} {\alpha_2}) }$

\vskip 1cm $HE_6^{(1)}\ :\qquad\qquad \xygraph{
    \circ ([]!{+(0,-.2)} {\alpha_1}) - [r]!{-(0.1,0)}
    \circ ([]!{+(0,-.2)} {\alpha_2}) - [r]!{-(0.1,0)}
    \circ ([]!{+(0,-.2)} {\alpha_3})
    (- [u]!{-(0,0.3)} \circ(- [u]!{-(0,0.3)} \circ(
        - [u]!{-(0,0.3)} \circ ([]!{+(.34,0)} {\alpha_{-1}}) )[]!{+(.25,0)} {\alpha_0}, []!{+(.25,0)} {\alpha_6} ),
    - [r]!{-(0.1,0)}
    \circ ([]!{+(0,-.2)} {\alpha_4}) - [r]!{-(0.1,0)}
    \circ ([]!{+(0,-.2)} {\alpha_5}) }$

\vskip 1cm $HE_7^{(1)} \ :\qquad \qquad\xygraph{
    \circ ([]!{+(0,-.2)} {\alpha_{-1}}) - [r]!{-(0.1,0)}
    \circ ([]!{+(0,-.2)} {\alpha_0}) - [r]!{-(0.1,0)}
    \circ ([]!{+(0,-.2)} {\alpha_1}) - [r]!{-(0.1,0)}
    \circ ([]!{+(0,-.2)} {\alpha_2}) - [r]!{-(0.1,0)}
    \circ ([]!{+(0,-.2)} {\alpha_3})(
        - [u]!{-(0,0.3)} \circ ([]!{+(.25,0)} {\alpha_7}), - [r]!{-(0.1,0)}
    \circ ([]!{+(0,-.2)} {\alpha_4}) - [r]!{-(0.1,0)}
    \circ ([]!{+(0,-.2)} {\alpha_5}) - [r]!{-(0.1,0)}
    \circ ([]!{+(0,-.2)} {\alpha_6})}$

\vskip 1cm $HE_8^{(1)} \ :\qquad \qquad\xygraph{
    \circ ([]!{+(0,-.2)} {\alpha_{-1}}) - [r]!{-(0.1,0)}
    \circ ([]!{+(0,-.2)} {\alpha_0}) - [r]!{-(0.1,0)}
    \circ ([]!{+(0,-.2)} {\alpha_1}) - [r]!{-(0.1,0)}
    \circ ([]!{+(0,-.2)} {\alpha_2}) - [r]!{-(0.1,0)}
    \circ ([]!{+(0,-.2)} {\alpha_3}) - [r]!{-(0.1,0)}
    \circ ([]!{+(0,-.2)} {\alpha_4}) - [r]!{-(0.1,0)}
    \circ ([]!{+(0,-.2)} {\alpha_5})(
        - [u]!{-(0,0.3)} \circ ([]!{+(.25,0)} {\alpha_8}), - [r]!{-(0.1,0)}
    \circ ([]!{+(0,-.2)} {\alpha_6}) - [r]!{-(0.1,0)}
    \circ ([]!{+(0,-.2)} {\alpha_7})}$

\vskip 1cm $HA_2^{(2)}\ :\qquad\qquad \xygraph{!~:{@3{<-}}
    \circ ([]!{+(0,-.2)} {\alpha_{-1}}) - [r]!{-(0.1,0)}
    \circ ([]!{+(0,-.2)} {\alpha_0}) : [r]!{-(0.1,0)}
    \circ ([]!{+(0,-.2)} {\alpha_1})}$

\vskip 1cm $HA_{2l}^{(2)},\ l\geq 2\ :\quad \xygraph{!~:{@{<=}}
    \circ ([]!{+(0,-.2)} {\alpha_{-1}}) - [r]!{-(0.1,0)}
    \circ ([]!{+(0,-.2)} {\alpha_0}) : [r]!{-(0.1,0)}
    \circ ([]!{+(0,-.2)} {\alpha_1}) - [r]!{-(0.1,0)}
    \circ ([]!{+(0,-.2)} {\alpha_2}) - [r]!{-(0.1,0)} {\,\cdots\,} - [r]!{-(0.1,0)}
    \circ ([]!{+(0,-.2)} {\alpha_{l-1}}) : [r]!{-(0.1,0)}
    \circ ([]!{+(0,-.2)} {\alpha_l})}$

\vskip 1cm $HA_{2l-1}^{(2)},\ l\geq 3 :\quad \xygraph{!~:{@{<=}}
    \circ ([]!{+(0,-.2)} {\alpha_{-1}}) - [r]!{-(0.1,0)}
    \circ ([]!{+(0,-.2)} {\alpha_0}) - [r]!{-(0.1,0)}
    \circ ([]!{+(0,-.2)} {\alpha_2}) (
        - [u]!{-(0,0.3)} \circ ([]!{+(.25,0)} {\alpha_1}), - [r]!{-(0.1,0)}
    \circ ([]!{+(0,-.2)} {\alpha_3}) - [r]!{-(0.1,0)} {\,\cdots\,} - [r]!{-(0.1,0)}
    \circ ([]!{+(0,-.2)} {\alpha_{l-1}}) : [r]!{-(0.1,0)}
    \circ ([]!{+(0,-.2)} {\alpha_l})}$

\vskip 1cm $HD_{l+1}^{(2)},\ l\geq 2\:\quad \xygraph{!~:{@{<=}}
    \circ ([]!{+(0,-.2)} {\alpha_{-1}}) - [r]!{-(0.1,0)}
    \circ ([]!{+(0,-.2)} {\alpha_0}) : [r]!{-(0.1,0)}
    \circ ([]!{+(0,-.2)} {\alpha_1}) - [r]!{-(0.1,0)}
    \circ ([]!{+(0,-.2)} {\alpha_2}) - [r]!{-(0.1,0)} {\,\cdots\,} - [r]!{-(0.1,0)}
    \circ ([]!{+(0,-.2)} {\alpha_{l-1}}) :@{=>} [r]!{-(0.1,0)}
    \circ ([]!{+(0,-.2)} {\alpha_l})}$

\vskip 1cm $HE_{6}^{(2)}\ :\qquad\qquad \xygraph{!~:{@{<=}}
    \circ ([]!{+(0,-.2)} {\alpha_{-1}}) - [r]!{-(0.1,0)}
    \circ ([]!{+(0,-.2)} {\alpha_0}) - [r]!{-(0.1,0)}
    \circ ([]!{+(0,-.2)} {\alpha_1}) - [r]!{-(0.1,0)}
    \circ ([]!{+(0,-.2)} {\alpha_2}) : [r]!{-(0.1,0)}
    \circ ([]!{+(0,-.2)} {\alpha_3}) - [r]!{-(0.1,0)}
    \circ ([]!{+(0,-.2)} {\alpha_4})}$

\vskip 1cm $HD_{4}^{(3)}\ :\qquad\qquad \xygraph{!~:{@3{<-}}
    \circ ([]!{+(0,-.2)} {\alpha_{-1}}) - [r]!{-(0.1,0)}
    \circ ([]!{+(0,-.2)} {\alpha_0}) - [r]!{-(0.1,0)}
    \circ ([]!{+(0,-.2)} {\alpha_1}) : [r]!{-(0.1,0)}
    \circ ([]!{+(0,-.2)} {\alpha_2}) }$

\par\vskip 1cm

Let $A=(a_{ij})_{i,j=-1,0,\cdots,\ell}$ be a hyperbolic generalized Cartan
matrix whose corresponding affine
submatrix of $A$ is given by $A_0=(a_{kl})_{k,l=0,1,\cdots,\ell}.$
We can realize the hyperbolic Kac-Moody Lie algebra $\g(A)$ as the minimal
graded Lie algebra $L=\oplus_{n\in\Bbb Z}L_n$ with local part
$V+\g(A_0)+V^*,$ where $V=L(-\a_{-1})$
is the basic representation of the affine Kac-Moody Lie algebra $\g(A_0)$
and $V^*$ is the contragredient of $V.$ Thus $L=G/I,$ and $L_n=G_n/I_n,$ where
$G=\oplus_{n\in\Bbb Z}G_n$ is the maximal graded
Lie algebra with local part $V+\g(A_0)+V^*$ and $I=\oplus_{n\in\Bbb Z}I_n$ is
the maximal graded ideal of $G$ intersecting
the local part trivially. Each $L_n\;(n\in\Bbb Z)$ is a $\g(A_0)$-module. (By
definition, $G=\oplus_{n\in\Bbb Z}G_n$ is called
a {\it graded Lie algebra} if $G$ is a Lie algebra and $[G_i,G_j]\subset G_{i+j}$
for all $i,j\in\Bbb Z.)$
A graded Lie algebra $G=\oplus_{n\in\Bbb Z}G_n$ is called {\it irreducible} if the
represantation $\phi_{-1}$ of $G_0$ on $G_{-1}$
defined by $\phi_{-1}(x_0)x_{-1}=[x_0,x_{-1}]$ for all $x_0\in G_0$ and
$x_{-1}\in G_{-1}$ is irreducible.
A graded Lie algebra $G=\oplus_{n\in\Bbb Z}G_n$ is said to be {\it maximal}
(resp. {\it minimal}\,) if for any other graded Lie algebra $G'=\oplus_{n\in\Bbb Z}G'_n,$
every isomorphism of the local parts of $G$ and $G'$ can be extended to an
epimorphism of $G$ onto $G'$ (resp. of $G'$ onto $G$). Kac (cf. [K])
proved that for any local Lie algebra $G_{-1}\oplus G_0\oplus G_1,$
there exist unique up to isomorphism maximal and minimal graded Lie algebras whose
local parts are isomorphic to a given Lie algebra $G_{-1}\oplus G_0\oplus G_1$.
\example{Example 2.3}\ Let
$$A=(a_{ij})_{i,j=-1,0,1}
:=\pmatrix 2&\;-1&\; 0\\
          -1&\; 2&\;-2\\
       0&\;-2&\; 2   \endpmatrix$$
be the hyperbolic generalized Cartan matrix. We can realize the corresponding
hyperbolic Kac-Moody Lie algebra $\g(A):=HA_1^{(1)}$ as the minimal graded Lie
algebra $L=\oplus_{n\in\Bbb Z}L_n$ with local part $V+\g(A_0)+V^*,$ where
$A_0=\pmatrix \;2&-2\\-2&\;2\endpmatrix$ and $V:=L(-\a_{-1})$ is
the basic representation of the affine Kac-Moody Lie algebra
$\g(A_0):=A_1^{(1)}.$ The dimensions $\dim(L_{-n})_\a$ for $0\leq n\leq 5$
were computed by A. J. Feingold, I. B. Frenkel, S.-J. Kang and
etc. For instance, $\dim(L_0)_\a=1$ and
$$\dim(L_{-1})_\a=p\left(1-\frac{(\a,\a)}{2}\right),$$
where $p$ is the partition function defined by
$$\sum_{n=0}^\infty p(n)q^n=\frac{1}{\phi(q)},\qquad\phi(q):=\prod_{n\geq 1}
(1-q^n).\tag2.8$$
\endexample
\example{Example 2.4}\ Let
$$A=(a_{ij})_{i,j=-1,0,1}
:=\pmatrix 2&\;-1&\;0\\
         - 1&\; 2&\;-4\\
       0&\;-1&\; 2   \endpmatrix$$
be the generalized Cartan matrix of hyperbolic type. We can realize
$\g(A):=HA_2^{(2)}$ as
the minimal graded Lie algebra $L=\oplus_{n\in\Bbb Z}L_n.$ The dimensions
$\dim L_{-n}\;(1\leq n\leq 3)$ were computed by A. J. Feingold
and S. J. Kang.
\endexample
\example{Example 2.5}
 \ Kac-Moody-Wakimoto (cf.[KMW]) considered the hyperbolic Kac-Moody
Lie algebra $HE_8^{(1)}=E_{10}.$ Using the modular invariant property of level
2 string functions, they computed root multiplicities of
$G_{-2}$ and $I_{-2}.$ Thus they obtained the formula
$$\dim(L_{-2})_\a=\xi\left( 3-\frac{(\a,\a)}2\right),$$
where $\xi(n)$ is defined by the relation
$$\sum_{n=0}^\infty\xi(n)q^n=\frac1{\phi(q)^8}\left(1-\frac{\phi(q^2)}
{\phi(q^4)}\right).\tag2.9 $$
\endexample
\remark{\smc Remark 2.6}\ In [Fr], Frenkel conjectured that for a hyperbolic
Kac-Moody Lie algebra $\g,$ we have
$$\dim\g_\a\leq p^{(\ell-2)}\left(1-\frac{(\a,\a)}2\right),$$
where $\ell$ is the size of the generalized Cartan matrix of $\g$ and the
function $p^{(\ell-2)}(n)$ is defined by
$$\sum_{n=0}^\infty p^{(\ell-2)}(n)q^n=\frac{1}{\phi(q)^{\ell-2}}=\prod_{n\geq 1}
(1-q^n)^{2-\ell}.\tag2.10$$
But this conjecture does not hold for $E_{10}$ (cf. [KMW]).
This conjecture is true for $HA_n^{(1)}.$ We observe that $HA_n^{(1)}$ is of
hyperbolic type for $n\leq 7$ and that $HA_n^{(1)}$ is of almost hyperbolic
type for $n\geq8.$
\endremark

\head  {\bf Appendix. \  Generalized Kac-Moody Algebras} \endhead
\par
Let $I$ be a countable index set. A real matrix $A=(a_{ij})_{i,j\in I}$
is called a {\it Borcherds-Cartan matrix} if it satisfies the following conditions:
\roster
\item"(BC1)" $a_{ii}=2$ or $a_{ii}\leq0$ for all $i\in I;$
\item"(BC2)" $a_{ij}\leq0$ if $i\neq j$ and $a_{ij}\in\Bbb Z$ if $a_{ii}=2;$
\item"(BC3)" $a_{ij}=0$ implies $a_{ji}=0.$
\endroster
Let $I^{re}:=\{i\in I\mid a_{ii}=2\}$ and $I^{im}:=\{i\in I\mid a_{ii}\leq0\}.$
Let $\m=(m_i\mid i\in I)$ be the {\it charge} of $A,$ i.e., $m_i=1$ for all
$i\in I^{re}$ and $m_j\in\Bbb Z^+$ for all $j\in I^{im}.$
A Borcherds-Cartan matrix $A$ is said to be {\it symmetrizable} if there
exists a diagonal matrix $D=\text{diag}(\delta_i\mid i\in I)$ with $\delta_i>0\;(i\in I)$
such that $DA$ is symmetric.

\definition{Definition 2.7} \ The {\it generalized Kac-Moody algebra} $\g=\g(A,\m)$
with a symmetrizable Borcherds-Cartan matrix $A$ of charge $\m=(m_i\mid i\in I)$
is the complex Lie algebra generated by the
elements $h_i,d_i\;(i\in I),\;e_{ik},f_{ik}\;(i\in I,\,k=1,\cdots,m_i)$ with
the defining relations:
$$\align
&[h_i,h_j]=[h_i,d_j]=[d_i,d_j]=0,\\
&[h_i,e_{j\ell}]=a_{ij}e_{j\ell},\;\;[h_i,f_{j\ell}]=-a_{ij}f_{j\ell},\\
&[d_i,e_{j\ell}]=\delta_{ij}e_{j\ell},\;\;[d_i,f_{j\ell}]=-\delta_{ij}f_{j\ell},\\
&[e_{ik},f_{j\ell}]=\delta_{ij}\delta_{k\ell}h_i,\\
&(ad\,e_{ik})^{1-a_{ij}}(e_{j\ell})=(ad\,f_{ik})^{1-a_{ij}}(f_{j\ell})=0\;
\quad\text{if}\;
a_{ii}=2\;\text{and}\;i\neq j,\\
&[e_{ik},e_{j\ell}]=[f_{ik},f_{j\ell}]=0\;\quad\text{if}\;a_{ii}=0
\endalign$$
for all $i,j\in I,\;k=1,\cdots,m_i,\;\ell=1,\cdots,m_j.$\par
The subalgebra $\h:=\left(\sum_{i\in I}\Bbb Ch_i\right)\oplus\left(\sum_{i\in I}
\Bbb C d_i\right)$
is called the {\it Cartan subalgebra} of $\g.$ For each $j\in I,$ we
define a {\it linear functional} $\a_j\in\h^*$ by
$$\a_j(h_i)=a_{ij},\quad\a_j(d_i)=\delta_{ij}\quad\text{for all}\;\;i,j\in I.$$
Let $\Pi:=\{\a_i\mid i\in I\}\subset\h^*$ and
$\Pi^\v:=\{h_i\mid i\in I\}\subset\h.$ The elements of $\Pi$ (resp. $\Pi^\v$)
are called the {\it simple roots} (resp. {\it simple coroots}) of $\g.$ We set
$$Q:=\sum_{i\in I}\Bbb Z\a_i,\quad Q^+:=\sum_{i\in I}\Bbb Z_+\a_i,
\quad Q^-:=-Q^+.$$
$Q$ is called the {\it root lattice} of $\g.$ We define a partial ordering
$\leq$ on $\h^*$ by $\l\leq\mu$ if $\mu-\l\in Q^+.$ For $\a\in\h^*,$ we put
$$\g_\a:=\{X\in\g\mid [h,X]=\a(h)X\,\;\;\text{for all } h\in\h\}.$$
If $\g_\a\neq0$ and $\a\neq0,\;\a$ is called {\it root} of $\g$ and $\g_\a$ is
called the {\it root space attached to the root $\a$}. The generalized
Kac-Moody algebra $\g$ has the root decomposition
$$\g=\sum_{\a\in Q}\g_\a\qquad(\text{direct sum}).\tag2.11$$
We observe that $\g_{\a_i}=\sum_{k=1}^{m_i}\Bbb Ce_{i,k}$ and
$\g_{-\a_i}=\sum_{k=1}^{m_i}\Bbb C f_{i,k}.$ The number $\text{mult}\,\a:=\dim\g_\a$
is called the {\it multiplicity} of $\a.$ A root $\a$ with $\a>0$ (with $\a<0$)
is said to be {\it positive} (resp. {\it negative}). We denote by
$\D,\D^+,\D^-$ the set of all roots, positive roots, and negative roots
respectively. We set
$$\frak n^+:=\sum_{\a\in\D^+}\g_\a,\qquad \frak n^-:=\sum_{\a\in\D^-}\g_\a.\tag2.12$$
Then we have the triangular decomposition:
$\g=\frak n^-\oplus\h\oplus\frak n^+.$\par
Since $A$ is symmetrizable, there exists a symmetric linear form $(\;|\;)$ on
$\h^*$ satisfying
the condition $(\a_i\mid \a_j)=\delta_i a_{ij}$ for all $i,j\in I.$ We say
that a root $\a$ is {\it real} if $(\a\mid\a)>0$ and {\it imaginary} if
$(\a\mid\a)\leq0.$ In particular, the simple root $\a_i$ is real if $a_{ii}=2$
and imaginary if $a_{ii}\leq0.$ We note that the imaginary simple roots may
have muliplicity $>1.$ For each $i\in I^{re},$ we let $\sigma_i\in \text{Aut}(\h^*)$
be the {\it simple reflection} on $\h^*$ defined by
$\sigma_i(\l):=\l-\l(h_i)\a_i$ for $\l\in\h^*.$ The subgroup $W$ of
$\text{Aut}(\h^*)$ generated by the $\sigma_i$'s $(i\in I^{re})$ is called the
{\it Weyl group} of $\g.$
\par
Let
$$ P^+_G:=\{\l\in\h^*\mid \l(h_i)\geq0\;\;\text{for all}\;i\in I,\;\,
\l(h_i)\in\Bbb Z_+\;\;\text{if}\;a_{ii}=2\}.$$
For $\l\in P^+_G,$ we let $V(\l)$ be the irreducible highest weight
module over $\g$ with highest weight $\l.$ We denote by $T$ the set of all
imaginary simple roots counted with multiplicities. We choose $\rho\in\h^*$
such that $\rho(h_i)=\frac12a_{ii}$ for all $i\in I.$ Then we have the
{\it Weyl-Kac-Borcherds character formula} [Bo1] :
$$\text{ch}\,V(\l)=\frac{\sum_{w\in W}\sum_{F\subset T,\,F\bot\l}(-1)^{\ell(w)+|F|}
e^{w(\l+\rho-s(F))-\rho}}{\prod_{\a\in\D^+}(1-e^{-\a})^{\text{mult}\,\a}},\tag2.13$$
where $F$ runs over all the finite subsets of $T$ such that any two distinct
elements of $F$ are mutually orthogonal, $\ell(w)$ denotes the length of $w\in W,$
 $|F|:=\text{Card}(F)$ and $s(F)$ denotes
the sum of elements in $F.$ For $\l=0,$ we obtain the {\it denominator identity}:
$$\prod_{\a\in\D^+}(1-e^{-\a})^{\text{mult}\,\a}=\sum\Sb w\in W\\ F\subset T\endSb
(-1)^{\ell(w)+|F|}e^{w(\rho-s(F))-\rho}.\tag2.14$$
\enddefinition
\remark{\smc Remark 2.8} The notion of a generalized Kac-Moody algebra was
introduced by Borcherds in his study of the vertex algebras and the
moonshine conjecture [Bo1-3]. As mentioned above, the structure and the
representation theory of generalized Kac-Moody algebras are very similar
to those of Kac-Moody algebras. The main difference is that a generalized
Kac-Moody algebra may have imaginary
{\it simple} roots.
\endremark

\def\a{\alpha}

\def\g{\gamma}
\def\G{\Gamma}

\def\l{\lambda}

\def\s{\sigma}

\def\lrt{\longrightarrow}

\def\M{{\cal M}}

\def\sh{\sharp}
\def\s{\sigma}

\def\lrt{\longrightarrow}

\def\M{{\Cal M}}

\def\BZ{\Bbb Z}
\redefine\BQ{\Bbb Q}
\def\BC{\Bbb C}
\def\BR{\Bbb R}
\def\G{\Gamma}

\head {\bf 3.\ The\ Moonshine\ Conjectures\ and\ The\ Monster\
Lie\ Alegebra}\endhead

\indent
In this section, we give a construction of the {\it Monster\ Lie\ algebra}
$M$ and a sketchy proof of the {\it Moonshine\ Conjectures} due to R. E.
Borcherds [Bo6].\par\medpagebreak
\indent
The {\it Fischer-Griess\ monster\ sporadic\ simple\ group} $G$, briefly
the MONSTER, is the largest among the $26$ sporadic finite simple groups of
order
$$2^{46}\cdot 3^{20}\cdot 5^9\cdot 7^6\cdot 11^2\cdot 13^3\cdot 17\cdot
19\cdot 23\cdot 29\cdot 31\cdot 41\cdot 47\cdot 59\cdot 71.$$
It is known that the dimension of the smallest nontrivial irreducible representation of
the MONSTER is $196883$\,(\,[FLT]\,). It was observed by John McKay\,(1939- ) that
$1+196883=196884,$ which is the first nontrivial coefficient of the elliptic
modular function $j_{\ast}(q):=j(q)-744$, where $j(q)$ is the
{\it modular\ invariant}\,:
$$j_{\ast}(q)=\sum_{n\geq -1}c(n) q^n=q^{-1}+196884 q+21493760 q^2+\cdots
.\tag 3.1$$
Later J.G. Thompson\,[Th2] found that the early coefficients of the elliptic
modular function $j_{\ast}(q)$ are simple linear combinations of the
irreducible character degrees of $G$. Motivated by these observations, J. H.
Conway and S. Norton [C-N] conjectured that there is an infinite dimensional
graded representation $V^{\sh}=\sum_{n\geq -1}\,V_n^{\sh}$ of the MONSTER
$G$ with $\text{dim}\,V_n^{\sh}=c(n)$ such that for any element $g\in G,$
the Thompson series
$$T_g(q):=\sum_{n\geq -1}\,tr(g|_{V_n^{\sh}})q^n,\ \ \ \ c_g(n):=\,
tr(g|_{V_n^{\sh}})\tag 3.2$$
is a {\it Hauptmodul} for a genus $0$ discrete subgroup of $SL(2,\BR).$
It is known that there are $194$ conjugacy classes of the MONSTER $G.$
Only $171$ of the Thompson series $T_g(q),\,g\in G$ are distinct.
Conway reports on this strange and remarkable phenomenon as follows\,:
``Because these new links are still completely unexplained, we refer to them
collectively as the `moonshine' properties of the MONSTER, intending the
word to convey our feelings that they are seen in a dim light, and that
the whole subject is rather vaguely illicit\,!\,". Therefore
the above-mentioned
conjectures had been called the {\it moonshine\ conjectures}. Recently these
conjectures were proved to be true by Borcherds [Bo6] by constructing
the so-called {\it monster\ Lie algebra}. In his proof, he made use of the
{\it natural} graded representation $V^{\sh}:=\sum_{n\geq -1}V_n^{\sh}$ of
the MONSTER $G$, called the {\it moonshine\ module} or the {\it Monster\
vertex\ algebra}, which was constructed by I.B. Frenkel, J. Lepowsky and
A. Meurman\,[FLM].\,(\,The vector space $V^{\sh}$ and $V_n^{\sh}$ are
denoted by $V^{\natural}$ and $V_{-n}^{\natural}$ respectively in [FLM].)
The {\it graded\ dimension} $\text{dim}_{\ast}\,V^{\sh}:=\sum_{n\geq -1}
(\text{dim}\,V_n^{\sh})q^n$ of the moonshine module $V^{\sh}$ is given by
$\text{dim}_{\ast}\,V^{\sh}=T_1(q)=j(q)-744.$\par\medpagebreak
\indent
Let $\Pi_{1,1}\cong \BZ^2$ be the $2$-dimensional even Lorentzian lattice
with the matrix
$\pmatrix 0 & -1\\ -1 & 0 \endpmatrix$. The {\it Monster\ Lie\ algebra}
$M$ constructed by Borcherds has the following properties\,:\par
\medpagebreak
\indent
{\bf (M1)} $M$ is a $\BZ^2$-graded generalized Kac-Moody Lie algebra with
Borcherds-Cartan matrix $A=(-(i+j))_{i,j\in I}$ of charge $\underline {m}=
((c(i)\vert\,i\in I),$ where $I=\{-1\}\cup\{ i\vert\,i\geq 1\}.$ The root
lattice of $M$ is $\Pi_{1,1}\cong \BZ^2.$\par\smallpagebreak
\indent
{\bf (M2)} $M$ is a $\BZ^2$-graded representation of the MONSTER $G$ such
that $M_{(0,0)}\cong \BR^2$ and $M_{(m,n)}\cong V_{mn}^{\sh}$ for all
$(m,n)\ne (0,0).$ In particular, $\text{dim}\,M_{(m,n)}=\text{dim}\,V_{mn}^{\sh}
=c(mn)$ for all $(m,n)\ne (0,0).$\par\smallpagebreak
\indent
{\bf (M3)} The only {\it real\ simple} root of $M$ is $(1,-1)$ and the
{\it imaginary\ simple} roots of $M$ are of the form $(1,i)$ for
$i\geq 1$ with multiplicity $c(i).$\par\smallpagebreak
\indent
{\bf (M4)} \ \ $tr(g|_{M_{(m,n)}})=tr(g|_{V_{mn}^{\sh}})=c_g(mn)$\ \ for all
$g\in G$ and $(m,n)\ne (0,0).$ \par\smallpagebreak
\indent
{\bf (M5)} $M$ has a {\it contravariant} bilinear form $(\,,\,)_0$ which
is positive definite on the piece $M_{(m,n)}$ of degree $(m,n)\ne (0,0).$
(By a contravariant bilinear form we mean that there is an involution
$\sigma$ on $M$ such that $(u,v):=-(u,\sigma(v))_0$ is invariant and
$(u,v)=0$ if $\text{deg}(u)+\text{deg}(v)\ne 0.$)\par\medpagebreak
\indent
We denote by $e_{-1}:=e_{1,-1},\,e_{i,k(i)}$ and $f_{-1}:=f_{-1,1},\,
f_{i,k(i)}\,(\,i\in I,\ k(i)=1,2,\cdots,c(i)\,)$ the positve and negative
simple root vectors of $M$ respectively. Then we have
$$\align
M_{(1,-1)} &=\BC e_{-1},\ \ \ \ M_{(-1,1)}=\BC f_{-1},\\
M_{(1,i)}  &=\BC e_{i,1}\oplus \BC e_{i,2}\oplus\cdots \oplus \BC e_{i,c(i)},\\
M_{(-1,-i)} &=\BC f_{i,1}\oplus \BC f_{i,2}\oplus\cdots \oplus
\BC f_{i,c(i)}\ \ \ (i\geq 1).\endalign$$
\indent Consider a basis of $M_{(1,i)}$ consisting of the eigenvectors
$v_{i,k(i)}(g)$ of an element $g\in G$ with eigenvalues
$\lambda_{i,k(i)}(g),$ where $k(i)=1,2,\cdots,c(i).$ Since $M_{(1,i)}\cong
V_i\,(i\geq 1)$ as $G$-modules, we have
$$\sum_{k(i)=1}^{c(i)}\lambda_{i,k(i)}(g)=\,tr(g|_{M_{(1,i)}})
=tr(g|_{V_i^{\sh}})=c_g(i) \tag 3.3$$
for all $g\in G$ and $i\geq 1.$ In addition, since $M_{(1,-1)}\cong
M_{(-1,1)}\cong V_{-1}^{\sh}\cong \BR^2$ is the trivial $G$-module, we have
$g\cdot e_{-1}=e_{-1},\ g\cdot f_{-1}=f_{-1}$ for all $g\in G.$\par
\medpagebreak
\indent For small degrees $M$ looks like Fig. 1.
\vskip 5mm
$$\matrix {}&\vdots&\vdots&\vdots&\vdots&\vdots&\vdots&\vdots&{} \\
      \cdots&  0   &  0   &  0   &   0  &  V^\sh_3 & V_6^\sh  &  V_9^\sh &\cdots \\
      \cdots&  0   &  0   &  0   &   0  &  V^\sh_2 & V^\sh_4  &  V^\sh_6 &\cdots \\
      \cdots&  0   &  0   &V^\sh_{-1}&   0  &  V^\sh_1 & V^\sh_2  &  V^\sh_3 &\cdots \\
      \cdots&  0   &  0   &  0   &\Bbb R^2& 0  &   0  &  0   &\cdots \\
      \cdots& V^\sh_3  & V^\sh_2  &  V^\sh_1 &   0    & V^\sh_{-1}& 0 &  0   &\cdots \\
      \cdots& V^\sh_6 &  V^\sh_4  &  V^\sh_2 &   0  &   0  &   0  &  0   & \cdots  \\
      \cdots& V^\sh_9 &  V^\sh_6  &  V^\sh_3 &   0  &   0  &   0  &  0   & \cdots \\
          {}&\vdots&\vdots&\vdots&\vdots&\vdots&\vdots&\vdots&{}
     \endmatrix$$
     $$ \text{\bf Fig. 1} $$
\par\vskip 5mm

\indent
Now we give a construction of the Monster Lie algebra. First of all we
define the notion of vertex algebras.\par\medpagebreak
\noindent {\bf Definition\ 3.1.} A {\it vertex\ algebra} $V$ over $\BR$ is
a real vector space with an infinite number of bilinear products, written
as $u_nv$ for $u,v\in V,\,n\in \BZ$, such that \par\smallpagebreak
\indent {\bf (V1)} \ $u_nv=0$ for $n$ sufficiently large (depending on
$u$ and $v$),\par\smallpagebreak
\indent {\bf (V2)} for all $u,v,w\in V$ and $m,n,q\in \BZ,$ we have
$$\sum_{i\in \BZ}\pmatrix m\\i\endpmatrix (u_{q+i}v)_{m+n-i}w=
\sum_{i\in \BZ}(-1)^i\pmatrix q\\i\endpmatrix \left(u_{m+q-i}(v_{n+i}w)\,-\,
(-1)^qv_{n+q-i}(u_{m+i}w)\right),$$
\indent
{\bf (V3)} there is an element $1\in V$ such that $v_n1=0$ if $n\geq 0$ and
$v_{-1}1=v.$\par\medpagebreak
\noindent
{\bf Example\ 3.2.} (1) For each even lattice $L$, there is a vertex
algebra $V_L$ associated with $L$ constructed by Borcherds [Bo1]. Let
$\hat{L}$ be the central extension of $L$ by $\BZ_2$, i.e., the double
cover of $L$. The underliying vector space of the vertex algebra $V_L$ is
given by $V_L=\BQ(\hat{L})\otimes S(\otimes_{i>0}L_i),$ where
$\BQ(\hat{L})$ is the twisted group ring of $\hat{L}$ and
$S(\oplus_{i>0}L_i)$ is the ring of polynomials over the sum of an infinite
number of copies $L_i$ of $L\otimes \BR.$\par\smallpagebreak
\noindent
(2) Let $V$ be a commutative algebra over $\BR$ with derivation $D$. Then
$V$ becomes a vertex algebra by defining
$$u_nv:=\cases {{D^{-n-1}(u)v}\over {(-n-1)!}}&\text{for $n<0$}\\
\ \ \ 0 &\text{for $n\geq 0$}.\endcases$$
Conversely any vertex algebra over $\BR$ for which $u_nv=0$ for $n\geq 0$
arises from a commutative algebra in this way. \par\smallpagebreak
\noindent
(3) Let $V$ and $W$ be two vertex algebras. Then the tensor product
$V\otimes W$ as vector spaces becomes a vertex algebra if we define
the multiplication by
$$(a\otimes b)_n(c\otimes d):=\sum_{i\in \BZ}(a_i c)\otimes (b_{n-1-i}d),
\ \ \ n\in \BZ.$$
We note that the identity element of $V\otimes W$ is given by $1_V\otimes
1_W.$\par\smallpagebreak
\noindent
(4) The moonshine module $V^{\sh}$ is a vertex algebra.
\par\medpagebreak
\noindent
{\bf Definition\ 3.3.} Let $V$ be a vertex algebra over $\BR$. A
{\it conformal\ vector} of {\it dimension} or {\it central\ charge}
$c\in \BR$ of $V$ is defined to be an element $\omega$ of $V$
satisfying the following three conditions:\par\smallpagebreak
\indent
(1) $\omega_0 v=D(v)$ \ \ for all $v\in V\,;$\par\smallpagebreak
\indent
(2) $\omega_1\omega=2\omega,\ \omega_3\omega=c/2,\ \omega_i\omega=0$\ \
for $i=2$ or $i>3$\,;\par\smallpagebreak
\indent
(3) any element of $V$ is a sum of eigenvectors of the operator
$L_0:=\omega_1$ with integral eigenvalues.\par\smallpagebreak
\noindent
Here $D$ is the operator of $V$ defined by $D(v):=v_{-2}1$ for all
$v\in V.$ If $v$ is an eigenvector of $L_0$, then its eigenvalue
$\lambda(v)$ is called the {\it conformal\ weight} of $v$ and $v$ is called
a {\it conformal\ vector} of conformal weight $\lambda(v).$ If
vertex algebras $V$ and $W$ have conformal vectors $\omega_V$ and
$\omega_W$ of dimension $m$ and $n$ respectively, then
$\omega_V\otimes \omega_W$ is a conformal vector of the vertex algebra
$V\otimes W$ of dimension $m+n.$ It is known that the vertex algebra $V$
associated with any $c$-dimensional even lattice  has a {\it canonical}
conformal vector $\omega$ of dimension $c$.
We define the operators $L_i\,(\,
i\in \BZ\,)$ of $V$ by
$$L_i:=\omega_{i+1},\ \ \ \ i\in \BZ.\tag 3.4$$
These operators satisfy the relations
$$[L_i,\,L_j]=(i-j)L_{i+j}+\pmatrix i+1\\ 3\endpmatrix {c\over 2}
\delta_{i+j,0},\ \ \ i,j\in \BZ\tag 3.5$$
and so make $V$ into a module over a {\it Virasoro\ algebra} spanned by
a central element $c$ and $L_i\,(i\in \BZ).$ We observe
that the operator $L_{-1}$ is equal to the operator $D$. We define
$$P^i=\left\{ v\in V\,|\ L_0(v)=\omega_1v=iv,\ L_k(v)=0\ \
\text{for } k>0\,\right\},\ \ i\in \BZ.\tag 3.6$$
Then the space $P^1/(DV\cap P^1)$ is a subalgebra of the Lie algebra
$V/DV$, which is equal to $P^1/DP^0$ for the vertex algebra $V_L$ or for
the Monster vertex algebra $V^{\sh}.$ Here $DV$ denotes the image of $V$
under $D$. It is known that the algebra $P^1/DP^0$ is a generalized
Kac-Moody algebra. The structure of a Lie algebra on $V/DV$ is given by
the bracket\,: $[u,v]:=u_0v (u,v\in V).$\par\medpagebreak
\indent
The vertex algebra $V_L$ associated with an even lattice $L$ has a real
valued symmetric bilinear form (\,,\,) such that if $u$ has
degree $k$, the adjoint $u_n^{\ast}$ of the operator $u_n$ with respect to
(\,,\,) is given by
$$u_n^{\ast}=(-1)^k\sum_{j\geq 0}
{{L_1^j(\s(u))_{2k-j-n-2}}\over {j!}},\tag 3.7$$
where $\s$ is the automorphism of $V_L$ defined by
$$\s(e^w):=(-1)^{(w,w)/2}(e^w)^{-1}\tag 3.8$$
for $e^w$ an element of the twisted group ring of $L$ corresponding to the
vector $w\in L.$ If a vertex algebra has a bilinear form with the
above properties, we say that {\it the bilinear\ form\ is\ compatible\
with\ the\ conformal\ vector}.\par\medpagebreak
\noindent
{\bf Definition\ 3.4.} A {\it vertex\ operator\ algebra} is defined to be
a vertex algebra with a conformal vector $\omega$ such that the
eigenspaces of the operator $L_0:=\omega_1$ are all finite dimensional and
their eigenvalues are all nonnegative integers.\par\medpagebreak
\indent
For example, the Monster vertex algebra $V^{\sh}$ is a vertex operator
algebra whose conformal vector spans the subspace $V_1^{\sh}$ fixed by
the MONSTER $G$. The vertex algebra $V_{\Pi_{1,1}}$ associated with the
$2$-dimensional even unimodular Lorentzian lattice
$\Pi_{1,1}$ is {\it not} a vertex operator algebra.
\par\medpagebreak
\indent
We recall the properties of the Monster vertex algebra $V^{\sh}.$
\par\smallpagebreak
\indent
{\bf $(V^{\sh}1)$} $V^{\sh}$ is a vertex algebra over $\BR$ with conformal
vector $\omega$ of dimension $24$ and a positive definite symmetric
bilinear form (\,,\,) such that the adjoint of $u_n$ is given by the
expression (3.7), where $\s$ is the trivial automorphism of $V^{\sh}.$
\par\smallpagebreak
\indent
{\bf $(V^{\sh}2)$} $V^{\sh}$ is a sum of eigenspaces $V_n^{\sh}$ of the
operator $L_0:=\omega_1,$ where $V_n^{\sh}$ is the eigenspace on which
$L_0$ has eigenvalue $n+1$ and $\text{dim}\,V_n^{\sh}=c(n).$ Thus
$V^{\sh}$ is a vertex operator algebra in the sense of Definition 3.4.
\par\smallpagebreak
\indent
{\bf $(V^{\sh}3)$} The MONSTER $G$ acts faithfully and homogeneously
on $V^{\sh}$ preserving the vertex
algebra structure, the conformal vector $\omega$ and the bilinear form
(\,,\,). The first few representations $V_n^{\sh}$ of the MONSTER $G$ are
decomposed as
$$\align
V_{-1}^{\sh}&=\chi_1,\ \ \ V_0^{\sh}=0,\\
V_1^{\sh}&=\chi_1+\chi_2,\\
V_2^{\sh}&=\chi_1+\chi_2+\chi_3,\\
V_3^{\sh}&=2\chi_1+2\chi_2+\chi_3+\chi_4,\\
V_5^{\sh}&=4\chi_1+5\chi_2+3\chi_3+2\chi_4+\chi_5+\chi_6+\chi_7,
\endalign $$
where $\chi_n\,(1\leq n\leq 7)$ are the first seven irreducible
representations of $G$, indexed in order of increasing dimension.
\par\smallpagebreak
\indent
{\bf $(V^{\sh}4)$} For $g\in G,$ the Thompson series $T_g(q)$ is a
{\it completely\ replicable\ function}, i.e., satisfies the identity
$$p^{-1}\text{exp}\left(-\sum_{i>0}\sum_{m\in \BZ^+,\,n\in \BZ}
tr(g^i|_{V_{mn}^{\sh}})p^{mi}q^{mi}/i\right)=
\sum_{m\in \BZ}tr(g|_{V_m^{\sh}})p^m\,-\,\sum_{n\in \BZ}
tr(g|_{V_n^{\sh}})q^n.\tag 3.9$$
We remark that the properties $(V^{\sh}1),\,(V^{\sh}2)$ and $(V^{\sh}3)$
characterize $V^{\sh}$ completely as a graded representation of $G.$
\par\medpagebreak
\indent
{\smc Construction\ of\ the\ Monster\ Lie\ algebra} $M$\,: The tensor
product $V:=V^{\sh}\otimes V_{\Pi_{1,1}}$ of $V^{\sh}$ and
$V_{\Pi_{1,1}}$ is also a vertex algebra. Then $P^1/DP^0$ is a Lie algebra
with an {\it invariant} bilinear form (\,,\,) and an involution $\tau$.
Here $P^1$ and $D:=L_1$ are defined by (3.4) and (3.6), and $\tau$ is the
involution on $V$ induced by the trivial automorphism of $V^{\sh}$ and
the involutuon $\omega$ of $V_{\Pi_{1,1}}.$ Let $R:=\{ v\in V\,|\
(u,v)=0\ \ \text{for}\ u\in V\,\}$ be the radical of (\,,\,). It is easy
to see that $DP$ is a proper subset of $R$. We define $M$ to be the
quotient of the Lie algebra $P^1/DP^0$ by $R$. The $\Pi_{1,1}$-grading
of $V_{\Pi_{1,1}}$ induces a $\Pi_{1,1}$-grading on $M$. According to the
{\it no-ghost\ theorem}, $M_{(m,n)}$ is isomorphic to the piece
$V_{mn}^{\sh}$ of degree $1-(m,n)^2/2=1+mn$ if $(m,n)\ne (0,0)$ and
$M_{(0,0)}\cong V_0^{\sh}\oplus \BR^2\cong \BR^2.$ And if $v\in M$ is
nonzero and homogeneous of nonzero degree in $\Pi_{1,1},$ then
$(v,\tau(v))>0.\ M$ satisfies the properties (M1)-(M5).
\par\medpagebreak
\noindent
{\smc Remark\ 3.5.} The construction of the Monster Lie algebra $M$ from
a vertex algebra can be carried out for any vertex algebra with a
conformal vector, but it is only when this vector has {\it dimension}
$24$ that we can apply the no-ghost theorem to identify the {\it homogeneous}
pieces of $M$ with those of $V^{\sh}.$ The important thing is that
the bilinear form (\,,\,) on $M$ is positive definite on any piece of
nonzero degree, and thus need not be true if the conformal vector has
dimension greater than $24$, even if the inner product is positive definite.
\par\medpagebreak\noindent
{\bf Problem.} Are there some other ways to construct the Monster Lie
algebra\,?\par\medpagebreak
\indent
{\smc Sketchy\ Proof\ of\ the\ Moonshine\ Conjectures\,:} The proof is
divided into two steps.\par\smallpagebreak
\noindent
{\smc Step\ I.} The Thompson series are determined by the first $5$
coefficients $c_g(i),\ 1\leq i\leq 5$ for all $g\in G$ because of the
identities (3.9).
\par\smallpagebreak
\noindent
{\smc Step\ II.} The Hauptmoduls listed in Conway and Norton [C-N, Table 2]
satisfy the identities (3.9) and have the same first $5$ coefficients of
the Thompson series.\par\medpagebreak\indent
The proof of step I is done by comparing the coefficients of $p^2$ and
$p^4$ of both sides of the identities (3.9) and so obtaining the
recursion formulas among $c_g(i).$ The proof of step II follows from
the result of Norton [No1] and Koike [Koi1]
that the modular functions associated with
elements of the MONSTER $G$ also satisfy the identites (3.9) and hence
satisfy the same recursion formulas. Roughly we explain how Conway and
Norton [C-N] associate to an element of $G$ a modular function of genus $0$.
Let $g$ be an element of $G$ corresponding to an element of odd order in
$\text{Aut}(\Lambda)$ with Leech lattice $\Lambda$ such that $g$ fixes
no nonzero vectors.
\def\ep{\epsilon}
Let $\ep_1,\cdots,\ep_{24}$ be eigenvalues of $g$ on the real vector space
$\Lambda\otimes \BR.$ We define
$$\eta_g(q):=\eta_g(\ep_1q)\cdots \eta(\ep_{24}q),\ \ \ q:=e^{2\pi i\tau},
\ \tau\in H_1,\tag 3.10$$
where $\eta(q):=q^{1/24}\prod_{n\geq 1}(1-q^n)$ is the Dedekind eta
function and $H_1:=\{ z\in \BC\vert\ \Im z >0\}$ is the Poincar${\acute e}$
upper half plane. We put
$$j_g(q):={1\over {\eta_g(q)}}-{1\over {\eta_g(0)}}.$$
Then $j_g(q)$ is the modular function of genus $0.\ j_g(q)$ is
the modular function associated with an element $g$ of $G$ by Conway and
Norton.
\vskip 1cm
\head \ {\bf Appendix\,:\ The\ No-Ghost\ Theorem} 
\endhead          
\indent Here we discuss the {\smc no-ghost\ theorem.} First of all, we
describe the concept of a {\it Virasoro\ algebra}.\par
\medpagebreak
\indent
Let ${\Bbb F}=\BR$ or $\BC.$ Let ${\Bbb F}[t,t^{-1}]$ be the commutative
associative algebra of Laurent polynomials in an indeterminate $t$, i.e.,
the algebra of finite linear combinations of integral powers of $t.$
\redefine\BF{\Bbb F}
Let $p(t)\in {\Bbb F}[t,t^{-1}]$ and we consider the derivation
$D_{p(t)}$ of $\BF[t,t^{-1}]$ defined by
$$D_{p(t)}:=\,p(t){ d\over {dt}}.\tag 1$$
The vector space $\frak a$ spanned by all the derivations of type (1) has
the Lie algebra structure with respect to the natural Lie bracket
$$[D_p,\,D_q]=D_{pq'-p'q}\ \ \ \text{for\ all}\ p,q\in \BF[t,t^{-1}].
\tag 2$$
We choose the following basis $\{ d_n|\,n\in \BZ\}$ of $\frak a$ defined by
$$d_n:=-t^{n+1}{d\over {dt}},\ \ \ n\in \BZ.\tag 3$$
By (2), we have the commutation relation
$$[d_m,d_n]=(m-n)d_{m+n},\ \ \ m,n\in \BZ.\tag 4$$
It is easily seen that $\frak a$ is precisely the Lie algebra consisting of
all derivations of $\BF[t,t^{-1}].$\par\medpagebreak
\indent
Now we consider the one-dimensional central extension $\frak b$ of
$\frak a$ by $\BF c$ with basis consisting of a central element $c$ and
elements $L_n,\,n\in \BZ,$ corresponding to the basis $d_n,\,n\in \BZ,$ of
$\frak a.$ We define the bilinear map $[\,,\,]_{\ast}:\frak b\times
\frak b\lrt \frak b$ by
$$[c,\frak b]_\ast=[\frak b,c]_\ast=[c,c]_\ast=0$$
and
$$[L_m,L_n]_\ast=(m-n)L_{m+n}+{1\over {12}}(m^3-m)\delta_{m+n,0}c\tag 5$$
for all $m,m\in \BZ.$ Then $[\,,\,]_{\ast}$ is
anti-symmetric and satisfy the
Jacobi identity. Thus $(\frak b,[\,,\,]_{\ast})$
has the Lie algebra structure.
The Lie algebra $\frak b$ is called a {\it Virasoro\ algebra.} It is not
difficult to see that the extension $\frak b$ of the Lie algebra $\frak a$
is the unique nontrivial one-dimensional central extension up to
isomorphism.\par\smallpagebreak
\indent
Now we state the {\it no-ghost\ theorem} and give its sketchy proof.
\par\smallpagebreak
\noindent
{\bf The\ No-Ghost\ Theorem.} Let $V$ be a vertex algebra
with a nondegenerate bilinear form $(\,,\,)_V.$ Suppose that $V$ is acted on
by a Virasoro algebra $\frak b$ in such a way that the adjoint of $L_k$
with respect to $(\,,\,)_V$  is $L_{-k}\,(k\in \BZ),$ the central
element of $\frak b$ {\it acts\ as\ multiplication\ by} $24$, any vector
of $V$ is a sum of eigenvectors of $L_0$ with nonnegative integral
eigenvalues, and all the eigenspaces of $L_0$ are finite dimensional. We
let $V^k:=\{ v\in V\,\vert\,L_0(v)=kv\}\,(k\in \BZ_+)$ be the
$k$-eigenspace of $L_0.$ Assume that $V$ is acted on by a group $G$ which
preserves all this strucutre. Let $V_{\Pi_{1,1}}$ be the vertex algebra
associated with the two-dimensional even unimodular Lorentzian lattice
$\Pi_{1,1}$ so that $V_{\Pi_{1,1}}$ is $\Pi_{1,1}$-graded, has a bilinear
form $(\,,\,)_{1,1}$ and is acted on by the Virasoro algebra $\frak b$ as
mentioned in this section. We let
$$P^1:=\{\,v\in V\otimes V_{\Pi_{1,1}}\,\vert\,L_0(v)=v,\ L_k(v)=0\ \
\text{for\ all}\ k>0\,\}$$
and let $P^1_r$ be the subspace of $P^1$ of degree $r\in \Pi_{1,1}.$
All these spaces inherit an action of $G$ from the action of $G$ on $V$
and the trivial action of $G$ on $V_{\Pi_{1,1}}.$ Let
$(\,,\,):=(\,,\,)_V\otimes (\,,\,)_{1,1}$ be the tensor product of
$(\,,\,)_V$ and $(\,,\,)_{1,1}$ and let
$$R:=\{\,v\in V\otimes V_{\Pi_{1,1}}\,\vert\ (u,v)=0\ \
\text{for\ all}\ u\in V\otimes V_{\Pi_{1,1}}\,\}$$
be the null space of (\,,\,). Then as $G$-modules with an invariant
bilinear form,
$$P_r^1/R\cong \cases V^{1-(r,r)/2}, &\text{for $r\ne 0$}\\
V^1\oplus \BR^2, &\text{for $r=0$}.\endcases $$
\noindent
\remark
{\smc Remark} (1) The name ``no-ghost theorem" comes from the fact that
in the original statement of the theorem in [G-T], $V$ was a part of the
underlying vector space of the vertex algebra associated with a positive
definite lattice so that the inner product on $V^i$ was positive definite,
and thus $P_r^1$ had no {\it ghosts}, i.e., {\it vectors\ of\ negative\
norm} for $r\ne 0.$\par\smallpagebreak\noindent
(2) If we take the the moonshine module $V^{\sh}$ as $V$, then $V_n^{\sh}$
corresponds to $V^{n+1}$ for all $n\in \BZ.$\par\medpagebreak
\endremark
\indent
{\smc A\ sketchy\ proof\,:} Fix a certain nonzero vector $r\in \Pi_{1,1}$
and some norm $0$ vector $w\in \Pi_{1,1}$ with $(r,w)\ne 0.$ We have
an action of the Virasoro algebra on $V\otimes V_{\Pi_{1,1}}$ generated by
its conformal vector. The operators $L_m$ of the Virasoro algebra $\frak b$
satisfy the relations
$$[L_m,L_n]_\ast=(m-n)L_{m+n}+{{26}\over {12}}(m^3-m)\delta_{m+n,0},\ \
m,n\in \BZ,\tag 6$$
and the adjoint of $L_m$ is $L_{-m}.$ Here $26$ comes from the fact that
the central element $c$ acts on $V$ as multiplication by $24$ and the
dimension of $\Pi_{1,1}$ is two. We define the operators $K_m,\,m\in \BZ$
by
$$K_m:=v_{m-1},\tag 7$$
where $v:=e_{-2}^{-w}e^w$ is an element of the vertex algebra $V_{\Pi_{1,1}}$
and $e^w$ is an element of the group ring of the double cover of
$\Pi_{1,1}$ corresponding to $w\in \Pi_{1,1}$ and $e^{-w}$ is its inverse.
Then these operators satisfy the relations
$$[L_m,K_n]_\ast=-nK_{m+n},\ \ \ [K_m,K_n]_\ast=0\tag 8$$
for all $m,n\in \BZ.$ (8) follows from the fact that the adjoint of $K_m$
is $K_{-m}$ and $(w,w)=0.$\par\smallpagebreak
\indent
Now we define the subspaces $T^1$ and $Ve^r$ of $V\otimes V_{\Pi_{1,1}}$ by
$$T^1:=\{\,v\in V\otimes V_{\Pi_{1,1}}\,\vert\,\text{deg}(v)=r,\,
L_0(v)=v,\,L_m(v)=K_m(v)=0\ \text{for\ all}\ m>0\,\}$$
and $Ve^r:=V\otimes e^r.$ Then we can prove that
$$T^1\cong V^{1-(r,r)/2}e^r\ \ \text{and}\ \ T^1\cong P^1/R.\tag 9$$
We leave the proof of (9) to the reader.
Consequently we have the desired result
$$P^1/R\cong V^{1-(r,r)/2}e^r\cong V^{1-(r,r)/2}.$$
For the case $r=0$, we leave the detail to the reader.
\par\smallpagebreak\indent
Finally we remark that in [FI] Frenkel uses the no-ghost theorem to
prove some results about Kac-Moody algebras.

\redefine\g{\frak g}
\redefine\h{\frak h}
\redefine\a{\alpha}
\redefine\v{\vee}
\redefine\l{\lambda}
\redefine\D{\Delta}
\redefine\m{\underline m}

\redefine\Sp{Sp(g,\Bbb R)}
\redefine\C{\Bbb C}
\redefine\R{\Bbb R}
\redefine\la{\lambda}
\redefine\ka{\kappa}

\head {\bf 4. \  Jacobi Forms}  
\endhead  
In this section, we discuss Jacobi forms associated to the symplectic group
$\Sp$ and those associated to the orthogonal group $O_{s+2,2}(\Bbb R)$
respectively. We also discuss the differences between them.
\newline\ex
{\bf I. Jacobi forms associated to $\Sp.$}
\par\ex
An exposition of the theory of Jacobi forms associated to the symplectic
group $\Sp$ can be found in [E-Z], [Y1]-[Y4] and [Zi].
\par
In this subsection, we establish the notations and define the concept of
Jacobi forms associated to the symplectic group. For any positive integer
$g\in \BZ^+,$ we let
$$\Sp = \{M \in \R^{(2g,2g)}\mid {}^t\! MJ_gM=J_g \}$$
be the symplectic group of degree $g,$ where
$$J_g:=\pmatrix 0 & E_g \\ -E_g & 0 \endpmatrix.$$
It is easy to see that $\Sp$ acts on $H_g$ transitively by
$$M<Z>:=(AZ+B)(CZ+D)^{-1},$$
where $M=\pmatrix A & B \\ C & D\endpmatrix \in \Sp$ and $Z\in H_g.$
For two positive integers $g$ and $h,$ we consider the {\it Heisenberg group}
$$H_{\R}^{(g,h)}:=\{[(\lambda,\mu),\kappa]\mid\la,\mu\in\R^{(h,g)},
\ka\in\R^{(h,h)},\ka+\mu{}^t\!\la\;\text{symmetric}\}$$
endowed with the following multiplication law
$$[(\lambda,\mu),\kappa]\circ[(\lambda',\mu'),\kappa']:=
[(\lambda+\la',\mu+\mu'),\kappa+\ka'+\la{}^t\!\mu'-\mu{}^t\!\la'].$$
We define the semidirect product of $\Sp$ and $H^{(g,h)}_{\R}$
$$G^J:=\Sp\ltimes H_{\R}^{(g,h)}$$
endowed with the following multiplication law
$$\align
(M,&[(\lambda,\mu),\kappa])\cdot(M',[(\lambda',\mu'),\kappa'])\\
&:=(MM',[(\widetilde\lambda+\la',\widetilde\mu+\mu'),\kappa+\ka'+\widetilde\la
{}^t\!\mu'-\widetilde\mu{}^t\!\la']),
\endalign$$
with $M,\;M'\in \Sp$ and $(\widetilde\la,\widetilde\mu):=
(\la,\mu)M'.$ It is easy to see that $G^J$ acts on
$H_g\times \Bbb C^{(h,g)}$ transitively by
$$(M,[(\lambda,\mu),\kappa])\cdot(Z,W):=(M<Z>,(W+\la Z+\mu)(CZ+D)^{-1}),\tag4.1$$
where $M=\pmatrix A & B \\ C & D\endpmatrix \in \Sp,\;
[(\lambda,\mu),\kappa]\in H_{\Bbb R}^{(g,h)}$ and $(Z,W)\in H_g\times \Bbb C^{(h,g)}.$
\par
Let $\rho$ be a rational representation of $GL(g,\Bbb C)$ on a finite dimensional complex vector
space $V_\rho.$ Let $\Cal M\in\R^{(h,h)}$ be a symmetric half-integral semipositive definite matrix
of degree $h.$ Let $C^\infty(H_g\times {\Bbb C}^{(h,g)},V_\rho)$
be the algebra of all $C^\infty$ functions on $H_g\times \Bbb C^{(h,g)}$
with values in $V_\rho.$
For $f\in C^\infty(H_g\times \Bbb C^{(h,g)},V_\rho),$ we define
$$\aligned
(f|_{\rho,\Cal M}&[(M,[(\lambda,\mu),\kappa])])(Z,W)\\
&:=e^{-2\pi i\sigma(\Cal M[W+\la Z+\mu](CZ+D)^{-1}C)}\cdot
e^{2\pi i\sigma(\Cal M(\la Z{}^t\!\la+2\la{}^t\!W+(\ka+\mu{}^t\!\la)))}\\
&\times\rho(CZ+D)^{-1}f(M<Z>,(W+\la Z+\mu)(CZ+D)^{-1}),
\endaligned\tag4.2$$
where $M=\pmatrix A & B \\ C & D\endpmatrix \in \Sp$ and $[(\la,\mu),\ka]\in
H_\R^{(g,h)}.$

\redefine\Z{\Bbb Z}

\definition {Definition 4.1} \ Let $\rho$ and $\Cal M$ be as above. Let
$$H_\Z^{(g,h)}:=\{[(\la,\mu),\ka]\in H_{\R}^{(g,h)}\mid \la,\mu\in\Z^{(h,g)}
,\ka\in\Z^{(h,h)}\}.$$
A {\it Jocobi form} of index $\Cal M$ with respect to $\rho$ on $\Gamma_g$ is
a holomorphic function $f\in C^\infty (H_g\times \Bbb C^{(h,g)},V_\rho)$
satisfying the following conditions (A) and (B):
\roster
\item"(A)" $f|_{\rho,\Cal M}[\widetilde\gamma]=f$ for all $\widetilde\gamma\in
\Gamma_g^J:=\Gamma_g\ltimes H_\Z^{(g,h)}.$
\item"(B)" $f$ has a Fourier expansion of the following form
$$f(Z,W)=\sum\Sb T\geq0\\ \text{half-integral}\endSb\sum_{R\in\Z^{(g,h)}}
c(T,R)\cdot e^{2\pi i\sigma(TZ)}\cdot e^{2\pi i\sigma(RW)}$$
with $c(T,R)\neq 0$ only if $\pmatrix T&\frac12 R\\ \frac12{}^t\!R&\Cal M
\endpmatrix\geq0.$ \par
Moreover if $c(T,R)\neq 0$ implies
$\pmatrix T&\frac12 R\\ \frac12{}^t\!R&\Cal M
\endpmatrix>0,$ $f$ is called a {\it cusp Jacobi form.}
\endroster
\enddefinition
\par
If $g\geq2,$ the condition (B) is superfluous by the Koecher principle (cf. [Zi],
Lemma 1.6). We denote by $J_{\rho,\Cal M} (\Gamma_g)$ the vector space of all
Jacobi forms of index $\Cal M$ with respect to  $\rho$ on $\Gamma_g.$
In the special case $V_\rho = \Bbb C,$
$\rho(A)=(\text{det}A)^k (k\in \Bbb Z,\,A\in GL(g,\Bbb C)),$ we write
$J_{k,\Cal M} (\Gamma_g)$ instead of $J_{\rho,\Cal M} (\Gamma_g)$
and call $k$ the {\it weight} of a Jacobi form
$f\in J_{k,\Cal M}(\Gamma_g).$\par
Ziegler (cf. [Zi], Theorem 1.8 or [E-Z], Theorem 1.1) proves that the vector space
$J_{\rho,\Cal M} (\Gamma_g)$ is finite dimensional.
\par
\definition { Definition 4.2} \ A Jacobi form $f\in J_{\rho,\Cal M} (\Gamma_g)$
is said to be {\it singular}
if it admits a Fourier expansion such that a Fourier coefficient $c(T,R)$
vanishes unless det\,$\pmatrix T&\frac12 R\\ \frac12{}^t\!R&\Cal M
\endpmatrix=0$.
\enddefinition
\example{\bf Example 4.3}
\ Let $S\in \Z^{(2k,2k)}$ be a symmetric, positive definite, unimodular even
integral matrix and $c\in \Z^{(2k,h)}.$ We define the theta series
$$\vartheta_{S,c}^{(g)}(Z,W):=\sum_{\la\in\Z^{(2k,g)}}e^{\pi\{\sigma(S\la
Z{}^t\!\la)+2\sigma({}^t\!cS\la{}^t\!W)\}},\quad Z\in H_g,\;W\in\Bbb C^{(h,g)}.
\tag4.3$$
We put $\Cal M:=\frac12{}^t\!cSc.$ We assume that
$2k < g + \text{rank}(\Cal M).$ Then it is easy
to see that $\vartheta_{S,c}^{(g)}$ is a singular Jacobi form in
$J_{k,\Cal M} (\Gamma_g)$ (cf. [Zi], p. 212).
\endexample
\remark {\smc Remark 4.4}
\ Without loss of generality, we may assume that $\Cal M$ is a
{\it positive definite} symmetric, half-integral matrix of degree $h$
(cf. [Zi], Theorem 2.4).
\endremark
>From now on, throughout this paper $\Cal M$ is assumed to be positive definite.
\definition { Definition 4.5}
\ An irreducible finite dimensional representation $\rho$ of $GL(g,\Bbb C)$
is determined uniquely by its highest weight
$(\lambda_1,\lambda_2,\cdots,\lambda_g)\in \Z^g$ with $\lambda_1 \geq
\lambda_2 \geq\cdots \geq\lambda_g.$
We denote this representation by $\rho = (\lambda_1,
\lambda_2, \cdots ,\lambda_g).$ The number $k(\rho):=\lambda_g$
is called the {\it weight} of $\rho.$
\enddefinition
\par
The author (cf. [Y3]) proved that singular Jacobi forms in
$J_{\rho,\Cal M} (\Gamma_g)$ are characterized by their singular weights.
\proclaim {Theorem 4.6 (Yang [Y3])}
\ Let $2\Cal M$ be a symmetric, positive definite even integral
matrix of degree $h.$ Assume that $\rho$ is an irreducible representation of
$GL(g,\Bbb C).$
Then a nonvanishing Jacobi form in $J_{\rho,\Cal M} (\Gamma_g)$ is singular
if and only if $2k(\rho) < g+h.$
Only the nonnegative integers $k$ with $0\leq k\leq\frac{g+h}2$ can be the
weights of singular Jacobi forms
in $J_{k,\Cal M} (\Gamma_g).$  These integers are called singular weights
in  $J_{k,\Cal M} (\Gamma_g).$
\endproclaim
\demo{Proof}
\ The proof can be found in [Y3], Theorem 4.5. \hfill $\square$
\enddemo
\par\medpagebreak
{\bf II. Jacobi forms associated to $O_{s+2,2}(\Bbb R)$} \par\ex
An exposition of the theory of Jacobi forms associated to the orthogonal group
can be found in [Bo7], [G1] and [G2].
\par
First we fix a positive integer $s.$ We let $L_0$ be a positive definite even
integral lattice with a quadratic form $Q_0$ and let $\Pi_{1,1}$ be the
nonsingular even integral lattice with its associated symmetric matrix
$I_2:= \pmatrix 0 & -1 \\-1 & 0\endpmatrix.$
We define the following lattices $L_1$ and $M$ by
$$L_1:=L_0\oplus \Pi_{1,1}\quad\text{and}\quad M:=\Pi_{1,1}\oplus L_1.\tag4.4$$
Then $L_1$ and $M$ are nonsingular even integral lattices of $(s+1,1)$
and $(s+2,2)$
respectively. From now on we denote by $Q_0,$ $Q_1,$ $Q_M$
(resp. $(\,,\,)_0,\;(\,,\,)_1,\;(\,,\,)_M$) the quadradic forms
(resp. the nondegenerate symmetric bilinar forms)
associated with the lattices $L_0,$ $L_1,$ $M$ respectively.
We also denote by $S_0,$ $S_1,$ and $S_M$ the nonsingular symmetric even
integral matrices associated with the lattices $L_0,$ $L$ and $M$ respectively.
Thus $S_1$ and $S_M$ are given by
$$S_1:=\pmatrix 0 & 0 & -1 \\ 0 & S_0 & 0 \\ -1  & 0 & 0 \endpmatrix\;\;\text{and}\;\;
S_M:=\pmatrix 0 & 0 & I_2 \\ 0 & S_0 & 0 \\ I_2  & 0 & 0 \endpmatrix,\tag4.5$$
where $I_2:=\pmatrix 0 & -1 \\ -1 & 0\endpmatrix.$
We let $M_{\R}:=M\otimes_\Z\R$ and $M_\Bbb C:=M\otimes_\Z\Bbb C$ be the quadratic
spaces over $\Bbb R$ and $\Bbb C$ respectively. We let
$$O(M_{\R},S_M):=\{g\in GL(M_\R)\mid {}^t\!gS_Mg=S_M\}\tag4.6$$
be the real orthogonal group of the quadratic space $(M_\Bbb R,Q_M).$
We denote by $O_M$ the isometry group of the lattice $(M,Q_M).$
Then $O_M$ is an algebraic group defined over $\Z.$
We observe that $S_M$ is congruent to $E_{s+2,2}$ over $\Bbb R,$ i.e.,
$S_M={}^t\!a E_{s+2,2}a$ for some $a\in GL(s+4,\Bbb R),$ where
$$E_{s+2,2}:=\pmatrix E_{s+2} & 0 \\ 0  & -E_2 \endpmatrix.\tag4.7$$
Then it is easy to see that $O(M_\Bbb R,S_M)=a^{-1} O(M_{\Bbb R},E_{s+2,2})a.$
Now for brevity we write $O(M_{\Bbb R})$ simply
instead of $O(M_\Bbb R,S_M).$ Obviously $O(M_{\Bbb R})$ is isomorphic to the
real orthogonal group
$$O_{s+2,2}(\Bbb R):=\{g\in GL(s+4,\Bbb R)\mid {}^t\! gE_{s+2,2} g=E_{s+2,2} \}.
\tag4.8$$
$O(M_{\Bbb R})$ has four connected components. Let $G^0_{\Bbb R}$ be the
identity component of $O(M_{\Bbb R})$ and let $K^0_{\Bbb R}$
be its maximal compact subgroup. Then the pair $(G^0_{\Bbb R},K^0_{\Bbb R})$
of the real semisimple groups isomorphic to the pair
$(SO(s+2,2)^0,SO(s+2,\Bbb R)\times SO(2,\Bbb R))$ is a symmetric pair of type (BDI)
(cf. [H] 445-446).
The homogeneous space $X:=G_\R^0/K_\R^0$ is a Hermitian symmetric space of
noncompact type of dimension $s+2$ (cf. see Appendix C). Indeed, $X$ is a
bounded symmetric domain of type IV in the Cartan classification.
It is known that $X$ is isomorphic to a $G^0_{\Bbb R}$-orbit
in the projective space $\Bbb P(M_{\Bbb C}).$ Precisely, if we let
$D:=\{z\in \Bbb P(M_{\Bbb C}) \mid (z,z)_M=0,\;(z,\overline z )_M < 0\},$ then
$$D\cong G_\R^0x_0\cup G_\R^0\overline{ x_0}=D^+\cup\overline{D^+},\;\;D^+:=
G_\R^0x_0,\tag4.9$$
where $\overline{ x_0}$ denotes the complex conjugation of $x_0$ in
$\Bbb P(M_{\Bbb C}).$
We shall denote by $G_{\Bbb R}$ the subgroup of $O(M_{\Bbb R})$ preserving
the domain $D^+.$ It is known that $D^+\cong G_\R^0/K_\R^0$ may be realized
as a tube domain in $\Bbb C^{s+2}$ given by
$$\Cal D:=\left\{ {}^t\!Z=(\omega,z,\tau)\in \Bbb C^{s+2} \mid
\omega\in H_1,\; \tau\in H_1, \;S_1[\text{Im}\, Z] < 0 \right\},\tag 4.10$$
where $\text{Im}\,Z$ denotes the imaginary part of the column vector $Z.$
An embedding of the tube domain $\Cal D$ into the projective space $\Bbb P(M_{\Bbb C}),$
called the Borel embedding, is of the following form
$$p(Z)=p({}^t\!(\omega,z_1,\cdots,z_s,\tau))={}^t(\frac12 S_1[Z]:\omega:z_1:\cdots
:z_s:\tau:1)\in\Bbb P(M_\Bbb C).\tag4.11$$
$G_{\Bbb R}$ acts on $\Cal D$ transitively as follows: if $g=(g_{kl})\in G_{\Bbb R}$
with $1\leq k,l\leq s+4$ and $Z={}^t\!(\omega,z_1,\cdots,z_s,\tau)\in \Cal D,$
then
$$g<Z>:=(\widetilde \omega,\widetilde z_1,\cdots,\widetilde z_s,\widetilde \tau),
\tag4.12$$
where
$$\align
&\widetilde \omega:=(\frac12g_{2,1}S_1[Z]+g_{2,2}\omega+\sum_{l=3}^{s+2}g_{2,l}z_{l-2}+
g_{2,s+3}\tau+g_{2,s+4})\,J(g,Z)^{-1},\\
&\widetilde z_k:=(\frac12g_{k+2,1}S_1[Z]+g_{k+2,2}\omega+\sum_{l=3}^{s+2}g_{k+2,l}z_{l-2}+
g_{k+2,s+3}\tau+g_{k+2,s+4})\,J(g,Z)^{-1},\;1\leq k\leq s,\\
&\widetilde \tau:=(\frac12g_{s+3,1}S_1[Z]+g_{s+3,2}\omega+
\sum_{l=3}^{s+2}g_{s+3,l}z_{l-2}+
g_{s+3,s+3}\tau+g_{s+3,s+4})\, J(g,Z)^{-1}.
\endalign$$
Here we put
$$J(g,Z):=\frac12g_{s+4,1}S_1[Z]+g_{s+4,2}\omega+\sum_{l=3}^{s+2}g_{s+4,l}
z_{l-2}+g_{s+4,s+3}\tau+g_{s+4,s+4}.\tag4.13$$
It is easily seen that
$$p(g<Z>)J(g,Z)=g\cdot p(Z)\quad(\,\cdot\,\;\text{is the matrix
multiplication})\tag4.14$$
and that $J:G_{\Bbb R}\times \Cal D\rightarrow GL(1,\Bbb C)=\Bbb C^\times$
is the automorphic factor, i.e.,
$$J(g_1g_2,Z)=J(g_1,g_2<Z>)J(g_2,Z)$$
for all $g_1,\;g_2\in G_{\Bbb R}$ and $Z\in \Cal D.$
\par
Let $O_M(\Bbb Z)$ be the isometry group of the lattice $M.$
Then $\Gamma_M:=G_{\Bbb R}\cap O_M(\Bbb Z)$ is an arithmetic subgroup
of $G_{\Bbb R}.$
\definition { Definition 4.7}
\ Let $k$ be an integer. A holomorphic function $f$ on $\Cal D$ is a modular
form of weight $k$ with respect to $\Gamma_M$ if it satisfies the following
transformation behaviour
$$(f|_{k}\gamma)(Z):=J(\gamma,Z)^{-k}f(\gamma<Z>)=f(Z)\tag4.15$$
for all $\gamma\in\Gamma_M$ and $Z\in\Cal D.$
For a subgroup $\Gamma$ of $\Gamma_M$ of finite index, a modular form
with respect to $\Gamma$ can be defined in the same way.
\enddefinition
We denote by $M_k(\Gamma)$ the vector space consisting of all modular forms of weight $k$
with respect to $\Gamma.$
We now introduce the concept of cusp forms for $\Gamma_M.$
First of all we note that the realization $D^+$ of our tube domain $\Cal D$ in the
projective space $\Bbb P(M_{\Bbb C})$ is obtained as a subset of the quadric
$D$ in $\Bbb P(M_{\Bbb C})$ (cf. see (4.10)).
A maximal connected complex analytic set $X$ in $\overline{D^+}\smallsetminus D^+$
is called a {\it boundary component} of $D^+,$ where $\overline{D^+}$
denotes the closure of $D^+$
in $\Bbb P(M_{\Bbb C}).$ The normalizer $N(X):=\{g\in G_{\Bbb R} \mid
g(X)=X\}$ of a boundary component $X$ of $D^+$ is a maximal parabolic
subgroup of $G_{\Bbb R}.$ $X$ is called a {\it rational boundary component}
if the normalizer $N(X)$ of $X$ is defined over $\Bbb Q.$
A modular form with respect to $\Gamma_M$ is called a {\it cusp form}
if it vanishes on every rational boundary component of $D^+.$ It is well
known that any rational boundary component $X$ of $D^+$ corresponds to a
primitive isotropic sublattice $S$ of $M$ via
$X=X_S:=\Bbb P(S\otimes\Bbb C)\cap \overline{D^+}.$ Since the lattice $M$
contains only isotropic lines and planes, there
exist two types of rational boundary components, which are points and curves.
\par
The orthogonal group $G_\Bbb R$ has the rank two and so there are two types of maximal
parabolic subgroups in $\Gamma_M.$ Therefore there are two types of Fourier expansions
of modular forms. A subgroup of $\Gamma_M$ fixing a null sublattice of $M$ of
rank one is called a {\it Fourier group.} A subgroup of $\Gamma_M$ fixing a null sublattice of
$M$ of rank two
is called a {\it Jacobi parabolic group}\footnote"*"{In [Bo7], this group was named
just a Jacobi group. The definition of a Jacobi group is different from ours.}.
Both a Fourier group and a Jacobi parabolic
group are maximal parabolic subgroups of $\Gamma_M.$
\par
Let $f\in M_k(\Gamma_M)$ be a modular form of weight $k$ with respect
to $\Gamma_M.$ Since the following $\gamma_\ell(\ell\in L_1\cong \Z^{s+2})$
defined by
$$\gamma_\ell:=\pmatrix 1 &{}^t\!\alpha&b\\ 0&E_{s+2}&\ell\\0&0&1\endpmatrix,
\quad b:=\frac12S_1[\ell],\;\alpha=S_1\ell\tag4.16$$
are elements of $\Gamma_M,\;f(Z+\ell)=f(Z)$ for all $\ell\in L_1.$
We note that $\gamma_\ell(Z)=Z+\ell$ for all $Z\in \Cal D$ and $J(\gamma_\ell,Z)=1.$
Hence we have a Fourier expansion
$$f(Z)=\sum_\ell a(\ell)e^{2\pi i({}^t\!\ell S_1Z)},\tag4.17$$
where $\ell$ runs over the set $\{\ell\in \widetilde{L_1}\mid i\ell\in \Cal
D,\;S_1[\ell]\geq 0\}.$
Here $\widehat{L_1}$ denotes the dual lattice of $L_1,$ that is,
$$\widehat{L_1}:=\{\ell\in L_1\otimes_\Z\Bbb Q\mid {}^t\!\ell S_1\alpha\in \Z
\;\text{for all}\; \alpha\in L_1\}.$$
We let
$$f(Z)=f(\omega,z,\tau)=\sum_{m\geq 0}\phi_m(\tau,z)e^{2\pi im\omega}\tag4.18$$
be the Fourier-Jacobi expansion of $f$ with respect to the variable $w.$
Obviously the Fourier-Jacobi coefficient
$$\phi_0(\tau,z)=\lim_{v\rightarrow\infty}f(iv,z,\tau)\tag4.19$$
depends only on $\tau.$ We can show that the Fourier-Jacobi coefficients
$\phi_m(\tau,z)\;(m\geq 0)$ satisfies the following functional equations
$$\phi_m\biggl(\frac{a\tau+b}{c\tau+d},\frac{z}{c\tau+d}\biggl)=
(c\tau+d)^ke^{\pi im\frac{cS_0[z]}{c\tau+d}}\phi_m(\tau,z)\tag4.20$$
and
$$\phi_m(\tau,z+x\tau+y)=e^{-2\pi im({}^t\!xS_0z+\frac12S_0[x]\tau)}\phi_m(\tau,z)
\tag4.21$$
for all $\pmatrix a&b\\c&d\endpmatrix\in\Gamma_1=SL(2,\Z)$ and all $x,\;y\in \Z^s.$
\par
Now we define the Jacobi forms associated to the orthogonal group. First we choose
the following basis of $M$ such that
$$M=\Z e_1\oplus\Z e_2\oplus L_0\oplus\Z e_{-2}\oplus\Z e_{-1},\tag4.22$$
where $e_1,e_2,e_{-1},e_{-2}$ are four isotropic vectors with $(e_i,e_j)=\delta_{i,-j}.$
Let $P_\Bbb R$ be the Jacobi parabolic subgroup of $G_\Bbb R$ preserving
the isotropic plane $\Bbb R e_1\oplus\Bbb R e_2.$ Then it is easily seen that
an element $g$ of $P_\Bbb R$ is given by the following form:
$$g=\pmatrix A^0&X_1&Y\\0&U&X\\0&0&A\endpmatrix,\;\ \ X_1\in \Bbb R^{(2,s)},
Y\in \Bbb R^{(2,2)},X\in \Bbb R^{(s,2)},\tag4.23$$
$$\align
&A\in GL_2(\Bbb R)^+,\; \;\; S_0[U]=S_0,\;\;\;A^0=I{}^t\!A^{-1}I,\\
&X_1=I{}^t\!A^{-1}{}^t\!XS_0U,\;\;\;{}^t\!YIA+{}^t\!AIY=S_0[X],
\endalign$$
where $I:=I_2=\pmatrix 0&-1\\-1&0\endpmatrix.$ We denote by $GO_L$ for the
general orthogonal group or conformal group
of the lattice $L$ consisting of linear transformations multiplying the quadratic
form by an invertible element of a lattice $L.$ We let $K:=
\Z e_1\oplus\Z e_2$ be the 2-dimensional
primitive null sublattice of $M.$ Then we have a homomorphism
$\pi_P:P_\Bbb R\rightarrow GO_K(\Bbb R)\times GO_{L_0}(\Bbb R)$ defined by
$$\pmatrix A^0&X_1&Y\\0&U&X\\0&0&A\endpmatrix\longmapsto(A^0,U).\tag4.24$$
The connected component of the kernel of $\pi_P$ is called a {\it Heisenberg group,}
denoted by Heis$(M_\Bbb R).$ It is easy to see that Heis$(M_\Bbb R)$
consists of the following elements
$$\{X;r\}=\{x,y;r\}=\pmatrix 1&0&{}^t\!yS_0&{}^t\!xS_0y-r&\frac12S_0[y]\\
0&1&{}^t\!xS_0&\frac12S_0[x]&r\\0&0&E_s&x&y\\ 0&0&0&1&0\\ 0&0&0&0&1
\endpmatrix,\tag4.25$$
where $X=(x,y),$ with $x,y\in\Bbb R^{(s,1)}$ and $r\in\Bbb R.$
The multiplication on Heis($M_\Bbb R)$ is given by
$$\{X_1;r_1\}\{X_2;r_2\}:=\{X_1+X_2;r_1+r_2+{}^t\!x_1S_0y_2\},\;X_1=(x_1,y_1),\;
X_2=(x_2,y_2).\tag4.26$$
We let $G^J_\Bbb R$ be the subgroup of $P_\Bbb R$ generated by the following elements
$$\{A\}:=\text{diag}(A^0,E_s,A),\;\;A\in SL(2,\Bbb R),\;A^0=I{}^t\!A^{-1}I\tag4.27$$
and $\{X;r\}$ in Heis$(M_\Bbb R).$ $G^J_\Bbb R$ is called the (real)
{\it Jacobi group} of
the lattice $M.$
We observe that $G^J_\Bbb R$ is isomorphic to the semidirect product of
$SL(2,\Bbb R)$ and Heis$(M_\Bbb R).$
It is easy to check that
$$\{X;r\}\{A\}=\{A\}\{XA;r+\frac12({}^t\!x_AS_0y_A-{}^t\!xS_0y)\},\tag4.28$$
where $x_A$ and $y_A$ are the columns of the matrix $XA.$
We see easily that Heis$(M_\Bbb R)$ is a
normal subgroup of the Jacobi group $G^J_\Bbb R$ and the center $C^J$ of
$G^J_\Bbb R$ consists of all
elements $\triangledown(r):=\{0,0;r\},\;r\in\Bbb R.$ According to (4.12),
the actions of $\{A\}$ and $\{x,y;r\}$ on $\Cal D$ are given
by as follows:
$$\align
&\{A\}<Z>={}^t\!\biggl(\omega-\frac{cS_0[z]}{2(c\tau+d)},\frac{{}^t\!z}{c\tau+d},
\frac{a\tau+b}{c\tau+d}\biggl);\tag4.29\\
&\{x,y;r\}<Z>={}^t\!(\omega+r+{}^t\!xS_0z+\frac12 S_0[x]\tau,{}^t\!(z+x\tau+y),\tau),
\tag4.30
\endalign$$
where $A=\pmatrix a&b\\c&d\endpmatrix\in SL(2,\Bbb R)$ and $Z={}^t\!(\omega,z,\tau)\in
\Cal D.$
From (4.29) and (4.30), we can define the action of the Jacobi group $G^J_\Bbb R$ on the
$(\tau,z)$-domain $H_1\times\Bbb C^{s},$ which we denote by $g<(\tau,z)>,\,g\in
G^J_\Bbb R.$
\par\smallpagebreak
Let $k$ and $m$ be two integers. For $g\in G^J_\Bbb R$ and $Z={}^t\!(\omega,
{}^t\!z,\tau)\in\Cal D,$ we denote by $\omega(g;Z)$ the $\omega$-component of $g<Z>.$
Now we define the mapping $J_{k,m}:G^J_\Bbb R\times(H_1\times\Bbb C^s)\rightarrow
GL(1,\Bbb C)=\Bbb C^\times$ by
$$J_{k,m}(g,(\tau,z)):=J(g,Z)^ke^{-2\pi im\omega(g;Z)}\cdot e^{2\pi im\omega},\tag4.31$$
where $g\in G^J_\Bbb R,\;Z={}^t\!(\omega,{}^t\!z,\tau)\in\Cal D$ and $J(g,Z)$
is the automorphic
factor defined by (4.13). $J_{k,m}$ is well-defined, i.e., it is independent of the choice
of $Z={}^t\!(w,{}^t\!z,\tau)\in\Cal D$ with given $(\tau,z)\in H_1\times\Bbb C^s.$
It is easy to check that $J_{k,m}$ is an automorphic factor for
the Jacobi group $G^J_\Bbb R.$ In particular, we have
$$
J_{k,m}(\{A\},(\tau,z))=e^{\pi i\frac{cmS_0[z]}{c\tau+d}}\cdot(c\tau+d)^k,\;A=
\pmatrix a&b\\c&d\endpmatrix\in SL(2,\Bbb R)\tag4.32$$
and
$$J_{k,m}(\{x,y;r\},(\tau,z))=e^{-2\pi im(r+{}^t\!xS_0z+\frac12S_0[x]\tau)}.
\tag4.33$$
We have a natural action of $G^J_\Bbb R$ on the algebra $C^\infty(H_1\times\Bbb
C^s)$ of all $C^\infty$ functions on $H_1\times\Bbb C^s$ given by
$$(\phi|_{k,m}g)(\tau,z):=J_{k,m}(g,(\tau,z))^{-1}\phi(g<(\tau,z)>),\tag4.34$$
where $\phi\in C^\infty(H_1\times C^s),\;g\in G^J_\Bbb R$ and $(\tau,z)\in H_1\times
\Bbb C^s.$
We let $\Gamma_M^J:=\Gamma_M\cap G_\Bbb R^J$ (cf. Definition 4.7).
Then $\Gamma_M^J$ is a discrete subgroup of $G^J_\Bbb R$ which acts
on $H_1\times\Bbb C^s$ properly discontinuously.
\definition{Definition 4.8}\ Let $k$ and $m$ be nonnegative integers.
A holomorphic function $\phi: H_1\times\Bbb C^s\rightarrow \Bbb C$
is called a Jacobi form of weight $k$ and index $m$ on $\Gamma_M^J$ if
$\phi$ satisfies the following functional equation
$$\phi|_{k,m}\gamma=\phi\quad \text{for all}\quad\gamma\in\Gamma_M^J\tag4.35$$
and $f(\tau,z)$ has a Fourier expansion
$$\phi(\tau,z)=\sum_{n\in\Z}\sum_{\ell\in\widehat{L_0}}c(n,\ell)
e^{2\pi i(n\tau+{}^t\!\ell S_0z)}\tag4.36$$
with $c(n,\ell)\neq0$ only if $2nm-S_0[\ell]\geq0.$ Here $\widehat{L_0}$
is the dual lattice of $L_0,$ i.e.,
$$\widehat{L_0}:=\{\ell\in L_0\otimes\Bbb Q\mid {}^t\!\ell S_0\alpha\in\Z\;\;
\text{for all}\;\alpha\in L_0\}.$$
A Jacobi form $\phi$ of weight $k$ and index $m$ is called a {\it cusp form}
if $c(n,\ell)\neq0$ implies $2nm-S_0[\ell]>0.$
We denote by $J_{k,m}(\Gamma^J_M)$ (resp. $J^{\text{cusp}}_{k,m}(\Gamma_M^J))$
the vector space of all Jacobi forms (resp. cusp forms)
of weight $k$ and index $m$ on $\Gamma_M^J.$
\enddefinition

\remark
{\smc Remark 4.9} \ (1) $J_{k,m}(\G_M^J)$ is finite dimensional. \newline
(2) The Fourier-Jacobi coefficients $\phi_m$ of a modular form $f$ (cf. (4.17))
are Jacobi forms of weight $k$ and index $m$ on $\G^J_M.$ (cf. (4.20) and (4.21)).
\newline
(3) If $\phi\in J_{k,m}(\G^J_M),$ the function $f_\phi(\omega,z,\tau) :=
\phi(\tau,z)e^{2\pi i m \omega}$ is a modular form with respect to the subgroup
of finite index of the integral Jacobi parabolic subgroup $P_{\Bbb Z} :=
P_{\Bbb R} \cap \G_M.$
\endremark
\par
Let $m$ be a nonnegative integer and let $G_m(L_0) := \widehat{L_0}/mL_0$ be the
{\it discriminant group} of the lattice $L_0.$ For each $ h\in G_m(L_0),$ we define
the theta function $\vartheta_{S_0,m,h}:= \vartheta_{L_0,m,h}$
$$\vartheta_{L_0,m,h}(\tau,z) := \sum_{\ell\in L_0} e^{\pi i m (S_0[\ell+\frac{h}{m}]
\tau + 2{}^t(\ell+\frac{h}{m})S_0z)}, \tag 4.37$$
where $(\tau,z) \in H_{1,s}:= H_1\times \Bbb C^s.$ Any Jacobi form $\phi\in
J_{k,m}(\G^J_M)$ can be written as
$$\phi(\tau,z)= \sum_{h\in G_m(L_0)} \,\phi_h(\tau)
\vartheta_{L_0,m,h}(\tau, z)\tag 4.38 $$
with
$$\phi_h(\tau):= \sum_{r\geq 0} c((2r+qS_0[h])(2qm)^{-1},h)e^{2\pi i \frac
{r\tau}{qm}},$$
where each $r\geq 0$ satisfies the condition
$2r\equiv -q S_0[h](\text{mod}\, 2qm),\;c(n,\ell)$ denotes the Fourier coefficients
of $\phi(\tau,z)$ and $q$ is the level of the quadratic form $S_0.$ We can
rewrite (4.38) as follows :
$$\phi(\tau,z) = {}^t\Phi(\tau)\cdot \Theta_{L_0,m}(\tau,z),\tag4.39$$
where
$$\Phi(\tau) :=\left( \phi_h(\tau)\right)_{h\in G_m(L_0)}\quad \text{and}
\quad\Theta_{L_0,m}:=\left(\vartheta_{L_0,m,h}\right)_{h\in G_m(L_0)}.
\tag4.40$$
\indent
Then we can show that for any $\gamma=\pmatrix a& b\\c&d\endpmatrix\in
SL(2,\Bbb Z),$ the theta function $\Theta_{L_0,m}$ satisfies the following
transformation formula
$$\Theta_{L_0,m}\left(\frac{a\tau +b}{c\tau +d}, \frac{z}{c\tau+d}\right)
=e^{\pi i m \frac{cS_0[z]}{c\tau +d}}\cdot (c\tau +d)^{\frac s 2}\cdot \chi(M)
\Theta_{L_0,m}(\tau,z),\tag 4.41$$
where $\chi(M)$ is a certain unitary matrix of degree $|G_m(L_0)|$
(cf. [G2], p.9 and [O], p.105).
And $\Phi(\tau)$ satisfies the following functional equations :
$$\Phi(\tau+1) = e^{-\pi i \frac{S_0[h]}{m}}\Phi(\tau),\quad
\Phi(- 1/\tau)= \tau^{k-\frac{s}{2}}\overline{U(J)}\Phi(\tau),
\tag4.42$$
where $J:=\pmatrix 0& 1 \\ -1& 0\endpmatrix$ and
$$U(J) := \left(\text{det}S_0\right)^{-\frac 12}
\left(\frac {i} {m} \right)^{\frac
s 2}\left(e^{-2 \pi i \frac{{}^t g S_0 h}{m}}\right)_{g,h\in G_m(L_0)}.\tag4.43$$
We note that the finite group $G_m(L_0)$ may be regarded as the quadratic space
equipped with the quadratic form $q_{m,L_0}$ defined by
$$q_{m,L_0}(h+mL_0):= (h+mL_0,h+mL_0)\in (h,h)_0+2\Bbb Z.\tag4.44$$
>From (4.41) and (4.42) it follows that $\Phi(\tau)$ is a vector-valued modular
form of a half-integral weight and that the vector space of Jacobi forms of
index $m$ depends only on the quadratic space $(G_m(L_0),q_{m,L_0}).$
\proclaim {\bf Lemma 4.10} \ Let $M_1$ and $M_2$ be two even integral lattices
of dimension $s_1$ ans $s_2.$ We assume that the quadratic spaces
$(G_{m_1}(M_1),q_{m_1,M_1})$ and $(G_{m_2}(M_2),q_{m_2,M_2})$ are isomorphic.
Then we have the isomorphism
$$J_{k,m_1}(\G^J_{M_1}) \cong J_{k+\frac{s_2 -s_1}{2},m_2} (\G^J_{M_2}).$$
\endproclaim
\demo{\it Proof} The proof is done if the map
$${}^t\Phi(\tau)\cdot \Theta_{M_1,m_1}(\tau,z) \longmapsto {}^t\Phi(\tau)
\cdot \Theta_{M_2,m_2}(\tau,z) \tag4.45$$
is an isomorphism of $J_{k,m_1}(\G^J_{M_1})$ onto $J_{k+\frac{s_2 - s_1}{2},m_2}
(\G^J_{M_2}).$ The isomorphism can be proved using (4.41) and $s_1 \equiv s_2\;
(\text{mod}\,8).$\hfill $\square$
\enddemo
\par
Now we discuss the concept of singular modular forms and singular Jacobi forms.
\demo{\bf Definition 4.11} \ A modular form $f$ with respect to $\G_M$ (or
a Jacobi form $\phi$ of index with respect to $\G^J_M$) is said to be
{\it singular} if its Fourier coefficients satisfy the following condition that
$$a(n,\ell,m) \ne 0\;(\text{or}\; c(n,\ell)\ne 0 )\;\text{implies}\;2nm-S_0[\ell]=0,$$
where $a(n,\ell,m)$ and $c(n,\ell)$ denote the Fourier coefficients of $f$ and $\phi$
in their Fourier expansions respectively.\enddemo
We consider the {\it differential operators} $D$ and $\widehat{D}$ defined by
$$D:= \frac{\partial^2}{\partial \omega \partial \tau} - \frac 12 S_0[\frac
{\partial}{\partial z}] \tag 4.46$$
and
$$\widehat{D} := \frac{\partial}{\partial z} -\frac{1}{4\pi i m} S_0
[\frac{\partial}{\partial z}]. \tag 4.47$$
Then it is easy to see that if $f$ if a singular modular form and if $\phi$ is
a Jacobi form of index $m,$ then $Df=0$ and $\widehat{D}\phi=0.$ We can also show
easily that any Jacobi form with respect to $\G_M^J$ has its weight $s/2$
and that any Jacobi form of weight $s/2$ with respect to $\G_M^J$ is
singular. From this fact, we see that a weight of a singular modular form with
respect to $\G_M$ is either $0$ or $s/2.$
\par For each positive integer $m,$ we let $M(m) := \Pi_{1,1} \oplus L_0
\oplus \Pi_{1,1}$ be the lattice with its associated symmetric matrix given
by
$$ S (m) := \pmatrix 0&0  & I_2 \\
                         0&mS_0 & 0 \\
             I_2&0& 0 \endpmatrix, \quad
             I_2 :=\pmatrix 0&-1\\ -1& 0\endpmatrix.\tag4.48$$
We let $\G^J_m := \G^J_{M(m)}$ be the integral Jacobi group of the lattice
$M(m).$ It follows immediately from the definitions that if
$\phi\in J_{k,m}(\G^J_M),$ then $ \phi\in J_{k,1}(\G^J_m).$
The existence of a nonconstant singular Jacobi form of index 1 with respect to
$\G_M^J$ guarantees the unimodularity of the lattice $M.$ Precisely, we have
\proclaim
{{\bf Proposition 4.12} ([G2], Lemma 4.5)} \ Let $M$ be a {\it maximal} even
integral lattice. This means that $M$ is not a sublattice of any even integral
lattice. Then a nonconstant singular Jacobi form of index 1 with respect to
$\G_M^J$ exists if and only if the lattice $M$ is unimodular.
\endproclaim
\demo
{\it Proof} \ The proof can be found in [G2], p.21. But we write his proof
here. Let $\phi(\tau,z)$ be a nonconstant singular Jacobi form of index 1 with
respect to $\G_M^J.$ According to (4.39), $\phi(\tau,z) =\Phi(\tau)\cdot
\Theta_{L_0,1}(\tau,z_0).$ Thus the components $\phi_h(\tau)\;(h\in G_1(L_0))$
of $\Phi$ are constants because their weights are all zero. By (4.42), we have
$$\phi_h(\tau+1) =e^{-\pi i S_0 [h]} \phi_h(\tau), \quad h\in G_1(L_0).$$
Therefore the components $\phi_h$ are not zero only for the isotrophic vectors
$h$ in the group $G_1(L_0)= \widehat{L_0}/{L_0}.$ Since $M= \Pi_{1,1}\oplus
L_0 \oplus \Pi_{1,1}$ is maximal, there exists only the trivial isotropic
element $h=0$ in $G_1(L_0).$ Again by (4.42), we obtain that $|G_1(L_0)| =1$
and so $L_0$ is unimodular. Hence the lattice $M$ is unimodular.\hfill
$\square$
\enddemo
\remark
{\bf Example 4.13} \ We assume that $M$ is a {\it unimodular} even integral
lattice of signature $(s+2,2).$ Then the theta series
$$\vartheta(\tau,z):= \sum_{\lambda\in L}
e^{\pi i(S_0[\lambda]\tau +2^t\!\lambda
S_0 z)},\quad (\tau,z)\in H_{1,s}\tag 4.49 $$
is a singular Jacobi form of weight $s/2$ and index 1 with respect to
$\G^J_M.$ The arithmetic lifting $f_{\vartheta}$ of $\vartheta(\tau,z)$ defined
by
$$f_\vartheta(\omega,z,\tau) := \frac{(s/2 -1)! \xi(s/2)}{(2 \pi i )^{s/2}}
+\sum\Sb n,m \ge 0, \lambda \in L \\ 2nm=S_0[\lambda] \\
(n,m) \ne (0,0)\endSb \sigma_{s/2 - 1}(n,m ;\lambda ) \,e^{2 \pi i (n\tau+
{}^t\lambda S_0 z + m \omega)}\tag 4.50$$
is a singular modular form of weight $s/2$ with respect to $\G_M^J,$
where $\sigma_{s/2-1}(n,m;\lambda)$ denotes the sum of $(s/2 -1)$-powers of all
common divisors of the numbers $n,m$ and the vector $\lambda \in L.$ For more
detail, we refer to Theorem 3.1 and Example 4.4 in [G2].
\endremark
\par
\bigpagebreak

As we have seen so far, automorphic forms on the real symplectic group and
those on the real orthogonal group have different geometric objects, different
automorphic factors (cf. (4.12)), and somewhat different properties. For
instance, in case of the orthogonal group $O_{s+2,2}(\Bbb R),$ there is a gap
between 0 and $s/2$ such that there exist no modular forms and no Jacobi
forms with weights in this gap. By the way, this phenomenon does not happen
for automorphic forms and Jacobi forms in the case of the symplectic group
$Sp(g,\Bbb R)$ because all integers less than half the largest singular weight
are also singular weights. For more detail, we refer to [F] for singular modular
forms and to [Y3] for singular Jacobi forms. In both cases the number of
singular weights is equal to the real rank of the corresponding Lie group.
Nonetheless the properties of Jacobi forms for the orthogonal group are
similar to those of Jacobi forms for the symplectic group. For example, the
Fourier coefficient $c(n,\ell)$ of a Jacobi form of weight $k$ and index $m$ for
$O_{s+2,2}(\Bbb R)$ depends only on the number $2mn - S_0[\ell]$ (which is the
norm of the vector $(n,\ell,m)$ in the lattice $\widehat{L_1}$) and the equivalence
class of $\ell$ in the discriminant group $G_m(L_0)=\widehat{L_0}/L_0$ (cf. compare
Theorem 2.2 in [E-Z] with our case). We observe that the automorphic factors
for the Jacobi groups for both cases are quite similar (cf. see (4.2) and
(4.31)-(4.33)). The expression of Jacobi forms in terms of (4.38) or (4.39)
are similar to that of Jacobi forms for the symplectic group (cf. [E-Z], [Y1],
and [Zi]).
\remark
{\smc Remark 4.14} \ In [Bo7], R. Borcherds investigates automorphic forms
and Jacobi forms for $O_{s+2,2}(\Bbb R)$ which are either nearly holomorphic
or meromorphic.
\endremark
Meromorphic functions with all poles at cusps are called {\it nearly holomorphic
ones.}
\par\ex
{\smc Borcherds' construction of Jacobi forms} : Let $K$ be a positive definite
integral lattice of dimension $s.$ A function
$c:K\rightarrow \Bbb Z^+ \cup \{0\}$ is said to be a {\it vector system} if
it satisfies the following three properties (1)--(3) :
\roster
\item"(1)" \ The set $\{v \in K | c(v)\ne 0 \}$ if finite.
\item"(2)" \ $c(v) = c(-v) $ for all $v\in K.$
\item"(3)" \ The function taking $\lambda$ to
$\sum_{v\in K} c(v) (\lambda,v)^2$ is constant on the sphere of norm 1 vectors
$\lambda \in K\otimes \Bbb R.$
\endroster
\indent
We will write $V$ for the {\it multiset} of vectors in a vector system and so
we think of $V$ as containing $c(v)$ copies of each vector $v\in K.$ And we
write $\sum_{v\in V} f(v)$ instead of $\sum_{v\in K} c(v)f(v).$ The vector system
is said to be {\it trivial} if it only contains vectors of zero norm.
\par
The hyperplanes orthogonal to the vectors of a vector system $V$ divides
$K_{\Bbb R}:= K\otimes \Bbb R$ into cones which we call the {\it Weyl chambers}
of $V.$ We note that unlike the case of root systems,
the Weyl chamber of $V$
need not be all the same type. If we choose a fixed Weyl chamber $W,$ then
we can define the {\it positive} and {\it negative} vectors of $V$ by saying
that $v$ is positive or negative, denoted by $v>0$ or $v<0$ if $(v,\lambda)>0$
or $(v,\lambda)<0$ for some vector $\lambda$ in the interior $W^0$ of $W.$
It is easy to check that the concept of positivity and negativity does not
depend on the choice of a vector $\lambda$ in $W^0.$ Obviously every nonzero
vector of the vector system $V$ is either positive or negative.
\par
We define the {\it Weyl vector} $\rho =\rho_W$ of $W$ by
$$\rho:= \frac 1 2 \sum\Sb v\in V \\ v>0\endSb v.$$
We define $d$ to be the number of vectors in $V$ and define $k:=\frac{c(0)}{2}.$
The rational number $k$ is called the {\it weight} of $V.$ We define the
{\it index} $m$ of $V$ by
$$m:= (2\, \text{dim} \,K)^{-1} \sum_{v\in V}(v,v)$$
\par
We can show that the index $m$ of $V$ is a nonnegative integer.
If $V$ is a
vector system in $K,$ we define the {\it (untwisted) affine vector system} of
$V$ to be the multiset of vectors $(v,n) \in K\oplus \Bbb Z$ with
$v\in V.$
We say that $(v,n)$ is {\it positive} if either $n>0$ or $n=0,\;v>0.$ It can be
seen that the Weyl vectors for different Weyl chambers differ by elements of
$K.$ \par
Borcherds (cf. [Bo7], p.183) define the function $\psi(\tau,z)$ on
$H_1\times K_{\Bbb C}$ with $K_{\Bbb C} = K\otimes \Bbb C \cong \Bbb C^s$ by
$$\psi(\tau,z) := q^{\frac{d}{24}}\zeta^{-\rho}\prod_{(v,n)>0}(1-q^n\zeta^v),\;\;(\tau,z)
\in H_1 \times K_{\Bbb C},\tag4.51$$
where $(v,n)$ runs over the set of all positive vectors in the affine vector
system of $V,\;q^a := e^{2 \pi i a \tau}$ and $\zeta^v=e^{2\pi i (z,v)}.$ Then
$\psi(\tau,z)$ is a {\it nearly holomorphic Jacobi form} of weight $k$ and
index $m.$ Thus $\psi$ can be written as a finite sum of theta functions times
nearly holomorphic modular forms. In fact, $\psi$ satisfies the following
transformation formulas:
$$\align
\psi(\tau+1,z)& =e^{\frac{\pi i d}{12}} \psi(\tau,z), \\
\psi(-1/\tau,z/\tau)& =(-i)^{d/2-k}(\tau/i)^k e^{\pi i m(z,z)/\tau}
\psi(\tau,z),\\
\psi(\tau,z+\mu)&=(-1)^{2(\rho,\mu)}\psi(\tau,z),\\
\psi(\tau, z+\lambda\tau) &=(-1)^{2(\rho,\lambda)}q^{-m(\lambda,\lambda)/2}
\zeta^{-m\lambda}\psi(\tau,z)
\endalign$$
for all $\lambda,\mu \in \widehat{K}$ (the dual of $K$).

\head {\bf 5.\ Infinite\ Products\ and\ Modular\ Forms} 
\endhead               
\indent
In [Bo7], R. Borcherds constructed automorphic forms on $O_{s+2,2}(\Bbb R)^0$
which are modular products and using the theory of these automorphic theory
expressed some meromorphic modular forms for $SL(2,{\Bbb Z})$ with
certain conditions as
infinite products. Roughly speaking, a modular product means an infinite
product whose exponents are the coefficients of a certain nearly holomorphic
modular form. For instance, he wrote modular forms as the modular invariant
$j$ and the Eisenstein series $E_4$ and $E_6$ as infinite products.
These results tell us implicitly that the denominator function of a
generalized Kac-Moody algebra is sometimes an automorphic form on
$O_{s+2,2}({\Bbb R})^0$ which is a modular product. In this section we
discuss Borcherds' results just mentioned in some detail.
\par\medpagebreak
\indent
We shall start by giving some well-known classical product identities.
First we give some of the product identities of L. Euler\,(1707-83) which
are
$$\sum_{n\geq 0}{{(-1)^nq^{n(n+1)/2}z^n}\over {(1-q)(1-q^2)\cdots (1-q^n)}}
=\prod_{n>0}(1-q^nz),\tag 5.1$$
$$\sum_{n\geq 0}{{z^n}\over {(1-q)(1-q^2)\cdots (1-q^n)}}=
\prod_{n>0}(1-q^nz)^{-1},\tag 5.2$$
$$\sum_{n\in {\Bbb Z}}(-1)^n q^{3(n+1/6)^2/2}=q^{1/{24}}\prod_{n>0}
(1-q^n).\tag 5.3$$
A similiar product identity due to C. F. Gauss\,(1777-1855) is
$$\sum_{n\in {\Bbb Z}}q^{n^2}=(1+q^2)(1-q^2)(1+q^3)^2(1-q^4)\cdots .
\tag 5.4$$
Both of (5.3) and (5.4) are special cases of the so-called Jacobi's triple
product identity\,[\,C.\,G.\,J. Jacobi\,(1804-51)\,]
$$\sum_{n\in {\Bbb Z}}(-1)^nq^{n^2}z^n = \prod_{n>0}(1-q^{2n})
(1-q^{2n-1}z)(1-q^{2n-1}z^{-1})\tag 5.5$$
if we choose $z$ to be some fixed power of $q$. In fact, if you replace
$q$ and $z$ in (5.5) by $q^{3/2}$ and $q^{1/2}$ respectively, you obtain
the identity (5.3), and if you replace $z$ in (5.5) by $-1$, you get the identity
(5.4).\par\smallpagebreak
\indent
The quintuple product identity derived by G.\,N. Watson\,(1886-1965) is
$$\multline
\sum_{n\in {\Bbb Z}}q^{(3n^2+n)/2}(z^{3n}-z^{-3n-1}) \\
=\prod_{n>0}
(1-q^n)(1-q^nz)(1-q^{n-1}z^{-1})(1-q^{2n-1}z^2)(1-q^{2n-1}z^{-2}).
\endmultline \tag 5.6$$
Historically speaking, in 1929 Watson\,(\,cf.\,[W1]\,) derived the identity
(5.6) in the course of proving some of Ramanujan's theorems on continued
fractions. In 1938, Watson\,(\,cf.\,[W2]\,) proved the following identity\,:
$$\multline
\sum_{n\in {\Bbb Z}}q^{n(3n+2)}(z^{-3n}-z^{3n+2})\\
=\prod_{n>0}
(1-q^{2n})(1-q^{2n-2}z^2)(1-q^{2n}z^{-2})(1+q^{2n-1}z)^{-1}
(1+q^{2n-1}z^{-1}).\endmultline \tag 5.7$$
\indent
Subbarao and Vidyasagar\,(\,cf.\,[S-V]\,) showed that the identities (5.6)
and (5.7) are equivalent. The two identities (5.1) and (5.2) of Euler are
easily established\,(\,cf.[Be],\,p.\,49\,). G.\,E. Andrews showed that
the Jacobi's triple product identity (5.5) can be obtained easily from
the identities (5.1) and (5.2) in his short paper\,[A]. Carlitz and
Subbarao\,(\,cf.\,[C-S]\,) gave a simple proof of
the quintuple product identity (5.7).
\par\medpagebreak
\indent
\redefine\a{\alpha}
The following denominator formula for a finite dimensional simple Lie
algebra $\frak g$
$$e^{\rho}\sum_{w\in W}\text{det}\,(w)\,e^{-w(\rho)}=
\prod_{\a > 0}(1-e^{\a})\tag 5.8$$
is due to Hermann Weyl\,(\,1885-1955\,), where $W$ is the Weyl group of
${\frak g},\ \rho$ is the Weyl vector and the product runs over the set of
all positive roots. Macdonald\,(\,cf.\ [Mac]\,) observed that the Weyl
denominator formula is just a statement about finite root systems, and then
generalized this formula to {\it affine} root systems producing the
so-called {\it Macdonald\ identities.} He noticed that the Jacobi's triple
product identity is just the Macdonald identity for the simplest affine
root system. Kac observed that the Macdonald identities are just the
denominator formulas for the Kac-Moody Lie algebras in the early 1970s.
Thereafter he obtained the so-called {\it Weyl-Kac\ character\ formulas}
for representations of the affine Kac-Moody algebras generalizing the
Weyl character formula\,(\,see\,(2.5)-(2.7) and [K],\,p.\,173\,). The
Weyl-Kac character formula for the affine Kac-Moody algebra is given as
follows\,:
$$e^{\rho}\sum_{w\in W}\text{det}(w)\,e^{-w(\rho)}=\prod_{\a>0}
(1-e^{-\a})^{\text{mult}(\a)},\tag 5.9$$
where $\text{mult}\,(\a)$ is the multiplicity of the root $\a$. For more
detail we refer to (2.6) and [K]. For instance, the Jacobi's triple
product identity is just the Weyl-Kac character formula for the affine
Kac-Moody algebra $SL_2({\Bbb R}[z,z^{-1}])$ and the Weyl-Kac character
formulas for the affine Kac-Moody algebras $SL_n({\Bbb R}[z,z^{-1}])$ are
just the Macdonald identities. It seems that the Weyl-Kac character
formula is true for non-affine Kac-Moody algebras. Borcherds obtained the
so-called {\it Weyl-Kac-Borcherds\ character\ formula} for a generalized
Kac-Moody algebra\,(\,cf.\,(2.13)-(2.14)\,). The Weyl-Kac character formula
is proved by the Euler-Poincar{\' e} principle applied to the cohomology
of the Lie subalgebra $E$ of $\frak g$ associated to the positive roots
of the Kac-Moody algebra $\frak g.$
\par\medpagebreak
\indent
It seems to the author that Borcherds was the first one that
discovered that the denominator functions
of the generalized Kac-Moody algebras which could be written as infinite
products are often automorphic forms on the orthogonal group
$O_{s+2,2}({\Bbb R})^0.$  Moreover he gave a method of constructing
automorphic forms on $O_{s+2,2}({\Bbb R})^0$ through modular forms of
weight $-s/2$ with integer coefficients and obtained the connection
between the
Kohnen's ``plus" space of weight $1/2$ and the space of modular forms on
$\Gamma_1$ satisfying some conditions.
\par\medpagebreak
\indent
Now we are in a position to describe his works on infinite products related
to automorphic forms on the orthogonal group $O_{s+2,2}({\Bbb R}).$
\par\medpagebreak
\indent
We let $L$ be the unimodular even integral Lorentzian lattice $\Pi_{s+1,1}$
of dimension $s+2$
and let $M:=L\oplus \Pi_{1,1},$ where $\Pi_{1,1}$ is the unique
$2$-dimensional unimodular even integral Lorentzian lattice with its inner
product matrix $\pmatrix 0 & -1\\ -1 & 0 \endpmatrix$.
We choose a negative norm vector $\a$ in $L_{\Bbb R}:=L\otimes {\Bbb R}.$
We say that a vector $v$ in $L_{\Bbb R}$ is {\it positive}, denoted by
$v > 0$ if $(v,\a)>0.$
\par \medpagebreak
\noindent
{\bf Theorem\ 5.1\,(\,Borcherds\,[Bo7],\,Theorem\ 10.1\,).} Let $f(\tau)=
\sum_n c(n)q^n$ be a nearly holomorphic modular form of weight $-s/2$ for
$\Gamma_1$ with integer coefficients, with $24|c(0)$ if $s=0.$ Then there is a
unique vector $\delta\in L$ such that
$$\Phi(v):=e^{-2\pi i(\delta,v)}\prod_{r>0,\,r\in L}
\left(1-e^{-2\pi i(r,v)}\right)^{c(-r^2/2)},\ \ \ v\in \Omega\tag 5.10$$
is a meromorphic automorphic form of weight $c(0)/2$ for $O_M({\Bbb Z})^0
\cong O_{s+2,2}({\Bbb Z})^0,$ where $r^2:=(r,r)$ and
$$\Omega:=\left\{\,z\in M\otimes {\Bbb C}\,\vert \ (z,z)=0,\
(z,{\bar z})>0\,\right\}.$$
\noindent
{\smc Remark\ 5.2.} Borcherds showed that all the zeros and poles of
$\Phi$ lie on the rational quadratic divisors and computed the
multiplicites of the zeros of $\Phi.$ Roughly speaking a rational quadratic
divisor means the zero set of $a(y,y)+(b,y)+c=0$ with $a,c\in {\Bbb Z}$
and $b\in L.$\par\medpagebreak
\noindent
\redefine\lrt{\longrightarrow}
{\bf Definition\ 5.3.} We define the function $H:{\Bbb Z}_+\lrt
{\Bbb Q}$ by
$$H(n):=\cases \text{the Hurwitz class number of the discriminant $-n$ if
$n>0$}\,;\\
-1/{12} \ \text{if $n=0$.}\endcases$$
We note that
$${\tilde H}(q):=\sum_{n\geq 0}H(n)q^n=-1/12+q^3/3+q^4/2+q^7+q^8+q^{11}+
(4/3)q^{12}+\cdots.$$
\indent
Now we state a very interesting result.

\redefine\BZ{\Bbb Z}
\redefine\BR{\Bbb R}
\redefine\BQ{\Bbb Q}
\redefine\BC{\Bbb C}
\par\medpagebreak
\noindent
{\bf Theorem\ 5.4\,(\,Borcherds\,[Bo7],\,Theorem\,14.1).}
Let ${\Cal A}$ be the additive group consisting of
nearly holomorphic modular forms of weight $1/2$ for $\Gamma_0(4)$ whose
coefficients are integers and satisfy the Kohnen's ``plus space" condition.
We also let ${\Cal B}$ be the multiplicative group consisting of meromorphic
modular forms for some characters of $\Gamma_1$ of integral weight with
leading coefficient $1$ whose coefficients are integers and all of whose
zeros and poles are either cusps or imaginary quadratic irrationals. To
each $f(\tau)=\sum_n c(n)q^n$ in ${\Cal A}$ we associate the function
$\Psi_f:H_1\lrt \BC$ defined by
$$\Psi_f(\tau):=q^{-h}\prod_{n>0}(1-q^n)^{c(n^2)},\tag 5.11$$
where $h$ is the constant term of $f(\tau){\tilde H}(q).$ Then we have
the following\,:
\par\medpagebreak
\noindent
(a) For each $f\in {\Cal A},\ \Psi_f$ is an element of ${\Cal B}$ whose
weight is $c(0)$\,;\par\smallpagebreak
\noindent
(b) the map $\Psi:{\Cal A}\lrt {\Cal B}$ given by $\Psi(f):=\Psi_f$ for
$f\in {\Cal A}$ is a group isomorphism of ${\Cal A}$ onto ${\Cal B}$\,;
\par\smallpagebreak
\noindent
(c) the multiplicity of the zero of $\Psi$ at a quadratic irrational
$\tau$ of discriminant $D<0$ is $\sum_{d>0}c(Dd^2).$
\par\medpagebreak
\noindent
{\smc Remark\ 5.5.} The product formula for the classical modular
polynomial\,(\,for discriminant $D<0$ whose degree is $H(-D)$\,)
$$\prod_{[\sigma]}(j(\tau)-j(\sigma))=q^{-H(-D)}\prod_{n>0}
(1-q^n)^{c(n^2)} \tag 5.12$$
holds, where $\sigma$ runs over a complete set of representatives modulo
$\Gamma_1$ for the imaginary quadratic irrationals which are roots of an
equation of the form $a\sigma^2+b\sigma+c=0\,(\,a,b,c\in \BZ\,)$ of the
discriminant $b^2-4ac=D<0$\,(\,except that $\sigma$ is a conjugate of one
of the elliptic fixed points $i$ or $(1+i\sqrt{3})/2$ we have to replace
the corresponding factor $j(\tau)-1728$ or $j(\tau)$ by
$(j(\tau)-1728)^{1/2}$ or $j(\tau)^{1/3})$ and the exponents $c(n^2)$ are
the coefficients of the uniquely determined
nearly holomorphic modular form in ${\Cal A}.$
It is easy to check that the classical modular polynomial on the left hand
side of (5.12) is contained in ${\Cal B}$ and that its corresponding
element in ${\Cal A}$ is of the form $q^D+O(q).$
\par\medpagebreak
\noindent
{\bf Examples\ 5.6.} (1) Let $f(\tau):=12\theta(\tau)=12\sum_{n\in\BZ}
q^{n^2}$. It is easy to check that $f(\tau)$ is an element of ${\Cal A}$ and
that $\Psi_f(\tau)=q\prod_{n>0}(1-q^n)^{24}$ is a cusp form for
$\Gamma_1$ of weight $12$ known as the discriminant function.
\par\medpagebreak\noindent
(2) We put
$$F(\tau):=\sum_{n>0,\,n:\text{odd}}
\sigma_1(n)q^n=q+4q^3+6q^5+8q^7+13q^9+\cdots.$$
We let
$$f(\tau)=3F(\tau)\theta(\tau)(\theta(\tau)^4-2F(\tau))\,(\theta(\tau)^4 -
16F(\tau))\,E_6(4\tau)/\Delta(4\tau)+168\theta(\tau),$$
where $\theta(\tau)=\sum_{n\in \BZ}q^{n^2},\ \Delta(\tau)$ and $E_4(\tau)$
denote the discriminant function and the Eisenstein series of weight $4$
respectively.\,(\,see Appendix A\,). It is easy to check that $f(\tau)$ is an element of
${\Cal A}$ and that $\Psi_f(\tau)=j(\tau)$ is the modular invariant. We also
check that $\Psi_f(\tau)=j(\tau)$ has order $3$ at the zero
${1+i\sqrt{3}}\over 2$ whose discriminant is $-3.$ Hence we obtain the
modular product
$$j(\tau)=q^{-1}(1-q)^{-744}(1-q^2)^{80256}(1-q^3)^{-12288744}\cdots.$$
\noindent
(3) The Eisenstein series $E_4$ and $E_6$ are elements of ${\Cal B}$. The
elements of ${\Cal A}$ corresponding to $E_4$ and $E_6$ are given by
$$f_4(\tau)=q^{-3}+4-240q+26760q^4-85995q^5+1707264q^8-4096240q^9+\cdots$$
and
$$\multline
f_6(\tau)=q^{-4}+6+504q+143388q^4\\+565760q^5+184373000q^8+
51180024q^9+O(q^{12})\qquad
\endmultline$$
respectively. Use the fact $E_4^3=j\cdot \Delta$ for $f_4.$ The function
$f_6(\tau)$ can be obtained from the theory of a generalized Kac-Moody algebra
of rank $1$ whose simple roots are all multiples of some root $\a$ of norm
$-2$ and the simple roots are $n\a$ with $n\in \BZ^+$ and multiplicity
$504\sigma_3(n).$ Precisely,
$$\multline
f_6(\tau)=(j(4\tau)-876)\theta(\tau)\\
-2F(\tau)\theta(\tau)
(\theta(\tau)^4-2F(\tau))(\theta(\tau)^4-16F(\tau))E_6(4\tau)/\Delta(4\tau),
\qquad
\endmultline$$
where $\theta(\tau)$ and $F(\tau)$ are defined in (2). Since $E_8=E_4^2,\
E_{10}=E_4\!\cdot \!E_6$
and $E_{14}=E_4^2E_6,$ their corresponding elements in
${\Cal A}$ are given by $2f_4,\ f_4+f_6$ and $2f_4+f_6$ respectively. The
remaining Eisenstein series\,(\,$k\ne 4,6,8,10,14,\ k:\text{even},\ k\geq 4\,)$
are not elements of ${\Cal B}$ and hence they cannot be written as modular
products. For instance, the modular products for $E_4,\ E_6,
\ E_8,\ E_{10}$ and
$E_{14}$ are given by
$$E_4(\tau)=(1-q)^{-240}(1-q^2)^{26760}(1-q^3)^{-4096240}\cdots,$$
$$E_6(\tau)=(1-q)^{504}(1-q^2)^{143388}(1-q^3)^{51180024}\cdots,$$
$$E_8(\tau)=(1-q)^{-480}(1-q^2)^{53520}(1-q^3)^{-8192480}\cdots,$$
$$E_{10}(\tau)=(1-q)^{264}(1-q^2)^{170148}(1-q^3)^{47083784}\cdots$$
and
$$E_{14}(\tau)=(1-q)^{24}(1-q^2)^{196908}(1-q^3)^{42987544}\cdots.$$
\noindent
(4) Using the above theorem, we can show that there exist precisely $14$
modular forms of weight 12 on $\Gamma_1$ which are contained in ${\Cal B}$.
Indeed, if
$$\align   \Xi:&=\{\,n\in \BZ\,|\ j(\tau)-n\in {\Cal B}\ \}\\
&=\{\,j(\tau)\in \BZ\,\vert\ \tau\in H_1,\ \tau\ \text{is\ imaginary\
quadratic}\,\},\endalign$$
only the modular forms $\Delta(\tau)(j(\tau)-n)$\,(\,where\ $n\in\Xi\,)$ and
$\Delta(\tau)$ are modular forms of weight $12$ in ${\Cal B}.$
It is well
known that the elements of $\Xi$ are
$$\align
&j\left(\frac {1+i\sqrt{3}}{2}\right)=0,\ \ \ j(i)=2^6\!\cdot\! 3^3,\ \ \
j\left(\frac {1+i\sqrt{7}}{2}\right)=-3^3\!\cdot\!5^3,\ \ \
j(i\sqrt{2})=2^6\!\cdot\!5^3,\\
&j\left(\frac {1+i\sqrt{11}}{2}\right)=-2^{15},\ \ \ \
j(i\sqrt{3})=2^4\!\cdot\!3^3\!\cdot\!
5^3, \ \ \ \ j(2i)=2^3\!\cdot\!3^3\!\cdot\!11^3,\\
& j\left(\frac {1+i\sqrt{19}}{2}\right)=
-2^{15}\!\cdot\!3^3,\ \ \
j\left(\frac {1+i\sqrt{27}}{2}\right)=-2^{15}\!\cdot\!3\!\cdot\!5^3, \ \ \
j(i\sqrt{7})=
3^3\!\cdot\!5^3\!\cdot\!17^3, \\
&j\left(\frac{1+i\sqrt{43}}{2}\right)=-2^{18}\!\cdot\!
3^3\!\cdot\!5^3, \ \ \ \
j\left(\frac{1+i\sqrt{67}}{2}\right)=-2^{15}\!\cdot\! 3^3\!\cdot\!
5^3\!\cdot\! 11^3,\\
&j\left(\frac{1+i\sqrt{163}}{2}\right)=-2^{18}\!\cdot\!3^3\!\cdot\!5^3\!\cdot\!
23^3\!\cdot\!29^3.
\endalign$$

\head {\bf 6.\ Final\ Remarks} 
\endhead               
\indent
In this final section we make some brief remarks on the fake monster Lie
algebras, generalized Kac-Moody algebras of the arithmetic type, hyperbolic
reflection groups and Jacobi forms. Finally we give some open problems.
\par\bigpagebreak
\head {\bf 6.1.\ The\ Fake\ Monster\ Lie\ Algebras} 
\endhead       
\indent
\def\L{\Lambda}
First of all we collect the properties of the
{\it fake\ monster\ Lie\ algebra}
$M_{\L}.$\,(\,In [Bo5], $M_{\L}$ was called just the monster Lie algebra because
the monster Lie algebra $M$ defined in section 3 had not been discovered
at that time yet.)\par\medpagebreak
\indent
Let $\L$ be the Leech lattice of dimension $24.$
$M_{\L}$ is the generalized Kac-Moody Lie algebra with
the following properties $(M_{\L}1)-(M_{\L}10)$\,:\par\medpagebreak
\indent
{\bf $(M_{\L}1)$} The root lattice $L$ of $M_{\L}$ is $\Pi_{25,1}:=
\L\oplus\Pi_{1,1}.$\par\smallpagebreak
\def\l{\lambda}
{\bf $(M_{\L}2)$} $\rho=(0,0,1)$ is the Weyl vector of $L$ with norm
$\rho^2=0.$ The real simple roots of $M_{\L}$ are the norm $2$ vectors of
the form $(\l,1,\l^2/2-1),\,\l\in \L,$ and the imaginary simple roots
are the positive multiples of $\rho$ each with multiplicity $24.$\,(\,We
observe that if $r$ is a real simple root, then $(\rho,r)=-1\,)$
\par\smallpagebreak
{\bf $(M_{\L}3)$} A nonzero vector $r\in L=\Pi_{25,1}$ is a root if and
only if $r^2\leq 2,$ in which case it has multiplicity $p_{24}(1-r^2/2),$
where $p_{24}(1-r^2/2)$ is the number of partitions of $1-r^2/2$ into
$24$ colours.
\par\smallpagebreak
{\bf $(M_{\L}4)$} $M_{\L}$ has a $\Pi_{25,1}$-grading. The piece
$M_{\L}(r)$ of degree $r\in \Pi_{25,1}, r\ne 0$ has dimension
$p_{24}(1-r^2/2).$\par\smallpagebreak
{\bf $(M_{\L}5)$} $M_{\L}$ has an involution $\omega$ which acts as $-1$
on $\Pi_{25,1}$ and also on the piece $M_{\L}(0)$ of degree
$0\in \Pi_{25,1}.$\par\smallpagebreak
{\bf $(M_{\L}6)$} $M_{\L}$ has a contravariant bilinear form $(\,,\,)$
such that $M_{\L}(k)$ is orthogonal to $M_{\L}(l)$ with respect to
$(\,,\,)$ if $k\ne -l,\,k,l\in \Pi_{25,1}$ and such that $(\,,\,)$ is
positive definite on $M_{\L}(k)$ for all $k\in \Pi_{25,1}$ with $k\ne 0.$
\par\smallpagebreak

{\bf $(M_{\L}7)$} The denominator formula for $M_{\L}$ is given by
$$e^{-\rho}\sum_{w\in W}\sum_{n\in \BZ}\text{det}(w)\,\tau(n)\, e^{w(n\rho)}=
\prod_{r\in L^+}(1-e^r)^{p_{24}(1-r^2/2)},\tag 6.1$$
where $W$ is the Weyl group, $L^+$ is the set of all positive roots of
$M_{\L},$
and $\tau(n)$ is the Ramanujan tau function.\,(\,The discriminant
function $\Delta(\tau)$ is the generating function of $\tau(n)$.\,)
Indeed, $L^+$ is given by
$$L^+=\left\{\,v\in \Pi_{25,1}\,\vert\ v^2\leq 2,\ (v,\rho)<0\,\right\}
\cup \left\{\,n\rho\,\vert\, n\in \BZ^+\,\right\}.$$
\indent
{\bf $(M_{\L}8)$} The universal central extension
${\hat M}_{\L}$ of $M_{\L}$ is a $\Pi_{25,1}$-graded Lie algebra. If
$0\ne r\in \Pi_{25,1},$ then the piece ${\hat M}_{\L}(r)$ of
${\hat M}_{\L}$ of degree $r$ is mapped isomorphically to $M_{\L}(r).$
The piece ${\hat M}_{\L}(0)$ of degree $0$, called the Cartan subalgebra of
${\hat M}_{\L},$ can be represented naturally as the sum of a
{\it one}-{\it dimensional} space for each vector of $\L$ and a space of
dimension $24^2=576$ for each positive integer.\par\smallpagebreak
\indent
{\bf $(M_{\L}9)$} For each $r\in L^+,$ we put
$$m(r):=\sum_{\scriptstyle n>0\atop
\scriptstyle r/n\in \Pi_{25,1}}\text{mult}(r/n)\!\cdot\!n.$$
Then for each $r\in L^+,$ we have the following formula
$$(r+\rho)^2m(r)=\sum_{\scriptstyle \alpha,\beta\in L^+\atop
\scriptstyle \alpha+\beta=r}(\alpha,\beta)\,m(\alpha)\,m(\beta).\tag 6.2$$
\par\smallpagebreak
\indent
{\bf $(M_{\L}10)$} $M_{\L}$ is a Aut(${\hat{\L}}$)-module. In fact,
Aut(${\hat{\L}}$) acts naturally on the vertex algebra of ${\hat \L}$
and hence on $M_{\L}.$\par\medpagebreak
\indent
The detail for all the properties $(M_{\L}1)$-$(M_{\L}10)$ can be found in
[Bo5].\par\medpagebreak
\noindent
{\smc Remark\ 6.1.} (a) $M_{\L}$ is essentially the space of
physical vectors of
the vertex algebra of ${\hat {\Pi}}_{25,1}$, where ${\hat{\Pi}}_{25,1}$ is
the unique central extension of $\Pi_{25,1}$ by $\BZ_2.$\par
\smallpagebreak\noindent
(b) $M_{\L}$ can be constructed from the vertex algebra of $V_{\L}$ of
the central extension ${\hat{\L}}$ of $\L$ by $\BZ_2$ in the same way that
the monster Lie algebra $M$ was constructed from the monster vertex algebra
$V$ in section 3.\par\smallbreak\noindent
(c) The multiplicities $p_{24}(1+n)$ of the roots of $M_{\L}$ is given by
the Rademacher's formula
$$p_{24}(1+n)=2\pi n^{-13/2}\sum_{k>0}{{I_{13}(4\pi\sqrt{n}/k)}\over k}
\cdot \sum_{\scriptstyle 0\leq h,h'\leq k\atop
\scriptstyle hh'\equiv -1\,(\text{mod}\,k)}e^{2\pi i(nh+h')/k},$$
where $I_{13}(z):=-iJ_{13}(iz)$ is the modified Bessel function of order
$13.$ In particular, $p_{24}(1+n)$ is asymptotic to
$2^{-1/2}n^{-27/4}e^{4\pi \sqrt{n}}$ for large $n.$
\par\medpagebreak\indent
In [Bo6], Borcherds constructed a family of Lie algebras and superalgebrs,
the so-called {\it monstrous\ Lie\ superalgebra} whose denominator formulas
are twisted denominator formulas of the monster Lie algebra $M.$ For each
element $g$ in the MONSTER $G$, we define the {\it monstrous\ Lie\ algebra}
of $g$ to be the genrealized Kac-Moody superalgebra which has its root
lattice $\Pi_{1,1}$ and simple roots $(1,n)$ with multiplicity
$tr(g\vert_{V^{\sharp}_n}).$ The denominator formula for the monstrous
Lie superalgebra $M_g$ of $g$ is given by
$$\multline
T_g(p)-T_g(q)=\sum_{m}tr(g|_{V_m^{\sharp}})p^m-
\sum_n tr(g|_{V_n^{\sharp}})q^n\\
=p^{-1}\prod_{m>0,\,n\in \BZ}
(1-p^mq^n)^{\text{mult}_g(m,n)}.\endmultline\tag 6.3$$
The multiplicity $\text{mult}_g(m,n)$ of the root $(m,n)\in \Pi_{1,1}$ is
$$\text{mult}_g(m,n)=\sum_{ds|(m,n,N)}{{\mu(s)}\over {ds}}
\,tr(g^d\vert_{V_{mn}^{\sharp}}).\tag 6.4$$
where $N$ is the order of $g$. We recall that the Thompson series
$T_g(q)$ of $g$ is the normalized generator for a genus zero function
field of a discrete subgroup of $SL(2,\BR)$ containing the Hecke subgroup
$\Gamma_0(hN)$, where $h$ is a positive integer with $h\vert (24,N).$
\par\medpagebreak\indent
Furthermore Borcherds\,(\,cf.\,[Bo6],) constructed a family of superalgebras
whose denominator formulas are twisted ones of the fake monster Lie algebra
$M_{\L}$ in the same way that he constructed a family of monstrous Lie
superalgebras from the monster
Lie algebra $M.$
\par\medpagebreak\indent
Let $g$ be an element of $\text{Aut}({\hat \L})\cong 2^{24}\cdot
\text{Aut}(\L)$ of order $N$. We let
$$L:=\left\{\,\l\in \L\,\vert\ g\l=\l\,\right\}\tag 6.5 $$
be the sublattice of $\L$ fixed by $g$. Then the dual $L'$ of $L$ is
equal to the projection of $\L$ into the vector space $L_{\BR}:=
L\otimes_{\BZ}\BR$ because $\L$ is unimodular. For simplicity we assume
that any power $g^n$ of $g$ fixes all elements of ${\hat \L}$ which are
in the inverse image of $\L^{g^n},$ where $\L^{g^n}$ is the set of
elements in $\L$
fixed by $g^n.$ According to [Bo3], there exists a reflection
group $W^g$ acting on $L$ with following properties (W1)-(W2)\,:
\par\medpagebreak\indent
(W1) The positive roots of $W^g$ are the sums of the conjugates of some
positive real roots of $\Pi_{25,1}.$\par\smallpagebreak\indent
(W2) Let $\rho$ be the Weyl vector of $W^g.$ The simple roots of $W^g$ are
the sums of orbits of simple roots of $W$ that have positive norms and
they are also the roots of $W^g$ such that $(r,\rho)=-r^2/2$ with
$\rho^2=0.$
\par\medpagebreak\indent
Let ${\frak g}_g$ be a generalized Kac-Moody superalgebra with the
following simple roots\,:\par\medpagebreak
\noindent
1. $L$ is the root lattice of ${\frak g}_g.$\par\smallpagebreak
\noindent
2. The real simple roots are the simple roots of the reflection group $W^g$,
which are the roots $r$ with $(r,\rho)=-r^2/2.$\par\smallpagebreak
\noindent
3. The imaginary simple roots are $n\rho\,(\,n\in\BZ^+\,)$ with multiplicity
$\text{mult}_g(n\rho)$ given by
$$\text{mult}_g(n\rho)=\sum_{ja_k=n}b_k$$
if $g$ has a generalized cycle shape $a_1^{b_1}a_2^{b_2}\cdots.$
\par\medpagebreak\noindent
Then the denominator formula for the {\it fake\ monstrous\ superalgebra}
${\frak g}_g$ is given by
$$e^{-\rho}\sum_{w\in W^g}\text{det}(w)\,w\left(\eta_g(e^{\rho})\right)=
\prod_{r\in L^+}(1-e^r)^{\text{mult}_g(r)},\tag 6.6$$
where $\eta_g(q)$ is the function defined by
\def\ve{\varepsilon}
$$\eta_g(q):=\eta(\ve_1q)\eta(\ve_2q)\cdots \eta(\ve_{24}q)\tag 6.7$$
if $g$ has eigenvalues $\ve_1,\cdots,\ve_{24}$ on $\L_{\BR}:=\L\otimes_{\BZ}
\BR.$ It is easy to check that if $g$ has a generalized cycle shape
$a_1^{b_1}a_2^{b_2}\cdots,$ then
$$\eta_g(q)=\eta(q^{a_1})^{b_1}\eta(q^{a_2})^{b_2}\cdots.$$
\noindent
{\bf Example\ 6.2.} Let $p=2,3,5,7,11,23$ be six prime numbers such that
$p+1$ divides $24$. We let $g$ be an element of $\text{Aut}({\hat \L})$ of
order $p$ corresponding to an element of $M_{24}\subset \text{Aut}(\L)$
of cycle shape $1^{24/(p+1)}p^{24/(p+1)},$ where $M_{24}$ is the Mathieu
group. Then the denominator formula for the fake monstrous superalgebra\,(\,in
fact, a Lie algebra\,) ${\frak g}_p:={\frak g}_g$ is given by
$$\multline e^{-\rho}\sum_{w\in W^g}\text{det}(w)\,w\left(e^{\rho}\prod_{n>0}
(1-e^{n\rho})^{24/(p+1)}(1-e^{pn\rho})^{24/(p+1)}\right) \\
=\prod_{r\in L^+}(1-e^r)^{p_g(1-r^2/2)}\prod_{r\in pL^+}
(1-e^r)^{p_g(1-r^2/2p)},\endmultline \tag 6.8$$
where $L^+$ denotes the set of all positive roots of ${\frak g}_p$ and
$p_g(1+n)$ is defined by
$$\sum_{n>0}p_g(1+n)q^n=1/\eta_g(q).\tag 6.9$$
${\frak g}_2,\,{\frak g}_3,\,{\frak g}_5,\,{\frak g}_7$ and ${\frak g}_{11}$
are called the {\it fake\ baby} monster Lie algebra, the {\it fake\ Fischer}
monster Lie algebra, the {\it fake Harada\,-\,Norton} monster Lie algebra,
the {\it fake\ Held} monster Lie algebra and the {\it fake\ Mathieu}
monster Lie algebra respectively. We observe that the dimension
of ${\frak g}_p\,(\,p=2,3,5,7,11,23\,)$ are $18,14,10,8,6,1$ respectively.
\par\medpagebreak \noindent
{\bf Example\ 6.3.} Let ${\frak g}_{fC}$ be the {\it fake\ Conway\ Lie\
superalgebra} of rank $10.\ {\frak g}_{fC}$ is the fake monstrous
Lie superalgebra associated with an element $g\in \text{Aut}({\hat\L})$
of order $2$ such that the descent $g_0$ of $g$ to $\text{Aut}(\L)$ is
of order $2$ and the lattice $\L^{g_0}$ of $\L$ fixed by $g_0$ is
isomorphic to the lattice $E_8$ with all norms doubled. The lattice $L$ of
${\frak g}_{fC}$ is the nonintegral lattice of determinant $1/4$ all whose
vectors have integral norm which is the dual lattice of the sublattice of
even vectors of $I_{9,1}.$ Here $I_{9,1}:=\left\{\,(v,m,n)\,\vert\
v\in E_8,\ m,n\in \BZ,\ m+n\ \text{is\ even}\,\right\}$ is the lattice of
dimension $10.$ Let $W$ be the Weyl group of ${\frak g}_{fC}.$
In other words, $W$ is the subgroup of $\text{Aut}(L)$ generated by
the reflection of norm $1$ vectors. The simple roots of $W$ are the norm 1
vectors with $(r,\rho)=-1/2.$ The simple roots of ${\frak g}_{fC}$ are
the simple roots of $W$ together with the positive multiple
$n\rho\,(\,n\in \BZ^+\,)$ of the Weyl vector $\rho=(0,0,1)$ each with
multiplicity $8(-1)^n.$ Here the multiplicity $-k<0$ means a superroot of
multiplicity $k$, so that the odd multiples of $\rho$ are superroots. The
multiplicity $\text{mult}(r)$ of the root $r=(v,m,n)\in L$ is given by
$$\text{mult}(r)=(-1)^{(m-1)(n-1)}p_g\left( (1-r^2)/2\right)=
(-1)^{m+n}\vert p_g((1-r^2)/2)\vert,\tag6.10$$
\noindent
where $p_g(n)$ is defined by
$$\sum\,p_g(n)=q^{-1/2}\prod_{n>0}(1-q^{n/2})^{-(-1)^n8}.$$
Finally the denominator formula for the fake Conway sueralgebra
${\frak g}_{fC}$ is given by
$$e^{-\rho}\sum_{w\in W}\text{det}(w)\,w\left( e^{\rho}\prod_{n>0}
(1-e^{n\rho})^{(-1)^n8}\right)=\prod_{r\in L^+}(1-e^r)^{\text{mult}(r)},
\tag6.11$$
where $L^+$ denotes the set of positive roots.
\par\bigpagebreak

\head {\bf 6.2.\ Kac-Moody\ Algebras\ of\ the\ Arithmetic\ Type} 
\endhead               
\indent
Let $A=(a_{ij})$ be a symmetrizable generalized
Cartan matrix of degree $n$ and
let ${\frak g}(A)$ its associated Kac-Moody Lie algebra\,(\,see section 2\,).
\def\ep{\epsilon}
Then there exist a diagonal matrix $D=\text{diag}(\ep_1,\cdots,\ep_n)$
with $\ep_i>0,\ \ep_i\in\BQ\,(\,1\leq i\leq n\,)$ and a symmetric integral
matrix $B=(b_{ij})$ such that
$$A=DB,\ \ \ \text{g.c.d.}(\left\{\ b_{ij}\,\vert\ 1\leq i,j\leq n\right\})
=1.\tag 6.12$$
\def\a{\alpha}
We note that such matrices $D$ and $B$ are uniquely determined. Let
$$Q:=\sum_{i=1}^n\BZ\,\a_i,\ \ \ Q_+:=\sum_{i=1}^n\BZ_+\a_i,\ \ \
Q_-:=-Q_+,$$
where $\a_1,\cdots,\a_n$ are simple roots of $A$ or ${\frak g}(A).$ Then
$Q=Q_+\cup Q_-$ is a root lattice of $A.$ \par\medpagebreak\indent
Now we have the {\it canonical\ symmetric} bilinear form
$$(\,,\,)\,:\,Q\times Q\longrightarrow \BZ,\ \ \ (\a_i,\a_j)=b_{ij}=
a_{ij}/\ep_i.\tag 6.13$$
Let $\Delta,\,\Delta^+$ and $\Delta^-$ be the set of all roots, positive
roots, and negative roots of ${\frak g}(A)$ respectively. We let $W$ be the
Weyl group of ${\frak g}(A)$ generated by the fundamental reflections
$$r_{\a_i}(\beta):=\beta-2{{(\beta,\a_i)}\over {(\a_i,\a_i)}}\a_i,
\ \ \ \beta\in Q,\ \ 1\leq i\leq n.\tag 6.14$$
It is clear that $\Delta$ is invariant under $W$. We let
$$K:=\{\,\a\in Q_+\,\vert\ \a\ne 0,\ (\a,\a_i)\leq 0\ \text{for\ all}\ i,
\ \text{and\ supp}(\a)\ \text{is\ connected}\,\},\tag 6.15$$
where for $\a=\sum_{i=1}^nk_i\a_i\in Q_+,\ \text{supp}(\a)$ is defined to
be the subset $\{\,\a_i\,\vert\,k_i>0\,\}$ of the set
$\{\,\a_1,\cdots,\a_n\,\},$ and $\text{supp}(\a)$ is said to be
{\it connected} if there do not exist nonempty two sets $A_1$ and $A_2$
such that $\text{supp}(\a)=A_1\cup\ A_2$ and $(\a,\beta)=0$ for all
$\a\in A_1$ and $\beta\in A_2.$ Let $\Delta^{re}$\,(\,resp.\,$\Delta_{im}$\,)
be the set of all real roots\,(\,resp.\,imaginary\,) roots of ${\frak g}(A)$.
Then it is easy to check that
$$\Delta^{re}=W(\a_1)\cup\cdots\cup W(\a_n)\tag 6.16$$
\noindent
and
$$\Delta^{im}\cap Q_+=W(K).\tag 6.17$$
\noindent
{\bf Definition\ 6.4.} A generalized Cartan matrix $A$ of degree $n$ or
its associated Kac-Moody Lie algebra ${\frak g}(A)$ is said to be
{\it of\ the\ arithmetic\ type} or {\it have\ the\ arithmetic\ type} if
it is symmetrizable and indecompsable and also if for each $\beta\in Q$
with the property $(\beta,\beta)<0$ there exist a positive integer
$n(\beta)\in \BZ^+$ and an imaginary root $\a\in \Delta^{im}$ such that
$$n(\beta)\beta\equiv \a\ \ \ \text{mod}\,Q_0\ \ \text{on}\ Q,\tag 6.18$$
where $Q_0:=\{\,\gamma\in Q\,\vert\ (\gamma,\delta)=0\ \text{for\ all}\
\delta\in Q\,\}$ denotes the kernel of $(\,,\,).$
\par\bigpagebreak\indent
If we set $M:=Q/Q_0,$
then $(,\,\,)$ induces the {\it canonical\ nondegenerate,\ symmetric\
integral} bilinear form on the free $\BZ$-module $M$ defined by
$$S\,:\,M\times M\longrightarrow \BZ.\tag 6.19$$
We let $\pi\,:\,Q\longrightarrow M$ be the projection of $Q$ onto $M$, and
we denote by ${\tilde x}=\pi(x)$ the image of $x\in Q$ under $\pi.$
We denote by $(t_+,t_-,t_0)$ the signature of a symmetric matrix $B.$
\par\medpagebreak\indent
The following theorem is due to V.\,V. Nikulin.\par\smallpagebreak
\noindent
{\bf Theorem\ 6.5\,(\,[N5],\,Theorem\ 2.1\,).} A symmetrizable
indecomposable generalized Cartan matrix $A$ or its associated Kac-Moody
Lie algebra ${\frak g}(A)$ has the {\it arithmetic\ type} if and only if
$A$ has one of the following types $(a),\,(b),\,(c)$ or $(d)$\,:
\par\medpagebreak\indent
(a) The finite type case\,:\ $B>0.$\par\smallpagebreak\indent
(b) The affine type case\,:\ $B\geq 0$ and $B$ has the signature
$(\ell,0,1)$.\par\smallpagebreak\indent
(c) The rank $2$ hyperbolic case\,: $B$ has the signature $(1,1,0).$
\par\smallpagebreak\indent
(d) The arithmetic hyperbolic type\,: $B$ is hyperbolic of rank $>\,2,$
equivalently, $B$ has the signature $(\ell-1,1,k)$ with $\ell\geq 3,$ and
the index $[O(S):{\tilde W}]$ is finite.\par\medpagebreak
\noindent
Here $O(S)$ and ${\tilde W}$ denote the orthogonal group of $S$ and the
image of the Weyl group $W$ under $\pi$ respectively.
\par\medpagebreak\indent
Now we assume that $A$ is of the arithmetic hyperbolic type and that
$B$ has the signature $(t_+,1,k)$ with $t_+\geq 2.$ We choose a subgroup
${\tilde W}$ of $W(S)$ of finite index generated by reflections. We choose
a fundamental polyhedron ${\Cal M}$ of ${\tilde W}$, and then let
$P({\Cal M})_{pr}$ be the set of primitive elements of ${\Cal M}$ which are
orthogonal to the faces of ${\Cal M}$ and directed outside.
\par\medpagebreak\noindent
{\bf Theorem\ 6.6\,(\,[N5],\,Theorem\ 4.5\,).} We assume that
$S\,:\,M\times M\longrightarrow \BZ$ is a reflexive primitive hyperbolic,
symmetric integral bilinear form and that ${\tilde W}\subset W(S)$
satisfies the following conditions (6.20) and (6.20)\,:
$$P({\Cal M})_{pr}\ \text{generates}\ M\,;\tag 6.20$$
$$P({\Cal M}_0)_{pr}\ \text{generates} \ M,\tag 6.21$$
where ${\Cal M}_0
\ \text{is\ the\ fundamental\ polyhedron\ of}\ W(S).$
\par\medpagebreak\indent
In additon, we assume that we have a function
$$\lambda\,:\,P({\Cal M})_{pr}\longrightarrow \BZ^+$$
satisfying the conditions (6.22) and (6.23)\.:
\def\l{\lambda}
$$S(\l(\a)\a,\l(\a)\a)\ \text{divides}\ 2\,S(\l(\beta)\beta,\l(\a)\a)\
\text{for\ all}\ \a,\beta\in P({\Cal M})_{pr}\,;\tag 6.22$$
$$\{\,\l(\a)\a\,\vert\ \a\in P({\Cal M})_{pr}\,\}\ \text{generates}\ M.
\tag 6.23$$
Then the data $(S,{\tilde W},\l)$ defines canonically a generalized Cartan
matrix of the arithmetic hyperbolic type
$$A(S,{\tilde W},\l)=\left(\,2\,S(\l(\beta)\beta,\,\l(\a)\a)/S(\l(\a)\a,
\l(\a)\a\right),\ \ \ \a,\,\beta\in P({\Cal M})_{pr}.$$
\noindent
{\smc Remark\ 6.7.} (a) According to Nikulin\,(\,cf.\,[N3],\,[N4]\,) and
Vinberg\,(\,cf.\,[V]\,), there exist only a finite number of isomorphism
classes of reflexive primitive hyperbolic symmetric integral bilinear
forms $S$ of rank $\geq 3$, and the rank of $S$ is less than $31.$
Therefore by Theorem 6.6, there are only {\it finite} Kac-Moody Lie
algebras of the arithmetic hyperbolic type.\par\medpagebreak\noindent
(2) In [K], a very special case of a generlized Cartan matrix $A$ is
considered. This matrix is called just {\it hyperbolic} there. This has
the property that the fundamental polyhedron ${\Cal M}$ of ${\tilde W}$
is a simplex. There exist only a finite list of these hyperbolic ones.
These are characterized by the property\,:\ $0\ne 0\in Q$ is an imaginary
root if and only if $(\a,\a)<0.$\par\medpagebreak\noindent
(c) The complete list of the bilinear forms mentioned in Theorem 6.6 is not
known yet.\par\medpagebreak
\noindent
{\bf Example\ 6.8.} We consider an example of s {\it symmetric}
generalized Cartan matrix $A$ of the arithmetic hyperbolic type given by
$$A=(a_{ij})=\pmatrix \ 2 & -2 & \ 0\\ -2 & \ 2 & -1\\ \ 0 & -1 & \ 2
\endpmatrix.\tag 6.22$$
Let ${\Cal F}:={\frak g}(A)$ be its associated Kac-Moody
Lie algebra of the arithmetic hyperbolic type.
\def\F{\Cal F}
Let $\F_0$ be the {\it affine} Kac-Moody Lie algebra of type $A_1^{(1)}$
with its Cartan matrix $A_0=\pmatrix\ 2 & -2\\ -2 & \ 2 \endpmatrix.$
Then it is known that
$$\F_0\cong {s\ell}_2(\BC)\otimes \BC [t,t^{-1}]\oplus \BC\!\cdot\!c,$$
which is a one-dimensional central extension of the loop algebra
${s\ell}_2(\BC)\otimes \BC [t,t^{-1}].$ We let $\F_0^e$ be the
semi-direct product of $\F_0$ and $\BC\!\cdot\! d,\ d:=-t{{d}
\over {dt}},$ whose bracket is defined as follows\,:
$$\align
&[x\otimes t^n,y\otimes t^m]=[x,y]\otimes t^{m+n}+n<x,y>\,\delta_{n,-m}c,
\ \ m,n\in \BZ,\\
&[d,\,x\otimes t^n]=-n(x\otimes t^n),\ \ n\in \BZ,\\
& [c,a]=0\ \text{for\ all}\ a\in \F_0^e,\ \text{i.e.,}\ c\
\text{acts\ centrally,}\endalign$$
where $x,y\in s\ell_2(\BC)$ and
$<x,y>:=\,\text{tr}\,(xy)\,
=\,1/4\,\text{tr\,(\,ad\,}x\,\text{ad}\,y\,)$ denotes the Cartan-Killing form
on the Lie algebra ${s\ell}_2(\BC).$ The fact that the Weyl group of $\F$
is isomorphic to $PGL(2,\BZ)$ implies that $\F$ is closely related to the
theory of classical modular forms. In [F-F], Feingold and Frenkel
constructed $\F$ concretely and computed the Weyl-Kac denominator formula
for $\F$ explicitly. The denominator formula for $\F$ is given by
\def\s{\sigma}
$$\multline
\sum_{g\in PGL(2,\BZ)}\text{det}\,(g)\,e^{2\pi i\s (gP\,^t\!gZ)}\\
=e^{2\pi i\s (PZ)}\,\prod_{0\leq N\in S_2(\BZ)}
\left(1-e^{2\pi i\s(NZ)}\right)^{\text{mult}(N)}\,\prod_{N\in R}
\left(1-e^{2\pi i\s (NZ)}\right),\endmultline \tag 6.23$$
where $P=\pmatrix 3 & 1/2 \\ 1/2 & 2 \endpmatrix,\ Z=\pmatrix \tau & z\\
z & \omega \endpmatrix \in H_2,\ S_2(\BZ)$ denotes the set of all symmetric
integral matrices of degree $2$ and
$$R:=\left\{\, N=\pmatrix n_1 & n_3 \\ n_3 & n_2 \endpmatrix\in S_2(\BZ)\,
\vert\ n_1n_2-n_3^2=-1,\ n_2\geq 0,\ n_3\leq n_1+n_2,\ 0\leq n_1+n_2\,
\right\}.$$
We note that the root lattice of $\F$ is isomorphic to $S_2(\BZ).$
\par\medpagebreak\indent
Let ${\frak h}:=\BC h_1\oplus \BC h_2\oplus \BC h_3$ be a Cartan subalgebra
of $\F$. We denote by $\a_1,\,\a_2,\,\a_3$ the elements of ${\frak h}^\ast$
defined by
$$\a_i(h_j)=a_{ij},\ \ \ 1\leq i,j\leq 3.\tag 6.24$$
We put
\def\g{\gamma}
$$ \g_1^{\ast}:=\a_1/2,\ \ \g_2^{\ast}:=-\a_1-\a_2-\a_3,\ \
\g_3^{\ast}:=-\a_1-\a_2$$
and
$$P^{++}:=\{\,n_1\g_1^{\ast}+n_2\g_2^{\ast}+n_3\g_3^{\ast}\,\vert\
n_1,n_2,n_3\in \BZ_+,\ n_3\geq n_2\geq n_1\geq 0\,\}.$$
\noindent
{\smc Definition.} (1) The number $m_1+m_2$ in the weight
$\l=m_1\g_1^{\ast}+(m_1+m_2)\g_2^{\ast}+m_3\g_3^{\ast}$ is called the
{\it level} of the weight $\l.$\par\smallpagebreak\noindent
(2) An irreducible standard $\F_0^{e}$-module or its character is called
$\F$-{\it dominant} if the highest weight of this module lies in $P^{++}$.
A $\F_0^e$-module or its character is called $\F$-{\it dominant}
if each irreducible
standard component is $\F$-dominant.\par\medpagebreak
\indent
\def\G{\Gamma}
We let $M_k$ be the complex vector space spanned by those $\F_0^e$-characters
of the form
$$\chi(\tau,z,\omega)=\sum_{m\geq 0}\chi_m(\tau,z,\omega),$$
where for each $m\geq 0$, $\chi_m$ is the function satisfying the condition
$$\chi_m(\tau,z,\omega)=(-\tau)^{-k}\chi_m(-1/\tau,-z/\tau,
\omega-z^2/\tau),\ \ \ \pmatrix \tau & z\\ z & \omega\endpmatrix\in
H_2.\tag 6.26$$
Let $M_k(m)$ be the subspace of $M_k$ spanned by the $\F_0^e$-characters
of level $m$ satisfying the condition (6.26). We recall the results of
J. Igusa\,(\,cf.\,[Ig1]\,) on Siegel modular forms of degree $2.$ We denote
by $[\G_2,k]\,(\,\text{resp.}\,[\G_2,k]_0\,)$ the complex vector space
of all Siegel modular forms\,(\,resp.\,cusp forms\,) of weight $k$ on
$\G_2.$ Let $E_k\,(\,k\geq 4,\ k:\text{even}\,)$ be the Eisenstein series
of weight $k$ on $\G_2$ defined by
$$E_k(Z):=\sum_{C,D}\,\text{det}\,(CZ+D)^{-k},\ \ Z\in H_2,\tag 6.27$$
where $(C,D)$ runs over the set of non-associated pairs of coprime
symmetric matrices in $\BZ^{(2,2)}.$ Igusa proved that $E_4,\,E_6,\,E_{10}$
and $E_{12}$ are algebraically independent over $\BC$ and that
$$\oplus_{k=0}^{\infty}[\G_2,k]=\BC [E_4,E_6,E_{10},E_{12}].\tag 6.28$$
We define two cusp forms $\chi_{10}$ and $\chi_{12}$
of weight $10$ and $12$ by
$$\chi_{10}:=-43876\cdot 2^{-12}\cdot 3^{-5}\cdot 5^{-2}\cdot 7^{-1}
\cdot 53^{-1}\,(E_4E_6-E_{10})\tag 6.29$$
and
$$\chi_{12}:=131\cdot 593\cdot 2^{-13}\cdot 3^{-7}\cdot 5^{-3}\cdot
7^{-2}\cdot 337^{-1}(3^2\cdot 7^2 E_4^3-2\cdot 5^3 E_6^2-691 E_{12}).
\tag 6.30$$
Then according to (6.28), we have
$$\oplus_{k=0}^{\infty}[\G_2,k]=\BC [E_4,E_6,\chi_{10},\chi_{12}].\tag 6.31$$
For two nonnegative integers $k,m\geq 0,$ we define the set
$$S(k,m):=\left\{\,(a,b,c,d)\in (\BZ_+)^4\,\vert\ k=4a+6b+10c+12d,\ c+d=m\,
\right\}.\tag 6.32$$
We define the subspace $[\G_2,k](m)$ of $[\G_2,k]$ by
$$[\G_2,k](m):=\sum_{(a,b,c,d)\in S(k,m)}\BC\,E_4^aE_6^b\chi_{10}^c
\chi_{12}^d.\tag 6.33$$
Obviously $[\G_2,k]=\sum_{m\geq 0}[\G_2,k](m).$\par\medpagebreak\indent
For $f(\tau,z,\omega)\in [\G_2,k],$ we let
$$f(Z):=f(\tau,z,\omega)=\sum_{m\geq 0}\phi_m(\tau,z)\,e^{2\pi im\omega},
\ \ Z=\pmatrix \tau & z \\ z & \omega \endpmatrix\in H_2 \tag 6.34$$
be the Fourier-Jacobi expansion of $f$. As noted in section $4,\
\phi_m(\tau,z)$ is a Jacobi form of weight $k$ and $m$. Now
for each non-negative integer $m\geq 0$ we define the linear map
$L_m\,:\,[\G_2,k]\longrightarrow M_k$ by
$$(L_mf)(Z)=(L_mf)(\tau,z,\omega):=\phi_m(\tau,z)\,e^{2\pi im\omega},
\ \ f\in [\G_2,k],\ \ Z=\pmatrix \tau & z\\ z & \omega\endpmatrix\in H_2,
\tag 6.35$$
where $\phi_m(\tau,z)\,(\,m\geq 0\,)$ is the Fourier-Jacobi coefficient of
the expansion (6.34) of $f.$\par\medpagebreak
\noindent
{\smc Definition.} Let $M_k'\,(\,\text{resp.}\,M_k'(m)\,)$ be the subspace
of $M_k\,(\,\text{resp.}\,M_k(m)\,)$ consisting of $PSL(2,\BZ)$-invariant
$\F_0^e$-characters which are $\F$-dominant.\par\medpagebreak
\indent
In [F-F],\,Theorem 7.9, Feingold and Frenkel showed that $L_m$ maps
$[\G_2,k](m)$ isomorphically onto $M_k'(m).$ Thus according to (6.33), we
obtain , for each $m\geq 0,$
$$\text{dim}_{\BC}\,[\G_2,k](m)=\text{dim}_{\BC}M_k'(m)=\sharp(S(k,m)),
\tag 6.36$$
where $\sharp(S)$ denotes the cardinality of the set $S$. Moreover we have the
following ring-isomorphism
$$M'=\sum_{k\geq 0}\sum_{m\geq 0}M_k'(m)\cong \BC\,[E_4,E_6,\chi_{10},
\chi_{12}].\tag 6.37$$
Let $[\G_2,k]^M$ be the Maass space.\,(\,See Appendix B.) Maass showed
that
$$[\G_2,k]^M=\BC\,E_4\oplus [\G_2,k]_0\ \ \text{and}\ \
\text{dim}_{\BC}\,[\G_2,k]_0=\sharp(S(k,1)).$$
Also Maass showed that
$$[\G_2,k]^M\cong M_k(1)\ \ \ \text{and}\ \ \ [\G_2,k]_0\cong M_k'(1).
\tag 6.38$$
The detail for (6.38) can be found in Appendix B,\,[E-Z] and [Ma2-4].
For $k\geq 4,\,\text{even}$, we have the simple dimensional formulas
$$\text{dim}_{\BC}\,M_k(1)=\left[\frac{k+2}6\right]\ \ \ \text{and}
\ \ \ \text{dim}_{\BC}\,M_k'(1)=\left[\frac{k-4}6\right]=\sharp(S(k,1)).
\tag 6.39$$

\par\bigpagebreak
\head {\bf 6.3.\ Open\ Problems} 
\endhead 
\indent
In this subsection, we give some open problems which should be
investigated and give some comments. Those of Problem 1-6 are due to
R. Borcherds\,(\,cf.\,[Bo6-7]\,).\par\medpagebreak
\indent
{\bf Problem\ 1.} Can the methods for constructing automorphic forms as
infinite products in section $5$ be used for semisimple Lie groups other
than $O_{s+2,2}(\BR)$\,?\par\medpagebreak
{\bf Problem\ 2.} Are there a finite or infinite number of singular
automorphic forms that can be written as modular products\,? Are there such
singular modular forms on $O_{s+2,2}(\BR)$ for $s>24$\,?
\par\medpagebreak
{\bf Problem\ 3.} Interprete the automorphic forms that are modular
products in terms of representation theory or the Langlands philosophy.
\par\medpagebreak
{\bf Problem\ 4.} Extend Theorem 5.4 to higher levels.
\par\medpagebreak
{\bf Problem\ 5.} Investigate the Lie algebras and the superalgebras
coming from other elements of the MONSTER $G$ or Aut\,$(\L)$ and
write down their denominator formulas explicitly in some nice form.
\par\medpagebreak
{\bf Problem\ 6.} Are there any generalized Kac-Moody algebras other than
the finite dimensional, affine, monstrous or fake monstrous ones, whose
simple roots and root multiplicities can both be described explicitly\,?
\par\medpagebreak
{\bf Problem\ 7.} Given a generalized Cartan matrix $A$ of the arithmetic
hyperbolic type, construct its associated Kac-Moody Lie algebra
${\frak g}(A)$ of the same type explicitly. Give a relationship between
the Kac-Moody algebras of the arithmetic hyperbolic type and classical
mathematics. For instance, M. Yoshida showed that the Weyl group $W(A)$ of
${\frak g}(A)$ of rank $3$ are all hyperbolic triangle groups and that
the semidirect product of the Weyl group $W(A)$ and the root lattice of
${\frak g}(A)$ is isomorphic to a discrete subgroup of a parabolic subgroup
of $Sp(2,\BR).$\par\medpagebreak
{\bf Problem\ 8.} Develope the theory of Kac-Moody Lie algebras of the
arithmetic hyperbolic type geometrically.
\par\medpagebreak
{\bf Problem\ 9.} Give an analytic proof of the denominator formula (6.23)
for $\F$ analogous to that of the Jacobi's triple product identity.
\par\medpagebreak
{\bf Problem\ 10.} Find the transformation behaviour of the denominator
formula (6.23) for $\F$ under the symplectic involution $Z\rightarrow
-Z^{-1}.$\par\medpagebreak
{\bf Problem\ 11.} Apply the theroy of the Kac-Moody Lie algebra $\F$ to
the study of the moduli space of principally polarized abelian surfaces.
\par\medpagebreak
{\bf Problem\ 12.} Generalize the Maass correspondence to the Kac-Moody
algebras of the arithmetic hyperbolic type other than $\F$\,?

\redefine\g{\frak g}
\redefine\h{\frak h}
\redefine\a{\alpha}
\redefine\v{\vee}
\redefine\l{\lambda}
\redefine\D{\Delta}
\redefine\m{\underline m}

\head {\bf Appendix A. \  Classical Modular Forms}  
\endhead  

\redefine\Z{\Bbb Z}
\redefine\G{\Gamma}

\redefine\M{\Cal M}

Here we present some well-known results on modular forms whose proofs can be
found in many references, e.g., [Kob], [Ma1], [S], and [T].
\par
Let $H_1$ be the upper half plane and let $\G:=SL(2,\Z)$ be the elliptic
modular group. For a positive integer $N\in\Z^+,$ we define
$$\G(N):=\left\{\pmatrix a&b\\c&d\endpmatrix\in\G\;\biggl|\;\pmatrix a&b\\
c&d\endpmatrix\equiv\pmatrix 1&0\\0&1\endpmatrix\;(\bmod N)\right\}$$
and
$$\G_0(N):=\left\{\pmatrix a&b\\c&d\endpmatrix\in\G\;\biggl|\;c\equiv0
\;(\bmod N)\right\}$$
$\G(N)$ (resp. $\G_0(N)$) is called the {\it principal congruence subgroup} of
 {\it level} $N$
(resp. the {\it Hecke subgroup} of level $N$). The subgroup $\G_\theta$ of $\G$ generated
by $\pm\pmatrix 1&2\\0&1\endpmatrix$ and $\pm\pmatrix 0&1\\-1&0\endpmatrix$
is called the {\it theta group}. A subgroup $\G_1$ of $\G$ is called a
{\it congruence subgroup} if $\G_1$ contains $\G(N)$ for some positive integer $N.$ For instance,
the Hecke subgroup $\G_0(N)$ is a congruence subgroup because $\G(N)\subset
\G_0(N)\subset\G.$ And $\G(N)$ is a normal subgroup
because it is the kernel of the reduction-modulo-$N$ homomorphism $SL(2,\Z)
\rightarrow SL(2,\Z/N\Z).$
It is well known that the index of $\G(N)$ in $\G$ is given by
$$[\G:\G(N)]=N^3\prod_{p|N}(1-p^{-2}).\tag A.1$$
The proof of (A.1) can be found in [Sh] pp.21-22. It was discovered
around the 1880s that there are an infinite number of examples of
noncongruence subgroups (cf. [Ma1] pp. 76-78). But $SL(n,\Z)$
behaves quite differently for $n\geq3.$
In fact, it has been proved that if $n\geq 3,$ every subgroup of $SL(n,\Z)$ of
finite index is
a congruence subgroup (cf. [Bas]). A similar result for the Siegel modular
group $Sp(n,\Z)$ for $n\geq2$ can be found in [Me].
\par
For an integer $k\in\Z,$ we denote by $[\G,k]$ (resp. $[\G,k]_0$) the vector
space of all modular forms (resp. cusp forms) of weight $k$ for the elliptic
modular group $\G.$ Only for $k\geq0,$ $k$ even, $[\G,k]$ does not vanish.
\par
For any positive integer $k$ with $k\geq 2,$ we put
$$G_{2k}(\tau):=\sum_{m,n}\,'\frac1{(m\tau+n)^{2k}},\qquad\tau\in H_1.\tag A.2$$
Here the symbol $\sum'$ means that the summation runs over all pair of
integers $(m,n)$ distinct from (0,0). Then $G_{2k}\in[\G,2k]$ and
$G_{2k}(\infty)=2\zeta(2k),$ where $\zeta(s)$ denotes the Riemann zeta
function. $G_{2k}\;(k\in\Z^+,\;k\geq 2)$ is called the {\it Eisenstein series} of
{\it index} $2k.$
The Fourier expansion of $G_{2k}\;(k\geq 2)$ is given by
$$G_{2k}(\tau)=2\zeta(2k)+2\frac{(2\pi i)^{2k}}{(2k-1)!}\sum_{n=1}^\infty\sigma
_{2k-1}(n)q^n,\qquad\tau\in H_1, \tag A.3$$
where $q=e^{2\pi i\tau}$ and $\sigma_\ell(n):=\sum_{0<d|n}d^\ell.$
\par
We consider the following parabolic subgroup $P$ of $\G$ given by
$$ P:=\left\{\pmatrix a&b\\0&d\endpmatrix\in\G\right\}.$$
Then we can see easily that
$$G_{2k}(\tau)=2\zeta(2k)\sum_{\gamma\in  P\backslash\G}
\left(\frac{d(\gamma<\tau>)}{d\tau}\right)^k.$$
Here for $\gamma=\pmatrix a&b\\ c&d\endpmatrix\in\G$ and $\tau\in H_1,$ we set
$$\gamma<\tau>:=(a\tau+b)(c\tau+d)^{-1}.$$
For a positive integer $k\geq 2,$ we can see easily that
$$E_{2k}(\tau):\overset \text{def}\to =\frac{G_{2k}(\tau)}{2\zeta(2k)}=1-
\frac{4k}{B_{2k}}\sum_{n=1}^\infty\sigma_{2k-1}(n)q^n,\tag A.4$$
where $B_k\;(k=0,1,2,\cdots )$ donotes the $k-$th Bernoulli number defined by the formal
power series expansion:
$$\frac{x}{e^x-1}=\sum_{k=0}^\infty B_k\frac{x^k}{k!}.$$
Then clearly $B_{2k+1}=0$ for $k\geq1.$
The first few $B_k$ are
$$\align
&B_0=1,\;B_1=-1/2,\;B_2=1/6,\;B_4=-1/30,\;B_6=1/42,\;B_8=-1/30,\;B_{10}=5/66,\\
&B_{12}=-691/2730,\;B_{14}=7/6,\;B_{16}=-3617/510,\;B_{18}=43867/798,
\cdots
\endalign$$
Indeed, (A.4) follows immediately from the relation
$$\zeta(2k)=(-1)^{k-1}\frac{2^{2k-1}\pi^{2k}}{(2k)!}B_{2k},\quad k=1,2,\cdots.$$
For example,
$$\align
&E_4(\tau)=1+240\sum_{n=1}^\infty\sigma_3(n)q^n,\quad (240=2^4\cdot3\cdot5)\\
&E_6(\tau)=1-504\sum_{n=1}^\infty\sigma_5(n)q^n,\quad (504=2^3\cdot3^2\cdot7)\\
&E_8(\tau)=1+480\sum_{n=1}^\infty\sigma_7(n)q^n, \quad (480=2^5\cdot3\cdot5)\\
&E_{10}(\tau)=1-264\sum_{n=1}^\infty\sigma_9(n)q^n,\quad (264=2^3\cdot3\cdot11)\\
&E_{12}(\tau)=1+\frac{65520}{691}\sum_{n=1}^\infty\sigma_{11}(n)q^n,\quad
(65520=2^4\cdot3^2\cdot5\cdot7\cdot13)\\
&E_{14}(\tau)=1-24\sum_{n=1}^\infty\sigma_{13}(n)q^n,\quad (24=2^3\cdot3).
\endalign$$
According to the argument on the dimension of $[\G,k],$ we obtain the relation
$$E_4^2=E_8,\quad E_4E_6=E_{10}.\tag A.5$$
These are equivalent to the identities:
$$\align
\sigma_7(n)&=\sigma_3(n)+120\sum_{n=1}^{n-1}\sigma_3(m)\sigma_3(n-m),\\
11\sigma_9(n)&=21\sigma_5(n)-10\sigma_3(n)+5040\sum_{m=1}^{n-1}\sigma_3(n)
\sigma_5(n-m).
\endalign$$
More generally, every $E_{2k}$ can be expressed as a polynomial in $E_4$ and
$E_6.$ For instance, $E_{14}=E_4^2E_6.$\par
We put
$$g_2:=60G_4,\;\;\text{and}\;\;g_3:=140G_6.\tag A.6$$
Then it is obvious that
$$g_2=\frac{(2\pi)^4}{2^2\cdot3}E_4,\;\;\text{and}\;\;g_3=\frac{(2\pi)^6}
{2^3\cdot3^3}E_6.$$
Since $g_2(\infty)=\frac43\pi^4$ and $g_3(\infty)=\frac8{27}\pi^6,$
we see that the {\it discriminant}
$$\tilde{\D}:=g_2^3-27g_3^2\tag A.7$$
is a cusp form of weight 12, that is, $\tilde{\D}\in[\G,12]_0.$
And we have
$$\aligned
\tilde{\D}(\tau)&=(2\pi)^{12}\cdot2^{-6}\cdot3^{-3}(E_4(\tau)^3-E_6(\tau)^2)\\
&=(2\pi)^{12}(q-24q^2+252q^3-1472q^4+\cdots)\\
&=(2\pi)^{12} q\prod_{n=1}^\infty(1-q^n)^{24}\qquad(\text{Jacobi's identity}).
\endaligned\tag A.8$$
In this article, we put $\D(\tau) := (2\pi)^{-12}\tilde{\D}(\tau) =
q\prod_{n=1}^\infty (1 - q^n)^{24}.$\par
Fix $\tau\in H_1.$ The Weierstrass $\wp$-function $\wp(z;\tau)$ is defined by
$$\wp(z;\tau):=\frac1{\tau^2}+\sum_{m,n}\,'\left\{
\frac{1}{(z-n-m\tau)^2}-\frac1{(n+m\tau)^2}\right\},\quad z\in\Bbb C.\tag A.9$$
Then $\wp(z;\tau)$ is a meromorphic function with respect to $1,\tau$ with
double poles at the points $n+m\tau,\;n,m\in\Z.$ The map $\varphi_\tau
:\Bbb C\rightarrow \Bbb P^2$ defined by
$$\varphi_\tau(z):=[1:\wp(z;\tau):\frac{d}{dz}\wp(z;\tau)],\quad z\in\Bbb C
\tag A.10$$
induces an isomorphism of $X=\Bbb C/L_\tau$ with a nonsingular plane curve of
the form
$$X_0X_2^2=4X_1^3+aX_0^2X_1+bX_0^3,\tag A.11$$
where $a$ and $b$ are suitable constants depending on $\tau$ and $L_\tau:=
\{m\tau+n\mid m,n\in\Bbb Z\}$ is the lattice in $\Bbb C$
generated by 1 and $\tau.$
If we put $x=\wp(z;\tau)$ and $y=\frac{d}{dz}\wp(z;\tau),$
we have the differential equation
$$y^2=4x^3-g_2(\tau)x-g_3(\tau).\tag A.12$$
Up to a numerical factor, $\tilde{\D}(\tau):=(g_2^3-g_3^2)(\tau)$
is the discriminant of the polynomial $4x^3-g_2(\tau)x-g_3(\tau).$
Since $\tilde{\D}(\tau)\neq0,\;X_\tau=\Bbb C^2/L_\tau$ is a {\it nonsingular} elliptic curve.
This story tells us as the reason why the function $\D$ is called the discrimant.
We observe that the differential equation (A.12) shows that it is the inverse
function for the elliptic integral in Weierstrass normal form, that is,
$$z-z_0=\int_{\wp(z_0;\tau)}^{\wp(z;\tau)}(4w^3-g_2(\tau)w-g_3(\tau))^{-1/2}dw.
\tag A.13$$
The Ramanujan tau function $\tau(n)\;(n\in\Z^+)$ is defined by
$$\D(\tau)=(2\pi)^{-12}\tilde{\D}(\tau)=q\prod_{n=1}^\infty(1-q^n)^{24}=\sum_{n=1}^\infty\tau(n)
q^n.\tag A.14$$
The Dedekind eta function $\eta(\tau)$ is defined by
$$\eta(\tau):=q^{1/24}\prod_{n\geq1}(1-q^n),\qquad\tau\in H_1.\tag A.15$$
Then $\eta(\tau)$ satisfies
$$\eta(\tau+1)=\eta(\tau)\qquad\text{and}\qquad\eta(-1/\tau)=(\tau/i)^{1/2}\eta
(\tau).$$
The Dedekind eta function $\eta(\tau)$ is related to the partition function
$p(n)$ as follows:
$$q^{1/24}\eta(\tau)^{-1}=\prod_{n\geq1}(1-q^n)^{-1}=\sum_{n\geq0} p(n)q^n,
\tag A.16$$
where $p(n)$ is the number of partitions of $n,$ i.e., the number of ways
of writing
$$n=n_1+\cdots+n_r,\quad n_j\in\Z^+\;\;(1\leq j\leq r).$$
\par
The {\it modular invariant} $J(\tau)$ is defined by
$$J:=(60G_4)^3/\tilde{\D}=g_2^3/\tilde{\D}=(2\pi)^{12}
\cdot 2^{-6}\cdot3^{-3}E_4^3/\tilde{\D},\tag A.17$$
The function $J(\tau)$ was first constructed by Julius Wilhelm Richard
Dedekind (1831-1916) in 1877 and Felix Klein (1849-1925) in 1878.
The modular invariant $J(\tau)$ has the following properties:
\roster
\item"(J1)" $J(\tau)$ is a modular function. $J$ is holomorphic in $H_1$
    with a simple pole at $\infty,$ $J(i)=1$ and $J\left(\frac{-1+\sqrt3i}2\right)=0.$
\item"(J2)" $J$ defines a conformal mapping which is one-to-one from
$H_1/\G$ onto $\Bbb C,$ and hence $J$ provides an identification of
$H_1/\G\cup\{\infty\}$ with the Riemann sphere $S^2=\Bbb C\cup\{\infty\}.$
\item"(J3)" The following are equivalent for a function $f$ which is
meromorphic on $H_1;$\newline
(a) $f$ is a modular function;\newline
(b) $f$ is a quotient of two modular forms of the same weight;\newline
(c) $f$ is a rational function of $J,$ i.e., a quotient of polynomials in $J.$
Thus $J$ is called the {\it Hauptmodul} or the {\it fundamental function.}
\endroster
\par
The $q$-expansion of $j(\tau):=1728J(\tau)=2^6\cdot3^3J(\tau),$
also called the modular invariant, is given by
$$j(\tau)=q^{-1}+744+\sum_{n=1}^\infty c(n)q^n=q^{-1}+744+196884q+21493760q^2+
\cdots.\tag A.19$$
We observe that $j(i)=1728=2^6\cdot3^3$ and $j\left(\frac{1+\sqrt3i}
{2}\right)=0.$ It was already mentioned that there is a surprising
connection of the coefficients in (A.19) with the representations of
the Fischer-Griess monster group. All of the early Fourier coefficients in
(A.19) are simple linear combinations of degrees of characters of the MONSTER.
This was first observed by John Mckay and John Thompson. The modular invariant $J(\tau)$ is used to prove the small Picard
theorem and to study an explicit reciprocity law for an imaginary quadratic
number field.
\par
For a positive definite symmetric real matrix $S$ of degree $n,$ we define the
theta series
$$\theta_S(\tau):=\sum_{x\in\Z^n}e^{\pi iS[x]\tau},\qquad\tau\in H_1,\tag A.20$$
where $S[x]:={}^t\!x Sx$ denotes the quadratic form associated to $S.$
We can prove the transformation formula
$$\theta_{S^{-1}}(-1/\tau)=(\det S)^{1/2}(\tau/i)^{n/2}\theta_S(\tau).\tag A.21$$
It is known that if $S$ is a positive definite symmetric even integral, unimodular
matrix of degree $n,$ then $n$ is divided by $8$ and $\theta_S(\tau)\in[\G,n/2].$
In fact, for $n=8,$ there is only one  positive definite symmetric even integral
unimodular matrix up to equivalence modulo $GL(8,\Z).$ For $n=16,$ there are
two nonequivalent examples modulo $GL(16,\Z).$ For $n=24,$ there are 24
nonequivalent examples modulo $GL(24,\Z).$
\par
We consider a Jacobi function
$$\theta(\tau,z):=\sum_{n\in\Z}e^{\pi i(n^2\tau+2nz)},\qquad
(\tau,z)\in H_1\times\Bbb C\tag A.22$$
Then $\theta(\tau,z)$ satisfies the following properties:
\roster
\item"($\theta$.1)" $\theta(\tau,z)$ is an entire function on $H_1\times\Bbb C.$
\item"($\theta$.2)" $\theta$ is quasi-periodic as a function of $z$
in the following sense:\newline
$\theta(\tau,z+n)=\theta(\tau,z)\quad\text{for all}\quad n\in\Bbb Z;$\newline
$\theta(\tau,z+n\tau)=e^{-\pi i(n^2\tau+2nz)}\theta(\tau,z)\quad\text{for all}
\quad n\in\Z.$
\item"($\theta$.3)" $\theta(\tau,z)$ satisfies the transformation formula
$$\theta(\tau,z)=(\tau/i)^{1/2}\sum_{n\in\Bbb Z}e^{-\pi i(n-z)^2/\tau}.$$
\item"($\theta$.4)" $\theta(\tau,(1+\tau)/2)=0.$
\item"($\theta$.5)" Fixing $\tau,$ the only zero of $\theta(z):=
\theta(\tau,z)$ as a function of $z$ in the period parallelogram on 1 and
$\tau$ is $z=(1+\tau)/2.$ Moreover, this zero is simple.
\endroster
For a proof of ($\theta.$3), use Poisson formula.

\head {\bf Appendix B. \ Kohnen Plus Space and Maass Space}  
\endhead   
\par
Here we review the Kohnen plus space and the Maass space. And then we give isomorphisms
of them with the vector spaces of Jacobi forms. For more detail we refer to
[Koh], [Ma2-4].\par
We fix two positive integers $n$ and $m.$ Let
$$H_n:=\{Z\in\Bbb C^{(n,n)}\mid Z={}^t\!Z,\;\text{Im}Z>0\}$$
be the Siegel upper half plane of degree $n$ and let $\G_n:=Sp(n,\Z)$
the Siegel modular group of degree $n.$ That is,
$$\G_n:=\{g\in\Z^{(2n,2n)}\mid {}^t\!gJg=J\},\quad J:=\pmatrix 0&E_n\\-E_n&0
\endpmatrix.$$
Here $E_n$ denotes the identity matrix of degree $n.$ Then the real
symplectic group $Sp(n,\Bbb R)$ acts on $H_n$ transitively. If $M=\pmatrix
A&B\\C&D\endpmatrix\in Sp(n,\Bbb R)$ and $Z\in H_n,$
$$M<Z>:=(AZ+B)(CZ+D)^{-1}.\tag B.1$$
\par
Let $\M$ be a positive definite, symmetric half integral matrix of degree $m.$
For a fixed element $Z\in H_n,$ we denote by $\Theta_{\M,Z}^{(n)}$
the vector space of all the functions $\theta:\Bbb C^{(m,n)}\rightarrow\Bbb C$
satisfying the condition
$$\theta(W+\l Z+\mu)=e^{-2\pi i\sigma(\M[\l]Z+2{}^t\!W\M\l)},\quad W\in\Bbb C^{(m,n)}
\tag B.2$$
for all $\l,\,\mu\in\Bbb Z^{(m,n)}.$ For brevity, we put $L:=\Z^{(m,n)}$ and
$\Cal L_\M:=L/(2\M)L.$ For each $\gamma\in\Cal L_\M:=L/(2\M)L,$
we define the theta series
$$\theta_\gamma(Z,W):=\sum_{\l\in L}e^{2\pi i\sigma(\M[\l+(2\M)^{-1}\gamma]
Z+2{}^t\!W\M(\l+(2\M)^{-1}\gamma))},\quad (Z,W)\in H_n\times\Bbb C^{(m,n)}.
\tag B.3$$
Then $\{\theta_\gamma(Z,W)\mid \gamma\in\Cal L_\M\}$ forms a basis for
$\Theta_{\M,Z}^{(n)}.$ For any Jacobi form $\phi(Z,W)\in J_{k,\M}(\G_n),$
the function $\phi(Z,\cdot)$ with fixed $Z$ is an element of $\Theta_{\M,Z}^{(n)}$
and $\phi(Z,W)$ can be written as a linear combination of theta series
$\theta_\gamma(Z,W) \;(\gamma\in\Cal L_\M):$
$$\phi(Z,W)=\sum_{\gamma\in\Cal L_\M}\phi_\gamma(Z)\theta_\gamma(Z,W).
\tag B.4$$
Here $\phi=(\phi_\gamma(Z))_{\gamma\in\Cal L_\M}$ is a vector valued
automorphic form with respect to theta mutiplier system.\newline\par\bigpagebreak

{\bf (I) Kohnen Plus Space} (cf. [Ib], [Koh])
\par\medpagebreak
We consider the case: $m=1,\;\M=E_m,\;L=\Z^{(1,n)}\cong\Z^n.$
We consider the theta series
$$\theta^{(n)}(Z):=\sum_{\l\in L}e^{2\pi i\sigma(\l Z{}^t\!\l)}=\theta_0(Z,0),
\quad Z\in H_n.\tag B.5$$
We put
$$\G_0^{(n)}(4):=\left\{\pmatrix A&B\\C&D\endpmatrix\in\G_n\;\biggl|\;
C\equiv0\,(\bmod\, 4)\right\}.\tag B.6$$
Then $\G_0^{(n)}(4)$ is a congruence subgroup of $\G_n.$ We define the
automorphic factor $j:\G_0^{(n)}(4)\times H_n\rightarrow\Bbb C^\times$ by
$$j(\gamma,Z):=\frac{\theta^{(n)}(\gamma<Z>)}{\theta^{(n)}(Z)},\quad
\gamma\in\G_0^{(n)}(4),\;Z\in H_n.\tag B.7$$
Then we obtain the relation
$$j(\gamma,Z)^2=\epsilon(\gamma)\cdot\det(CZ+D),\quad\epsilon(\gamma)^2=1
\tag B.8$$
for any $\gamma=\pmatrix A&B\\C&D\endpmatrix \in\G_0^{(n)}(4).$
\par
Now we define the Kohnen plus space $M_{k-\frac12}^+(\G_0^{(n)}(4))$
introduced by W. Kohnen (cf. [Koh]).
$M_{k-\frac12}^+(\G_0^{(n)}(4))$ is the vector space consisting of holomorphic
functions $f: H_n\rightarrow\Bbb C$ satisfying the following conditions:
\roster
\item"(a)" $f(\gamma<Z>)=j(\gamma,Z)^{2k-1}f(Z)$ for all $\gamma\in
\G_0^{(n)}(4);$
\item"(b)" $f$ has the Fourier expansion
$$f(Z)=\sum_{T\geq0}a(T)e^{2\pi i\sigma(TZ)},$$
\endroster
where $T$ runs over the set of semi-positive, half-integral symmetric matrices
of degree $n$ and $a(T)=0$ unless $T\equiv-\mu{}^t\!\mu\bmod 4S_n^*(\Z)$
for some $\mu\in\Z^{(n,1)}.$ Here we put
$$S^*_n(\Z):=\{T\in\Bbb R^{(n,n)}\mid T={}^t\!T,\;\sigma(TS)\in\Z\;\text{for
all}\; S={}^t\!S\in\Z^{(n,n)}\}.$$
For $\phi\in J_{k,1}(\G_n),$ according to (B.4), we have
$$\phi(Z,W)=\sum_{\gamma\in L/2L}f_\gamma(Z)\theta_\gamma(Z,W),\quad
Z\in H_n,\;W\in\Bbb C^n.\tag B.9$$
Now we put
$$f_\phi(Z):=\sum_{\gamma\in L/2L}f_\gamma(4Z),\quad Z\in H_n.\tag B.10$$
Then $f_\phi\in M_{k-\frac12}^+(\G_0^{(n)}(4)).$
\proclaim{Theorem 1 (Kohnen-Zagier $(n=1)$, Ibukiyama $(n>1)$)} \
Suppose $k$ is an even positive integer. We have
the isomorphism
$$\align
J_{k,1}(\G_n)&\cong M_{k-\frac12}^+(\G_0^{(n)}(4))\\
\phi&\mapsto f_\phi.
\endalign$$
Furthermore the isomorphism is compatible with the action of Hecke operators.
\endproclaim
\par\bigpagebreak
{\bf (II) \ Maass Space} \newline
\par\medpagebreak
The Maass space or the Maa\ss's Spezialschar was introduced by
H. Maass (1911-1993) to solve the Saito-Kurokawa conjecture.
Let $k\in\Z^+.$ We denote by $[\G_2,k]$ the vector space of all Siegel modular
forms of weight $k$ and degree 2.
We denote by $[\G_2,k]^M$ the vector space of all Siegel modular forms $F: H_2
\rightarrow\Bbb C,\;F(Z)=\sum_{T\geq0}a_F(T)e^{2\pi i\sigma(TZ)}$ in
$[\G_2,k]$ satisfying the following condition:
$$\aligned
a_F\pmatrix n&\frac r2\\ \frac r2&m\endpmatrix=\sum_{d|(n,r,m),\,d>0}&d^{k-1}
a_F\pmatrix \frac{mn}{d^2}&\frac{r}{2d}\\ \frac{r}{2d}&1\endpmatrix \\
&\text{for all}\;T=\pmatrix n&\frac{r}2\\ \frac{r}2&m\endpmatrix\geq 0
\;\text{with}\;n,r,m\in\Z.\endaligned\tag B.11$$
The vector space $[\G_2,k]^M$ is called the {\it Maass space} or the
{\it Maa\ss's Spezialschar.}\par
For any $F$ in $[\G_2,k],$ we let
$$F(Z) = \sum_{m\geq 0} \phi_m(\tau,z) e^{2 \pi i m \tau'}, \; Z=\pmatrix
\tau & z \\ z &\tau' \endpmatrix \in H_2 \tag B.12$$
be the Fourier-Jacobi expansion of $F.$ Then for any $m\in \Bbb Z_+$
we obtain the linear map
$$\rho_m : [\G_2,k] \longrightarrow J_{k,m}(\G_1), \quad F\longmapsto \phi_m.
\tag B.13$$
We denote that $\rho_0$ is nothing but the Siegel $\Phi-$operator.\par
Maass (cf. [Ma 2-3]) showed that for $k$ even, there exists a natural map
$ V: J_{k,1}(\G_1) \longrightarrow [\G_2,k]$ such that $\rho_1 \circ V$ is the
identity. More precisely, we let $\phi\in J_{k,1}(\G_1)$ with Fourier
coefficients $c(n,r) \;(n,r\in \Bbb Z,\; r^2 \leq 4n)$ and we define for any
$m\in \Bbb Z^+_{\geq 0}$
$$ (V_m\,\phi)(\tau,z) := \sum_{n,r\in \Bbb Z,\,\,r^2\leq 4mn}
\left( \sum_{d|(n,r,m)} \,d^{k-1}c\left(\frac{mn}{d^2}, \frac{r}{d}\right)
\right)\,e^{2\pi i(n\tau +rz)}.\tag B.13$$
It is easy to see that $V_1\phi =\phi$ and $V_m\phi \in J_{k,m}(\G_1).$
We define
$$(V\phi)\pmatrix \tau &z \\ z&\tau'\endpmatrix := \sum_{m\geq 0}
(V_m\phi)
(\tau,z) e^{2\pi i m \tau'},\quad \pmatrix \tau & z \\ z&\tau'\endpmatrix\in H_2.
\tag B.15$$
\par
We denote by $T_n\;(n\in \Bbb Z^+)$ the usual Hecke operators on $[\G_2,k]$ resp.
$[\G_2,k]_0.$ Here $[\G_2,k]_0$ denote the vector subspace consisting of all cusp
forms in $[\G_2,k].$ For instance, if $p$ is a prime, $T_p$ and $T_{p^2}$ are
the Hecke operators corresponding to the two generators
$\G_2\,\text{diag}(1,1,p,p)\,\G_2$ and $\G_2\,\text{diag}(1,p,p^2,p)\,\G_2$
of the
local Hecke algebra of $\G_2$ at $p$ respectively. We denote by $T_{J,n}\;(n\in
\Bbb Z^+)$ the Hecke operators on $J_{k,m}(\G_1)$ resp. $J_{k,m}^{\text{cusp}}
(\G_1)$ (cf. [E-Z]).
\proclaim
{{\bf Theorem 2} (Maass [Ma 2-4], Eichler-Zagier [E-Z], Theorem 6.3)} \ Suppose
$k$ is an even positive integer. Then the map $\phi \mapsto V \phi$
gives an injection of $J_{k,1}(\G_1)$ into $[\G_2,k]$ which sends cusp forms to
cusp forms and is compatible with the action of Hecke operators. The image of
the map $V $ is equal to the Maass space $[\G_2,k]^M.$ If $p$ is a prime,
one has
$$T_p\circ V = V \circ (T_{J,p} +p^{k-2}(p+1))$$
and
$$T_{p^2}\circ V = V \circ (T^2_{J,p} +p^{k-2}(p+1)T_{J,p} +p^{2k-2}).$$
\endproclaim
In summary, we have the following isomorphisms
$$
\matrix
[\G_2,k]^M &\cong & J_{k,1}(\G_1) & \cong &
M^+_{k-\frac12}(\G^{(1)}_0(4)) & \cong & [\G_1,2k-2], \\
V\phi  & \leftarrow &     \phi     & \rightarrow &
f_{\phi}                       &   {  }         &  { }
\endmatrix
$$
where the last isomorphism is the Shimura correspondence. And all the above
isomorphisms are compatible with the action of Hecke operators.
\remark{\smc Remark} \ (1) $[\G_2,k]^M = {\Bbb C}E^{(2)}_k \oplus
[\G_2,k]^M_0,$ where $E^{(2)}_k$ is the Siegel-Eisenstein series of weight $k$
on $\G_2$ given by
$$E_k^{(2)}(Z) := \sum_{\{C,D\}}\,\text{det}(CZ+D)^{-k},\quad Z\in H_2$$
(sum over non-associated pairs of coprime symmetric matrices $C,D\in\Bbb Z^{(2,2)}$)
and $[\G_2,k]^M_0 := [\G_2,k]^M \cap [\G_2,k]_0.$ \newline
(2) Maass proved that $\text{dim} [\G_2,k]^M =\left[\frac{k-4}{6}\right]$ for
$k\geq 4$ even. It is known that $\text{dim}[\G_2,k]\sim 2^{-6}\cdot 3^{-3}
\cdot 5^{-1}\cdot k^3$ as $k\rightarrow \infty.$
\endremark
\par
We observe that Theorem 2 implies that $[\G_2,k]^M$ is invariant under all the
Hecke operators and that it is annihilated by the operator
$$\Cal C_p := T_p^2 - p^{k-2}(p+1)T_p - T_{p^2}+p^{2k-2} \tag B.16$$
for every prime $p.$ We let $F\in [\G_2,k]$ be a nonzero Hecke eigenform with
$T_n F = \lambda_n F$ for $n\in \Bbb Z^+.$ For a prime $p$, we put
$$Z_{F,p}(X) := 1-\lambda_p X +(\lambda^2_p - \lambda_{p^2} - p^{2k-4}) X^2
-\lambda_p p^{2k-3} X^3 + p^{4k-6} X^4$$
so that $Z_{F,p}(p^{-s})\;(s\in \Bbb C)$ is the local spinor zeta function of
$F$ at $p.$ We put
$$Z_F(s) := \prod_p Z_{F,p}(p^{-s}),\quad \text{Re} \,s \gg 0. \tag B.17$$
Then we have
$$Z_F(s) =\zeta(2s-2k+4) \sum_{n\geq 1}\frac{\lambda_n}{n^s},\quad
\text{Re}\,s \gg 0.\tag B.18$$
\proclaim{Theorem 3 (Saito-Kurokawa conjecture ; Andrianov [An], Maass [Ma 2-3],
Zagier [Za])} \ Let $k\in \Bbb Z^+$ be even and let $F$ be a nonzero Hecke eigenform
in $[\G_2,k]^M.$ Then there exists a unique normalized Hecke eigenform $f$ in
$[\G_1,2k-2]$ such that
$$Z_F(s) = \zeta(s-k+1)\zeta(s-k+2)L_f(s),$$
where $L_f(s)$ is the Hecke $L$-function attached to $f.$
\endproclaim
\indent
Theorem 2 implies that $Z_F(s)$ has a pole at $s=k$ if $F$ is a Hecke eigenform
in $[\G_2,k]^M_0.$ If $F\in [\G_2,k]_0$ is an eigenform, it was
proved by Andrianov that
$Z_F(s)$ has an analytic continuation to $\Bbb C$ which is holomorphic
everywhere if $k$ is odd and is holomorphic except for a possible simple pole
at $s=k$ if $k$ is even. Moreover, the global function
$$Z_F^*(s) := (2\pi)^{-s} \Gamma(s) \Gamma(s-k+2)Z_F(s) \tag B.19$$
is $(-1)^k$-invariant under $s\mapsto 2k-2-s.$ It was proved by Evdokimov and
Oda that $Z_F(s)$ is holomorphic everywhere if and only if $F$ is contained in
the orthogonal complement of $[\G_2,k]^M_0$ in $[\G_2,M].$
\par
So far a generalization of the Maass space to higher genus $n>2$ has not been
given. There is a partial negative result by Ziegler (cf [Zi], Theorem 4.2).
We will describe his result roughly. Let $F\in [ \G_{g+1},k]\;(g\in \Bbb Z^+,
\;k\,:\,\text{even})$ be a Siegel modular form on $H_{g+1}$ of weight $k$ and
let
$$F\pmatrix Z_1 & {}^t W \\ W & z_2 \endpmatrix = \sum_{m\geq 0} \Phi_{F, m}
(Z_1,W)\,e^{2\pi i m z_2}, \quad \pmatrix Z_1 & {}^t W \\ W& z_2 \endpmatrix
\in H_{g+1},\;\text{with}\; Z_1\in H_g,\;z_2\in H_1$$
be the Fourier-Jacobi expansion of $F.$ For any nonnegative integer $m,$ we
consider the linear mapping
$$\rho_{g,m,k}\, :\, [\G_{g+1},k] \longrightarrow J_{k,m}(\G_g)$$
defined by
$$\rho_{g,m,k}(F) \,:\, =\, \Phi_{F,m},\;\;F\in [\G_{g+1},k].$$
Ziegler showed that for $g\geq 32,$ the mapping
$$\rho_{g,1,16}\, :\, [\G_{g+1},16] \longrightarrow J_{16,1}(\G_g)$$
is not surjective. \par
{\bf Question}: Is $\rho_{g,1,k}$ surjective for an integer $k\neq 16$\,?

\head {\bf Appendix\ C. \ The\ Orthogonal\ Group\ $O_{s+2,2}({\Bbb R})$} 
\endhead               

\indent
A lattice is a free $\BZ$-module of finite rank with a nondegerate symmetric
bilinear form with values in $\BQ$. Let $K$ be a positive definite
unimodular even integral lattice of rank $s$ with its associated symmetric
matrix $S_0.$ Let $\Pi_{1,1}$ be the unique unimodular even integral
Lorentzian lattice of rank $2$ with its associated symmetric matrix
$\pmatrix \ 0 & -1 \\ -1 &\ 0\endpmatrix.$ If there is no confusion, we
write $\Pi_{1,1}=\pmatrix \ 0 & -1 \\ -1 & 0\endpmatrix.$
\par\medpagebreak\indent
We define the lattices $L$ and $M$ by
$$L:=K\oplus \Pi_{1,1}\ \ \ \text{and}\ \ \ M:=\Pi_{1,1}\oplus L=
\Pi_{1,1}\oplus K\oplus \Pi_{1,1}.\tag C.1$$
We put
$$K_{\BR}:=K\otimes_{\BZ}\BR,\ \ \ \ L_{\BR}:=L\otimes_{\BZ}\BR,\ \ \ \
M_{\BR}:=M\otimes_{\BZ}\BC.$$
We let
$$Q_K:=S_0,\ \ \ Q_L:=\pmatrix \ 0 & 0 & -1\\ \ 0 & S_0 &\ 0\\
-1 & 0 & \ 0 \endpmatrix,\ \ \ Q_M:=\pmatrix 0 & 0 & \Pi_{1,1}\\
0 & S_0 & 0\\ \Pi_{1,1} & 0 & 0 \endpmatrix$$
be the unimodular even integral symmetric matrices associated with the
lattices $K,\,L$ and $M$ respectively. The isometry group $O_M(\BR)$ of
the quadratic space $(M_{\BR},Q_M)$ is defined by
$$O_M(\BR):=\left\{\,g\in GL(M_{\BR})\cong GL(s+4,\BR)\,\vert\
^t\!gQ_Mg=Q_M\ \right\}.\tag C.2$$
Then it is easy to see that $O_M(\BR)$ is isomorphic to the orthogonal
group $O_{s+2,2}(\BR)$. Here for two nonnegative integers $p$ and $q$ with
$p+q=n,\ O_{p,q}(\BR)$ is defined by
$$O_{p,q}(\BR):=\left\{\,g\in GL(p+q,\BR)\,\vert\ ^t\!gE_{p,q}g=E_{p,q}\
\right\},\tag C.3$$
where
$$E_{p,q}:=\pmatrix E_p & 0\\ 0 & -E_q\endpmatrix.$$
Indeed, $Q_M$ is congruent to $E_{s+2,2}$ over $\BR$, that is,
$Q_M=\,^ta\!E_{s+2,2}a$ for some $a\in GL(s+4,\BR)$ and hence
$O_M(\BR)=a^{-1}O_{s+2,2}(\BR)a.$ For brevity, we write
$O(p,\BR):=O_{p,0}(\BR)$ and $SO(p,\BR):=SO_{p,0}(\BR).$ Similiarly, we have
$O_L(\BR)\cong O_{s+1,1}(\BR)$ and $O_K(\BR)\cong O(s,\BR).$ We denote
by $(\,,\,)_K,\,(\,,\,)_L$ and $(\,,\,)_M$  the nondegenerate symmetric
bilinear forms on $K_{\BR},\,L_{\BR}$ and $M_{\BR}$ corresponding to
$Q_K,\,Q_L$ and $Q_M$ respectively. \par\medpagebreak\indent
We let
$$D=D(M_{\BR}):=\left\{\,z\subset M_{\BR}\,\vert\
\text{dim}_{\BR}\,z=2,\ z\ \text{is\ oriented\ and}\ (\,,\,)_M\vert_z <0
\right\}\tag C.4$$
be the space of oriented negative two dimensional planes in $M_{\BR}.$
We observe that a negative two dimensional plane in $M_{\BR}$ occurs twice
in $D$ with opposite orientation. Thus $D$ may be regarded as a space
consisting of two copies of the space of negative two dimensional planes
in $M_{\BR}.$ For $z\in D,$ the majorant associated to $z$ is defined by
$$(\,,\,)_z:=\cases (\,,\,)_M &\ \ \text{on}\ z^{\perp}\,;\\
-(\,,\,)_M &\ \ \text{on}\ z.\endcases\tag C.5$$
Then $(M_{\BR},\,(\,,\,)_z)$ is a positive definite quadratic space. It is
easy to see that we have the orthogonal decomposition $M_{\BR}=z^{\perp}
\oplus z$ with respect to $(\,,\,)_z$ and that $(\,,\,)_M$ has the signature
$(s+2,0)$ on $z^{\perp}$ and $(0,2)$ on $z.$\par\medpagebreak
\indent
According to Witt's theorem, $O_{s+2,2}(\BR)$ acts on $D$ transitively.
For a fixed element $z_0\in D,$ we denote by $K_{\infty}$ the stabilizer of
$O_{s+2,2}(\BR)$ at $z_0$. Then
$$D\cong O_{s+2,2}(\BR)/K_{\infty}\tag C.6$$
is realized as a homogeneous space. It is easily seen that $K_{\infty}$
is isomorphic to $O(s+2,\BR)\times SO(2,\BR),$ which is a subgroup of the
maximal compact subgroup $O(s+2,\BR)\times O(2,\BR)$ of $O_M(\BR)\cong
O_{s+2,2}(\BR)$. It is also easy to check that $O_{s+2,2}(\BR)$ has four
connected components. We denote by $SO_{s+2,2}(\Bbb R)^0$ the identity component of
$O_{s+2,2}(\BR).$ In fact, $SO_{s+2,2}(\BR)^0$ is the kernel of the spinor
norm mapping
$$\rho\,:\,SO_{s+2,2}(\BR)\longrightarrow \BR^{\times}/(\BR^{\times})^2.
\tag C.7$$
Now we know that
$$D\cong O_{s+2,2}(\BR)/O(s+2,\BR)\times SO(2,\BR)\tag C.8$$
has two connected components and the connected component $D^0$ containing
the origin $o:=z_0$ is realized as the homogeneous space as follows\,:
$$D^0\cong O_{s+2,2}(\BR)/O(s,\BR)\times O(2,\BR)\cong
SO_{s+2,2}(\BR)^0/SO(4,\BR)\times SO(2,\BR).\tag C.9$$
It is known that $D^0$ is a Hermitian symmetric space of noncompact type with
complex dimension $s+2.$ Let us describe a Hermitian structure on $D^0$
explicitly. For brevity, we write $G_{\BR}^0:=SO_{s+2,2}(\BR)^0$ and
$K_{\BR}^0:=SO(s+2,\BR)\times SO(2,\BR).$ Obviously $G_{\BR}^0$ is the
identity component of $O_{s+2,2}(\BR)\cong O_M(\BR)$ and $K_{\BR}^0$
is the identity component of $O(s+2,\BR)\times O(2,\BR)\cong K_{\infty}.$
For a positive integer $n$,
the Lie algebra ${\frak {so}}(n,\BR)$ of $SO(n,\BR)$ has dimension
$(n-1)n/2$ and
$${\frak {so}}(n,\BR)=\{\,X\in \BR^{(n,n)}\,\vert\ \sigma(X)=0,\ ^tX+X=0\,
\}.\tag C.10$$
Then the Lie algebra ${\frak g}$ of $G_{\BR}^0$ is given by
$${\frak g}=\left\{\,\pmatrix A & C\\ ^t\!C & B\endpmatrix\in
\BR^{(s+4,s+4)}\,\vert\ A\in {\frak {so}}(s+2,\BR),\ B\in {\frak {so}}
(2,\BR),\ C\in \BR^{(s+2,2)}\ \right\}.\tag C.11$$
Let $\theta$ be the Cartan involution of $G_{\BR}^0$ defined by
$$\theta(g):=E_{s+2,2}gE_{s+2,2},\ \ \ g\in G_{\BR}^0.\tag C.12$$
Then $K_{\BR}^0$ is the subgroup of $G_{\BR}^0$ consisting of elements in
$G_{\BR}^0$ fixed by $\theta.$ We also denote by $\theta$ the differential
of $\theta$ which is given by
$$\theta(X)=E_{s+2,2}XE_{s+2,2},\ \ \ X\in {\frak g}.\tag C.13$$
Then ${\frak g}$ has the Cartan decomposition
$${\frak g}={\frak k}+{\frak p},\tag C.14$$
where ${\frak k}$ and ${\frak p}$ denote the $(+1)$-eigenspace and
$(-1)$-eigenspace of $\theta$ respectively. More explicitly,
$${\frak k}=\left\{\,\pmatrix A & 0\\ 0& B\endpmatrix\in {\frak g}\,\vert
\ A\in {\frak {so}}(s+2,\BR),\ B\in {\frak {so}}(2,\BR)\ \right\}$$
and
$${\frak p}=\left\{\,\pmatrix 0 & C \\ ^t\!C & 0\endpmatrix\ \vert\
C\in \BR^{(s+2,2)}\ \right\}.$$
The real dimension of ${\frak g},\,{\frak k}$ and ${\frak p}$ are
$(s+3)(s+4)/2,\ (s^2+3s+4)/2$ and $2(s+2)$ respectively. Thus the real
dimension of $D^0$ is $2(s+2).$ Since ${\frak p}$ is stable under the
adjoint action of $K_{\BR}^0$, i.e., $\text{Ad}(k){\frak p}={\frak p}$
for all
$k\in K_{\BR}^0,\ ({\frak g},{\frak k},\theta)$ is reductive. Thus the
tangent space $T_o(D^0)$ of $D^0\subset D$ at $o:=z_0$ can be canonically
identified with $\BR^{2(s+2)}$ via
$$\pmatrix 0 & C\\ ^tC & 0 \endpmatrix \longmapsto (x,y),\ \ \
C=(x,y)\in \BR^{(s+2,2)}\cong \BR^{2(s+2)}.$$
Then the adjoint action of $K_{\BR}^0$ on ${\frak p}\cong \BR^{2(s+2)}$ is
expressed as
$$\text{Ad}\left(\pmatrix k_1 & 0 \\ 0 & k_2\endpmatrix\right)(x,y)=
\pmatrix 0 & k_1C\,^t\!k_2\\ k_2\,^tC\,^t\!k_1 & 0 \endpmatrix,\tag C.15$$
where $k_1\in SO(s+2,\BR),\ k_2\in SO(2,\BR)$ and $C=(x,y)\in
\BR^{(s+2,2)}.$ The Cartan-Killing form $B$ of ${\frak g}$ is given by
$$B(X,Y)=(s+2)\,\sigma(XY),\ \ \ X,Y\in {\frak g}.\tag C.16$$
The restriction $B_0$ of $B$ to ${\frak p}$ is given by
$$B_0((x,y),(x',y'))=2(s+2)(<x,x'>+<y,y'>),$$
where $<\,,\,>$ is the standard inner product on $\BR^{s+2}.$ The restriction
$B_0$ induces a $G_{\BR}^0$-invariant Riemannian metric $g_0$ on $D^0$
defined by
$$g_0(X,Y):=B_0(X,Y),\ \ \ X,Y\in {\frak p}.$$
It is easy to check that $g_0$ is invariant under the adjoint action of
$K_{\BR}.$\par\medpagebreak\indent
Now let $J_0$ be the complex structure on the real vector space
${\frak p}$ defined by
$$J_0((x,y)):=(-y,x),\ \ \ (x,y)\in {\frak p}.\tag C.17$$
We note that
$$J_0=\text{Ad}\pmatrix E_{s+2} & 0 \\ 0 & I_\ast \endpmatrix,\ \ \
I_\ast:=\pmatrix 0 & -1 \\ 1 & \ 0\endpmatrix,\ \ \ J_0^2=-\text{Id}_{\frak p}.$$
It is easy to check that $J_0$ is Ad$(K_{\BR})$-invariant, i.e.,
$$J_0(\text{Ad}(k)X)=J_0(X)\ \ \text{for\ all}\ k\in K_{\BR}\ \text{and}\
X\in {\frak p}.$$
Hence $J_0$ induces an almost complex structure $J$ on $D^0$ and also on
$D.\ J$ becomes a complex structure on $D^0$ via the natural identification
$$T_0(D^0)\cong {\frak p}\cong \BR^{2(s+2)}\cong \BC^{s+2},\ \ \
(x,y)\longmapsto x+iy,\ \ x,y\in \BR^{s+2}.\tag C.18$$
Indeed, $J$ is the pull-back of the standard complex structure on
$\BC^{s+2}$. The complexification ${\frak p}_{\BC}:={\frak p}\otimes_{\BR}
\BC$ has a canonical decomposition
$${\frak p}_{\BC}={\frak p}_+\oplus {\frak p}_- ,$$
where ${\frak p}_+\,(\,\text{resp.}\,{\frak p}_-\,)$ denotes the
$(+i)$-eigenspace\,(\,resp.\,$(-i)$-eigenspace\,) of $J_0.$ Precisely,
${\frak p}_+$ and ${\frak p}_-$ are given by
$${\frak p}_+=\{\,(x,-ix)\,\vert\ x\in \BC^{s+2}\,\}\ \ \ \text{and}\ \ \
{\frak p}_-=\{\,(x,ix)\,\vert\ x\in \BC^{s+2}\,\}.$$
Usually ${\frak p}_+$ and ${\frak p}_-$ are called the {\it holomorphic}
tangent space and the {\it anti}-{\it holomorphic} tangent space
respectively. Moreover, it is easy to check that the Riemannian metric
$g_0$ on $D^0$ is Hermitian with respect to the complex structure $J$, i.e.,
$g_0(JX,JY)=g_0(X,Y)$ for all smooth vector fields $X$ and $Y$ on $D^0.$ And
$D^0$ has the canonical orientation induced by its complex structure.
\par\medpagebreak\indent
In summary, we have\par\medpagebreak\noindent
{\bf Theorem\ 1.} $D^0$ is a Hermitian symmetric space of noncompact type
with dimension $s+2.\ D^0$ is realized as a bounded symmetric domain in
$\BC^{s+2}$ and hence $D$ is a union of two bounded symmetric domains in
$\BC^{s+2}.$\par\medpagebreak\noindent
{\smc Remark\ 2.} We choose an orthogonal basis of $z_0^{\perp}.$ We also
choose a basis of $z_0$ which is properly oriented. Let
$$\tau_{z_0^{\perp}}:=\text{diag}\,(E_{s+1},-1)\ \ \ \text{and}\ \ \
\tau_{z_0}:=\text{diag}\,(1,-1)$$
be the symmetries in the isometry groups $O(z_0^{\perp})$ and $O(z_0)$
with respect to the last coordinates of $z_0^{\perp}$ and $z_0$ respectively.
We observe that $\tau_{z_0}$ reverses the orientation of $z_0$ and lies in
$O(z_0)-SO(z_0).$ It is easy to check that
$$\rho(1_{z_0^{\perp}}\times \tau_{z_0})=-1,\ \ \ \rho(\tau_{z_0^{\perp}}
\times 1_{z_0})=1,$$
where $\rho$ is the spinor norm mapping defined by (C.7). It is easy to see
that the set
$$\left\{\,1_{M_{\BR}},\ 1_{z_0^{\perp}}\times \tau_{z_0},\
\tau_{z_0^{\perp}}\times 1_{z_0},\ \tau_{z_0^{\perp}}\times \tau_{z_0}\
\right\}\tag C.19$$
is a complete set of coset representatives of $O_{s+2,2}(\BR)/
SO_{s+2,2}(\BR)^0.$ We note that the set (C.19) is contained in
$O(s+2,\BR)\times O(2,\BR)$ and so that (C.19) is a complete set of coset
representatives of $O(s+2,\BR)\times O(2,\BR)/(SO(s+2,\BR)\times SO(2,\BR))$.
It is easy to see that the set $\{\,1_{M_{\BR}},\ 1_{z_0^{\perp}}\times
\tau_{z_0}\,\}$ is a complete set of coset representatives of
$(O(s+2,\BR)\times O(2,\BR))/(O(s+2,\BR)\times SO(2,\BR)).$ Thus we have
$$D=D^0\cup (1_{z_0^{\perp}}\times \tau_{z_0})D^0.\tag C.20$$
The complex structure $-J_0$ on ${\frak p}$ determines the
opposite almost complex structure on $D^0$ and the almost complex structure
on the connected component $D-D^0$ is the one on $D^0$ carried by the
element $1_{z_0^{\perp}}\times \tau_{z_0}.$ The ground manifolds $D^0$ and
$D-D^0$ may be regarded as the same one, but each carries the opposite
almost complex structure.
\par\medpagebreak\noindent
\remark
{\smc Remark\ 3} $D^0$ may be regarded as an open orbit of $G_{\BR}^0$
in the complex projective quadratic space $^tzQ_Mz=0$ via the Borel
embedding.\,(\,See [Bai] for detail.) $D^0$ is realized as a tube domain
in $\BC^{s+2}$ given by (4.10) in section 4. For the explicit realization
of $D^0$ as a bounded symmetric domain in $\BC^{s+2},$ we refer to [Bai],\,
[H] and [O].\endremark
\par\medpagebreak\indent
Finally we present the useful equations for $g$ to belong to
$O_{s+2,2}(\BR).$ For $g\in O_{s+2,2}(\BR),$ we write
$$g=\pmatrix A_{11} & A_{12} & A_{13} \\ A_{21} & A_{22} & A_{23} \\
A_{31} & A_{32} & A_{33} \endpmatrix,$$
where $A_{11},\,A_{13},\,A_{31},\,A_{33}\in \BR^{(2,2)},\ A_{22}\in
\BR^{(s,s)},\ A_{12},\,A_{32}\in \BR^{(2,s)}$, and $A_{21},\,A_{23}\in
\BR^{(s,2)}.$ Then the condition $^tgQ_Mg=Q_M$ is equivalent to the
following equations given by
$$^tA_{11}\Pi_{1,1}A_{31}+\,^tA_{21}S_0A_{32}+\,^tA_{31}\Pi_{1,1}
A_{12}=0,\tag C.21$$
$$^tA_{11}\Pi_{1,1}A_{32}+\,^tA_{21}S_0A_{22}+\,^tA_{31}\Pi_{1,1}
A_{12}=0,\tag C.22$$
$$^tA_{11}\Pi_{1,1}A_{33}+\,^tA_{21}S_0A_{23}+\,^tA_{31}\Pi_{1,1}
A_{13}=\Pi_{1,1},\tag C.23$$
$$^tA_{12}\Pi_{1,1}A_{32}+\,^tA_{22}S_0A_{22}+\,^tA_{32}\Pi_{1,1}
A_{12}=S_0,\tag C.24$$
$$^tA_{12}\Pi_{1,1}A_{33}+\,^tA_{22}S_0A_{23}+\,^tA_{32}\Pi_{1,1}
A_{13}=0,\tag C.25$$
and
$$^tA_{13}\Pi_{1,1}A_{33}+\,^tA_{23}S_0A_{23}+\,^tA_{33}\Pi_{1,1}
A_{13}=0.\tag C.26$$
\par\bigpagebreak

\head {\bf Appendix\ D. \ The\ Leech\ Lattice\ $\Lambda$} 
\endhead               

\redefine\L{\Lambda}
\indent
Here we collect some properties of the Leech lattice $\L$. Most of the
materials in this appendix can be found in [C-S].\par\medpagebreak
\indent
The Leech lattice $L$ is the unique positive definite unimodular even
integral lattice of rank $24$ with minimal norm $4.\ \L$ was discovered
by J. Leech in $1965$. (\,cf.\,Notes on sphere packing, Can. J. Math. 19
(1967), 251-267.\,) It was realized by Conway, Parker and Sloane
that the Leech lattice $\L$ has
many strange geometric properties.
Past three decades more than 20 constructions of $\L$ were found.
\par\medpagebreak\indent
The following properties of $\L$ are well known\,:\par\medpagebreak
\noindent
$(\L1)$ The determinant of $\L$ is $\text{det}\,\L=1.$ The kissing number
is $\tau=196560$ and the packing radius is $\rho=1.$ The density is
$\Delta=\pi^{12}/(12!)=0.001930\cdots.$ The covering radius is
$R=\sqrt{2}$ and the thickness is $\Theta=(2\pi)^{12}/(12!)=7.9035\cdots.$
\par\smallpagebreak\noindent
$(\L2)$ There are $23$ different types of deep hole one of which is the
octahedral hole $8^{-1/2}(4,0^{23})$ surrounded by $48$ lattice points.
\par\smallpagebreak\noindent
$(\L3)$ The Veronoi cell has $16969680$ faces, $196560$ corresponding to
the minimal vectors and $16773120$ to those of the next layer.
\par\smallpagebreak\noindent
$(\L4)$ The automorphism group $\text{Aut}(\L)$ of $\L$ has order
$$2^{22}\cdot 3^9\cdot 5^4\cdot 7^2\cdot 11\cdot 13\cdot 23=
8315553613086720000.$$
$\text{Aut}(\L)$ has the Mathieu group $M_{24}$ as a subgroup. The
automorphism group $\text{Aut}(\L)$ is often denoted by $Co_0$ or
$\cdot\!0$ because J.H. Conway first discovered this group.
\par\bigpagebreak\indent
For a given lattice $L$, we denote $N_m(L)$ by the number of vectors of
norm $m$. Conway characterized the Leech lattice as follows\,(\,cf.\,
A characterization of Leech's lattice, Invent. Math. 7 (1969), 137-142 or
Chapter 12 in [C-S]\,)\,:\par\medpagebreak\noindent
\proclaim
{\bf Theorem\ 1\,(\,Conway\,)} $\L$ is the unique positive definite
unimodular even integral lattice $L$ with rank $<32$ that satisfies any one of
the following\par\medpagebreak\noindent
(a) $L$ is not directly congruent to its mirror-image.\par\smallpagebreak
\noindent
(b) No reflection leaves $L$ invariant.\par\smallpagebreak\noindent
(c) $N_2(L)=0.$\par\smallpagebreak\noindent
(d) $N_{2m}(L)=0$ for some $m\geq 0.$
\endproclaim
\par\bigpagebreak\noindent
\proclaim
{\bf Theorem\ 2\,(\,Conway\,)} If $L$ is a unimodular even integral lattice
with rank $<32$ and $N_2(L)=0,$ then $L=\L$ and $N_4(L)=196560,\
N_6(L)=16773120,\ N_8(L)=398034000.$
\endproclaim
\par\medpagebreak\indent
Now we review the Jacobi theta functions.
For the present time being, we put $q=e^{i\pi\tau}$ and $\zeta=e^{i\pi z}.$
(\,We note that we set $q=e^{2\pi i\tau}$ and $\zeta=e^{2\pi iz}$ at other
places.) We define the Jacobi theta functions
$$
\theta_1(\tau,z):=i^{-1}\sum_{n\in \BZ}(-1)^n
q^{(n+1/2)^2}\zeta^{2n+1},\tag D.1$$
$$\theta_2(\tau,z):=\sum_{n\in \BZ}q^{(n+1/2)^2}\zeta^{2n+1},
\tag D.2$$
$$\theta_3(\tau,z):=\sum_{n\in \BZ}q^{n^2}\zeta^{2n},\tag D.3$$
$$\theta_4(\tau,z):=\sum_{n\in \BZ}(-1)^nq^{n^2}\zeta^{2n},\tag D.4$$
where $\tau\in H_1$ and $z\in \BC.$ We also define the theta functions
$\theta_k(\tau):=\theta_k(\tau,0)$ for $k=2,3,4.$ Then it is easy to see
that
$$\theta_2(\tau)=\sum_{n\in\BZ}q^{(n+1/2)^2}=e^{\pi i\tau/4}
\theta_3(\tau,\tau/2),\tag D.5$$
$$\theta_4(\tau)=\sum_{n\in \BZ}(-1)^nq^{n^2}=\theta_3(\tau,\pi/2)=
\theta_3(\tau+1).\tag D.6$$
\indent
According to the Poisson summation formula, we obtain
$$\theta_3(-1/\tau,z/\tau)=(-i\tau)^{1/2}e^{\pi iz^2/\tau}\theta_3(\tau,z),
\ \ \ (\tau,z)\in H_1\times \BC.\tag D.7$$
And these theta functions can be written as infinite products as follows\,:
$$\theta_1(\tau,z)=2\,\text{sin}\,\pi z\cdot q^{1/4}\prod_{n>0}
(1-q^{2n})(1-q^{2n}\zeta^2)(1-q^{2n}\zeta^{-2}),\tag D.8$$
$$\theta_2(\tau,z)=q^{1/4}\zeta\prod_{n>0}(1-q^{2n})(1+q^{2n}\zeta^2)
(1+q^{2n}\zeta^{-2}),\tag D.9$$
$$\theta_3(\tau,z)=\prod_{n>0}(1-q^{2n})(1+q^{2n}\zeta^2)(1+q^{2n-1}
\zeta^{-2}),\tag D.10$$
$$\theta_4(\tau,z)=\prod_{n>0}(1-q^{2n})(1-q^{2n-1}\zeta^2)
(1-q^{2n-1}\zeta^{-2}).\tag D.11$$
Thus the theta functions $\theta_k(\tau)\,(\,2\leq k\leq 4\,)$ are written
as infinite products\,:
$$\theta_2(\tau)=q^{1/4}\prod_{n>0}(1-q^{2n})(1+q^{2n})(1+q^{2n-2})
=2q^{1/4}\prod_{n>0}(1-q^{2n})(1+q^{2n})^2,\tag D.9'$$
$$\theta_3(\tau)=\prod_{n>0}(1-q^{2n})(1+q^{2n-1})^2,\tag D.10'$$
$$\theta_4(\tau)=\prod_{n>0}(1-q^{2n})(1-q^{2n-1})^2.\tag D.11'$$
We note that the discriminant function $\Delta(\tau)$ is written as
$$\Delta(\tau)=q^2\prod_{n>0}(1-q^{2n})^{24}=\left\{\,1/2\theta_2(\tau)
\theta_3(\tau)\theta_4(\tau)\,\right\}^8.\tag D.12$$
We observe that the theta function $\theta_3(\tau,z)$ is annihilated by
the heat operator $H:={{\partial^2}\over {\partial z^2}}-4\pi i
{{\partial}\over {\partial \tau}}.$ It is easy to check that
$\theta_1(\tau,z)$ has zeros only at $m_1+m_2\tau\,(\,m_1,m_2\in \BZ\,)$
and satisfies the equations
$$\theta_1(\tau,z+1)=-\theta_1(\tau,z),\ \ \ \theta_1(\tau,\tau+z)=
-q^{-1}e^{-2\pi iz}\theta_1(\tau,z).\tag D.13$$
\def\a{\alpha}
Now for a given positive definite lattice $L$, we define the {\it theta\
series} $\Theta_L(\tau)$ of a lattice $L$ by
$$\Theta_{L}(\tau):=\sum_{\a\in L}q^{N(\a)}=
\sum_{m\geq 0}N_m(L)q^m,\ \ \tau\in H_1,\tag D.14$$
where $N(\a):=(\a,\a)$ denotes the norm of a vector $\a\in L.$ We can also
use (D.14) to define the theta series of a nonlattice packing $L$. The
commonest examples of this appear when $L$ is a translate of a lattice or
a union of translates. Clearly $\theta_2(\tau)=\Theta_{\BZ+1/2}(\tau),\
\theta_3(\tau)=\Theta_{\BZ}(\tau)$ and $\Theta_{\BZ^n}(\tau)=
\Theta_{\BZ}(\tau)^n=\theta_3(\tau)^n.$\par\medpagebreak\indent
Returing to the Leech lattice $L$,
$$\align\Theta_{\L}(\tau)&=\Theta_{E_8}(\tau)^3-720\,\Delta(\tau)\\
&=1/8\left\{\,\theta_2(\tau)^8+\theta_3(\tau)^8+\theta_4(\tau)^8\,\right\}^3
-45/16\left\{\,\theta_2(\tau)\theta_3(\tau)\theta_4(\tau)\,\right\}^8\\
&=1/2\left\{\,\theta_2(\tau)^{24}+\theta_3(\tau)^{24}+\theta_4(\tau)^{24}\,
\right\}-69/16\left\{\,\theta_2(\tau)\theta_3(\tau)\theta_4(\tau)\,
\right\}^8\\
&=\sum_{m\geq 0}N_m(\L)q^m=1+196560q^4+16773120q^6+\cdots,\endalign
$$
where $\Theta_{E_8}(\tau)$ is the theta series of the exceptional lattice
$E_8$ of rank $8.$ It is known that
$$N_m(\L)=65520/691\,\left(\,\sigma_{11}(m/2)-\tau(m/2)\right).\tag D.15$$
The values of $N_m(\L)$ for $0\leq m\leq 100,\,m:\text{even}$ can be found in
[C-S],\,p.\,135.\par\medpagebreak\indent
In the middle of 1980s, M. Koike, T. Kondo and T. Tasaka solved a special
part of the Moonshine Conjectures for the Mathieu group $M_{24}.$ For
$g\in M_{24},$ we write
$$g=(n_1)(n_2)\cdots (n_s),\ \ \ n_1\geq\cdots \geq n_s\geq 1,\tag D.16$$
where $(n_i)$ is a cycle of length $n_i\,(\,1\leq i\leq s\,).$ To each
$g\in M_{24}$ of the form (D.16), we associate modular forms
$\eta_g(\tau)$ and $\theta_g(\tau)$ defined by
$$\eta_g(\tau):=\eta(n_1\tau)\,\eta(n_2\tau)\cdots\eta(n_s\tau),\ \ \
\tau\in H_1\tag D.17$$
and
$$\theta_g(\tau):=\sum_{\a\in \L_g}e^{\pi iN(\a)\tau},\ \ \ \tau\in H_1,
\tag D.18$$
where $\eta(\tau)$ is the Dedekind eta function and
$$\L_g:=\{\,\a\in \L\,\vert\ g\cdot \a=\a\,\}\tag D.19$$
is the positive definite even integral lattice of rank $s$. We observe
that $\theta_g(\tau)=\Theta_{\L_g}(\tau).$\par\medpagebreak\noindent
\proclaim
{\bf Theorem\ 3\,(\,[Koi2]\,)} For any element $g\in M_{24}$
with $g\ne 12^2,\ 4^6,\ 2^{12},\ 10^2\cdot 2^2,\ 12\cdot 6\cdot 4\cdot 2,
\ 4^4\cdot 2^4,$ there exists a unique modular form $f_g(\tau)=1+
\sum_{n\geq 0}a_g(n)q^{2n},\ a_g(n)\in \BZ$ satisfying the following
conditions\,:\par\medpagebreak\indent
(K1) There exists an element $g_1\in G$ such that $f_g(\tau)\eta_g(\tau)^{-1}
=T_{g_1}(\tau)+c$ for some constant $c$, where $G$ is the MONSTER and
$T_{g_1}(\tau)$ denotes the Thompson series of $g_1\in G.$
\par\smallpagebreak\indent
(K2) $a_g(1)=0,$ and $a_g(n)$ are nonnegative even integers for all
$n\geq 1.$\par\smallpagebreak\indent
(K3) If $g'=g^r$ for some $r\in \BZ$, then $a_g(n)\leq a_{g'}(n)$ for all
$n$.\par\smallpagebreak\indent
(K4) $a_g(2)$ is equal to the cardinality of the set
$\{\,\a\in \L_g\,\vert\ N(\a)=(\a,\a)=4\,\}.$
\endproclaim
\par\medpagebreak\noindent
\proclaim
{\bf Theorem\ 4\,(\,Kondo\ and\ Tasaka\,[K-T]\,)} Let $g\in M_{24}$ be any
element of the Mathieu group $M_{24}.$ Then the function $\theta_g(\tau)
\,\eta_g(\tau)^{-1}$ is a Hauptmodul for a genus $0$ discrete subgroup of
$SL(2,\BR).$ The function $\theta_g(\tau)$ is the unique modular form
satisfying the conditions (K1)-(K4).
\endproclaim
\par\bigpagebreak
\indent
For more detail on the Leech lattice $\L$ we refer to [Bo3],\,[C-S] and
[Kon].

\vskip 1cm

{\smc Department\ of\ Mathematics\par Inha\ University\par Incheon
402-751
\par
Republic of Korea} \vskip 0.3cm {\tt e-mail\ address\,:\,\
jhyang\@inha.ac.kr}

\end